\documentclass[reqno]{amsart}

\usepackage {fullpage}						% set all margins to 1.5cm
\usepackage{csquotes}
\usepackage [utf8]{inputenc}
\usepackage {graphicx}
\usepackage [bookmarks=false, pdfstartview={FitH}, colorlinks=true, linkcolor=blue, citecolor=blue, urlcolor=blue]{hyperref} % create the .pdf links
\usepackage {amsfonts,amsmath,amsthm,amssymb}
\usepackage {latexsym}
\usepackage {graphics, graphicx}
\usepackage {fancyhdr}
\usepackage[table]{xcolor}
\usepackage {epstopdf} 						% to automatically transform .eps figs to .pdf
\usepackage {multicol}						% to locally have two columns
\usepackage {multirow}  					% merge rows in a tabular
\usepackage {soul} 							% to break \underline to several lines using \ul
\usepackage {datetime} 						% date and time in latex
\usepackage {subfigure}						% to have enumerated subfigures in one figure
\usepackage {parskip}						% no indentation and space between paragraphs
\usepackage {setspace}						% to set the vertical spacing
\usepackage[font=small,textfont=it,width=0.97\textwidth]{caption}		% to set the fonts for the captions, to use caption*
\usepackage{mathtools}						% to have dcases environment
\usepackage{tikz,pgfplots}
\usepackage{nicefrac}
\usepackage{bm}
\usepackage{cite}
\usepackage{morefloats}
\usepackage[section]{placeins}
\usepackage{siunitx}

\usetikzlibrary{external}
\usetikzlibrary{matrix}
\newlength\figureheight
\newlength\figurewidth

\theoremstyle{remark}
\newtheorem{remark}{Remark}[section]
\newtheorem{example}{Example}[section]
\newenvironment{acknowledgment}{{\flushleft \bf Acknowledgment:}}{}

\numberwithin{equation}{section}

\allowdisplaybreaks[1]

\renewcommand{\(}{\left(}
\renewcommand{\)}{\right)}

\newcommand{\w}{\bm w}
\newcommand{\x}{\bm x}
\newcommand{\e}{\bm e}
\newcommand{\F}{\bm F}

\newcommand{\D}{\mathfrak D}
\newcommand{\R}{\mathfrak R}

\def\softd{{\leavevmode\setbox1=\hbox{d}%
		\hbox to 1.05\wd1{d\kern-0.4ex{\char039}\hss}}}%cstocs

\newcommand{\rhoprime}{\rho^\prime}
\newcommand{\pprime}{p^\prime}
\newcommand{\uprime}{\bm{u}^\prime}
\newcommand{\thetaprime}{\theta^\prime}
\newcommand{\rhothetaprime}{(\rho\theta)^\prime}
\newcommand{\rhobar}{\bar{\rho}}
\newcommand{\pbar}{\bar{p}}
\newcommand{\ubar}{\bar{\bm{u}}}
\newcommand{\thetabar}{\bar{\theta}}
\newcommand{\rhothetabar}{\overline{\rho\theta}}
\newcommand{\uvec}{\bm{u}}

\definecolor{darkgreen}{rgb}{0,0.5,0.1}

\graphicspath{{figures/}}

\title{Stochastic Galerkin method for cloud simulation. Part II: a fully random Navier-Stokes-cloud model}

\author{A. Chertock}
\address[A. Chertock]{\newline Department of Mathematics, North Carolina State University}
\email{chertock@math.ncsu.edu}

\author{A. Kurganov}
\address[A. Kurganov]{\newline Department of Mathematics, SUSTech International Center for Mathematics and Guangdong Provincial Key
Laboratory of Computational Science and Material Design, Southern University of Science and Technology, Shenzhen, 518055, China}
\email{alexander@sustech.edu.cn}

\author{M. Luk\'a\v{c}ov\'a-Medvi\v{d}ov\'a}
\address[M. Luk\'a\v{c}ov\'a-Medvi{\softd}ov\'a]{\newline Institute of Mathematics, Johannes Gutenberg-University Mainz}
\email{lukacova@uni-mainz.de}

\author{P. Spichtinger}
\address[P. Spichtinger]{\newline Institute of Atmospheric Physics, Johannes Gutenberg-University Mainz}
\email{spichtin@uni-mainz.de}

\author{B. Wiebe}
\address[B. Wiebe]{\newline Institute of Mathematics, University of Zurich}
\email{bettina.wiebe@math.uzh.ch}

\begin{document}
\maketitle
	
\vspace{-1cm}
	
\begin{abstract}
This paper is a continuation of the work presented in [Chertock et al., Math. Cli. Weather Forecast. 5, 1 (2019), 65--106]. We study
uncertainty propagation in warm cloud dynamics of weakly compressible fluids. The mathematical model is governed by a multiscale system of
PDEs in which the macroscopic fluid dynamics is described by a weakly compressible Navier-Stokes system and the microscopic cloud dynamics
is modeled by a convection-diffusion-reaction system. In order to quantify uncertainties present in the system, we derive and implement a
generalized polynomial chaos stochastic Galerkin method. Unlike the first part of this work, where we restricted our consideration to the
partially stochastic case in which the uncertainties were only present in the cloud physics equations, we now study a fully random
Navier-Stokes-cloud model in which we include randomness in the macroscopic fluid dynamics as well. We conduct a series of numerical
experiments illustrating the accuracy and efficiency of the developed approach.
\end{abstract}
	
\medskip
\section{Introduction}
In this paper we continue the study of generalized polynomial chaos (gPC) stochastic Galerkin methods for multiscale cloud dynamics with
uncertainty. The gPC expansion based methods, such as stochastic Galerkin \cite{SG2019,XK02,PDL,TLNE10a,WK06,DPL13} and stochastic
collocation \cite{MZ09,WLB09,XH05}, are well-known techniques for uncertainty quantification in physical and engineering applications. In
these methods the solution is approximated by a truncated generalized Fourier series in terms of orthogonal polynomials corresponding to the
underlying probability density function. The stochastic collocation method belongs to the class of {\em non-intrusive methods} in which the
original deterministic system is solved at certain collocation points in the stochastic space and then the gPC coefficients are obtained by
using a polynomial interpolation and a suitable numerical quadrature. The stochastic Galerkin method on the other hand is an {\em intrusive
method} in which the Galerkin projection is employed after substituting the gPC expansions into the original system. This leads to a larger
coupled system of deterministic PDEs for the gPC coefficients for which an accurate and efficient numerical method is to be developed.

The multiscale cloud model we study is governed by a coupled macro-microscopic Navier-Stokes-cloud system  with random data and parameters.
In the first part of this work \cite{SG2019}, we restricted our consideration to the partially stochastic case in which the uncertainties
were only present in the cloud physics equations keeping the fluid dynamics deterministic. For that model we have developed an operator
splitting approach in which the deterministic macroscopic Navier-Stokes equations were solved numerically by an implicit-explicit (IMEX)
asymptotically preserving (AP) finite volume method, while the stochastic microscopic cloud dynamics was solved by a gPC expansion. The
latter was realized by the stochastic Galerkin method combined with an explicit finite volume approximation for the system of gPC
coefficients. The coupling of the deterministic Navier-Stokes equations with the stochastic cloud dynamics was realized by using expected
values for the coupling source terms.

In the present work we extend the gPC stochastic Galerkin (gPC-SG) method from \cite{SG2019} to a fully coupled random multiscale
Navier-Stokes-cloud model where we include randomness in the macroscopic fluid dynamics as well. Consequently, the random coupling source terms will be considered. This will result in deterministic PDEs for the gPC coefficients representing both Navier-Stokes and cloud dynamics. The gPC coefficients for the fluid dynamics will be approximated by using an IMEX AP finite volume method similar to the one used for the deterministic Navier-Stokes system. The gPC coefficients for the cloud dynamics will be computed by an adaptation of the explicit finite volume developed in \cite{SG2019}.

The present paper is organized as follows. For the consistency of presentation, we start in \S\ref{mathematical_model} with a short recap of
the deterministic multiscale Navier-Stokes-cloud model and then introduce a fully random cloud model. The numerical method for the latter is
presented in \S\ref{numerical_scheme}, which is followed by numerical experiments in \S\ref{numerical_experiments}. We demonstrate the
experimental convergence of the developed method as well as its applicability for uncertainty quantification in atmospheric flows through
well-known meteorological benchmark tests of a rising warm and moist bubble and Rayleigh-B\'{e}nard convection.

\section{Mathematical model}\label{mathematical_model}
We consider a mathematical model of cloud dynamics, which is based on the weakly compressible nonhydrostatic Navier-Stokes equations for
moist atmosphere (that is, a mixture of ideal gases such as dry air and water vapor),
\begin{equation}
\begin{aligned}
\rho_t+\nabla\cdot\(\rho\bm u\)&=0,\\	
(\rho\bm u)_t+\nabla\cdot\left(\rho\bm u\otimes\bm u+p\,{\rm Id}-\mu_m\rho\left(\nabla\bm u+(\nabla\bm u)^\top\right)\right)&=
-\rho g\bm{e_d},\\
(\rho\theta)_t+\nabla\cdot\({\rho\theta}\bm u-\mu_h\rho\nabla\theta\)&=S_\theta,
\end{aligned}
\label{NS_equations}
\end{equation}
and evolution equations for cloud variables,
\begin{equation}
\begin{aligned}
(\rho q_v)_t+\nabla\cdot\(\rho q_v\bm u-\mu_q\rho\nabla q_v\)&=\rho(-C+E),\\
(\rho q_c)_t+\nabla\cdot\(\rho q_c\bm u-\mu_q\rho\nabla q_c\)&=\rho(C-A_1-A_2),\\
(\rho q_r)_t+\nabla\cdot\(-v_q\rho q_r\bm{e_d}+\rho q_r\bm u-\mu_q\rho\nabla q_r\)&=\rho(A_1+A_2-E).
\end{aligned}
\label{cloud_equations}
\end{equation}
Here, $t$ is the time variable, $\x=(x_1,\dots,x_d)\in\mathbb{R}^d$ is the space vector, $\rho$ is the density, $\bm u=(u_1,\dots,u_d)^\top$
is the velocity vector, $\theta$ is the moist potential temperature, $p$ is the pressure, $g$ is the acceleration due to gravity, $\mu_m$ is
the dynamic viscosity, $\mu_h$ is the thermal conductivity, and $\mu_q$ is the cloud diffusivity. Furthermore,
$\bm{e_d}=\bm{e_3}=(0,0,1)^\top$ and $\bm{e_d}=\bm{e_2}=(0,1)^\top$ in the three-dimensional (3-D) and two-dimensional (2-D) cases,
respectively. The cloud variables representing the mass concentration of water vapor, cloud droplets and rain drops, $q_v$, $q_c$ and $q_r$,
respectively, are given by
\begin{equation*}
q_\ell=\frac{\mbox{mass of the respective phase}}{\mbox{mass of dry air}}\quad\mbox{for}\quad\ell\in\{v,c,r\}.
\end{equation*}
The terms $C$ and $E$ represent phase changes between vapor and cloud water (droplets), and $A_1$ and $A_2$ represent collision processes,
which lead to the formation of large droplets and thus precipitation.

Note that the systems \eqref{NS_equations} and \eqref{cloud_equations} are coupled through the source term $S_\theta$, which represents the
impact of phase changes and is given by
\begin{equation*}
S_\theta=\rho\frac{L\theta}{c_p T}(C-E).
\end{equation*}
For a detailed description of $S_\theta$ and the terms $E$, $C$, $A_1$ and $A_2$ see \cite{SG2019}. The temperature $T$ can be obtained from
the moist adiabatic ideal gas equation
\begin{equation*}
T=\frac{R}{R_m}\theta\(\frac{p}{p_0}\)^{\nicefrac{R_m}{c_p}},
\end{equation*}
where $p_0=10^5\,\mathrm{Pa}$ is the reference pressure at the sea level. In addition to the usual definition of a potential temperature, we
use $R_m=(1-q_v-q_c-q_r)R+q_v R_v$ with the ideal gas constant of dry air $R=287.05\,\mathrm{\nicefrac{J}{(kg\cdot K)}}$, the gas constant
of water vapor $R_v=461.51\,\mathrm{\nicefrac{J}{(kg\cdot K)}}$ and the specific heat capacity of dry air for constant pressure
$c_p=1005\,\mathrm{\nicefrac{J}{(kg \cdot K)}}$. In order to close the system, we determine the pressure from the equation of state that
includes moisture
\begin{equation}
p=p_0\left(\frac{R\rho\theta}{p_0}\right)^{\gamma_m}\quad\mbox{with}\quad\gamma_m=\frac{c_p}{c_p-R_m}.
\label{state_equation}
\end{equation}
We note that in the dry case, that is when $q_v=q_c=q_r=0$, $R_m$ reduces to $R$, $S_\theta=0$ and the moist ideal gas equation as well as
the moist equation of state become their dry analogue.

Solving the Navier-Stokes equations \eqref{NS_equations} in a weakly compressible regime is known to trigger numerical instabilities due to
the multiscale effects. We follow the approach typically used in meteorological models, where the dynamics of interest is described by a
perturbation of a background state, which is the hydrostatic equilibrium. The latter expresses a balance between the gravity and pressure
forces. Denoting by $\pbar$, $\rhobar$, $\ubar=\bm{0}$, $\thetabar$ and $\rhothetabar$ the respective background state, the hydrostatic
equilibrium satisfies
\begin{equation*}
\frac{\partial\pbar}{\partial x_d}=-\rhobar g,\quad S_\theta=0,
\end{equation*}
where $\pbar$ is obtained from the equation of state \eqref{state_equation}
\begin{equation}
\pbar=p(\rhothetabar)=p_0\left(\frac{R\rhothetabar}{p_0}\right)^{\gamma_m}.
\label{p_equ}
\end{equation}
Let $\pprime$, $\rhoprime$, $\uprime$, $\thetaprime$ and $\rhothetaprime$ stand for the corresponding perturbations of the equilibrium
state, then
\begin{equation*}
p=\pbar+\pprime,~~\rho=\rhobar+\rhoprime,~~\theta=\thetabar+\thetaprime,~~\uvec=\uprime,~~
\rho\theta=\rhobar\thetabar+\rhobar\thetaprime+\rhoprime\thetabar+\rhoprime\thetaprime=\rhothetabar+\rhothetaprime.
\end{equation*}
The pressure perturbation $\pprime$ is derived from \eqref{state_equation} and \eqref{p_equ} using the following Taylor expansion
\begin{equation*}
p(\rho\theta)\approx p(\rhothetabar)+\frac{\partial p}{\partial(\rho\theta)}\(\rho\theta-\rhothetabar\)=
\pbar+\gamma_mp_0\(\frac{R\rhothetabar}{p_0}\)^{\gamma_m}\frac{\rhothetaprime}{\rhothetabar},
\end{equation*}	
which results in
\begin{equation*}
\pprime\approx\gamma_mp_0\(\frac{R\rhothetabar}{p_0}\)^{\gamma_m}\frac{\rhothetaprime}{\rhothetabar}.
\end{equation*}
Thus, a numerically preferable perturbation formulation of the Navier-Stokes equations \eqref{NS_equations} reads
\begin{equation}
\begin{aligned}
\rhoprime_t+\nabla\cdot\(\rho\bm u\)&=0,\\	
(\rho\bm u)_t+\nabla\cdot\left(\rho\bm u\otimes\bm u+\pprime\,{\rm Id}-\mu_m\rho\left(\nabla\bm u+(\nabla\bm u)^\top\right)\right)&=
-\rhoprime g\bm{e_d},\\
\rhothetaprime_t+\nabla\cdot\left({\rho\theta}\bm u-\mu_h\rho\nabla\theta\right)&=S_\theta.
\end{aligned}
\label{NS_equations_pert}
\end{equation}

Meteorological applications typically inherit several sources of uncertainty, such as model parameters, initial and boundary conditions.
Consequently, stochastic models need to be designed to analyze the influence of uncertainties on the fluid and cloud dynamics. In general,
there are different ways to represent and take into account model uncertainty. In this paper, we choose a widely used approach in which the
uncertainty is described by random fields. To this end we define an abstract probability space $(\Gamma,\Sigma,P)$ and denote by $\omega$
an event $\omega\in\Gamma$. We assume that the initial data depend on $\omega$, that is,
\begin{equation*}
\rhoprime\big|_{t=0}=\rhoprime(\x,0,\omega),\quad(\rho\uvec)\big|_{t=0}=(\rho\uvec)(\x,0,\omega)\quad\mbox{and}\quad
\rhothetaprime\big|_{t=0}=\rhothetaprime(\x,0,\omega)
\end{equation*}
for the fluid variables and
\begin{equation*}
(\rho q_\ell)\big|_{t=0}=(\rho q_\ell)(\x,0,\omega)\quad\mbox{with}\quad\ell\in\{v,c,r\}
\end{equation*}
for the cloud variables. Consequently, the solution at later time will also depend on $\omega$, that is, we will have
$\rhoprime(\x,t,\omega)$, $(\rho\uvec)(\x,t,\omega)$, $\rhothetaprime(\x,t,\omega)$ and $(\rho q_\ell)(\x,t,\omega)$ for $\ell\in\{v,c,r\}$.
\begin{remark}
The system parameters and boundary conditions could also depend on the random variable. Some of that cases were considered in our previous
work \cite{SG2019}, but in this paper we restrict our attention to the situation in which randomness arises in the initial data only.
\end{remark}

\section{Numerical scheme}\label{numerical_scheme}
In this section, we describe a gPC-SG method for the coupled system \eqref{NS_equations_pert}, \eqref{cloud_equations}. First, in
\S\ref{SUBSEC_system_gpc} we derive a system for the gPC coefficients and then in \S\ref{SUBSEC_discr_gpc} we describe a method used to
numerically solve the resulting system.

\subsection{System of gPC coefficients}\label{SUBSEC_system_gpc}
In the gPC-SG setup, the solution is sought in the form of polynomial expansions
\begin{equation}
\begin{aligned}
&\rhoprime(\x,t,\omega)=\sum_{k=0}^M(\widehat{\rhoprime})_k(\x,t)\Phi_k(\omega),\quad
\rho\uvec(\x,t,\omega)=\sum_{k=0}^M(\widehat{\rho\uvec})_k(\x,t)\Phi_k(\omega),\\
&\rhothetaprime(\x,t,\omega)=\sum_{k=0}^M(\widehat{\rhothetaprime})_k(\x,t)\Phi_k(\omega)
\end{aligned}
\label{expansion_NS_variables}
\end{equation}
and
\begin{equation}
\rho q_\ell(\x,t,\omega)=\sum_{k=0}^M(\widehat{\rho q_\ell})_k(\x,t)\Phi_k(\omega)\quad\mbox{with}\quad\ell\in\{v,c,r\},
\label{expansion_rhoq}
\end{equation}
where $\Phi_k(\omega)$, $\,k=0,\ldots,M$, are polynomials of degree $k$ that are orthogonal with respect to the probability density function
$\mu(\omega)$. Assuming that $\Gamma$ is a compact metric event space, the corresponding Riemann integrals can be defined (see
\cite{Taylor}), and then the orthogonality property implies
\begin{equation}
\int\limits_\Gamma\Phi_k(\omega)\Phi_{k'}(\omega)\mu(\omega)\,{\rm d}\omega=c_k\delta_{kk'}\quad\mbox{for}\quad0\le k,k'\le M,
\label{orthogonal_property}
\end{equation}
where $\delta_{kk'}$ is the Kronecker symbol and $c_k$ are constants depending on the probability density function $\mu$. In this work, we
will focus on two distributions that are important for meteorological applications:
\begin{enumerate}
\item A \textit{uniform distribution} $\mathcal{U}(\Gamma)$ with $\Gamma=[-1,1]$, which corresponds to the Legendre polynomials
\begin{equation*}
\Phi_k(\omega)=\sum_{j=0}^{\lfloor\frac{k}{2}\rfloor}(-1)^j\frac{(2k-2j)!}{(k-j)!(k-2j)!j!2^k}\omega^{k-2j},\quad
\left\lfloor\frac{k}{2}\right\rfloor=\begin{cases}\frac{n}{2},&\mbox{if $n$ is even},\\\frac{n-1}{2},&\mbox{if $n$ is odd},
\end{cases}
\end{equation*}
which are orthogonal with respect to the probability density function $\mu(\omega)=\nicefrac{1}{2}$ and the constants in
\eqref{orthogonal_property} are
\begin{equation}
c_k=\frac{1}{2k+1}.
\label{c_mu_legendre}
\end{equation}
\item A \textit{normal distribution} $\mathcal{N}(\mu_H,\sigma_H^2)$ with $\Gamma=(-\infty,\infty)$ which corresponds to the Hermite
polynomials
\begin{equation*}
\Phi_k(\omega)=2^{-\frac{k}{2}}H_k\(\frac{\omega-\mu_H}{\sqrt{2}\sigma_H}\)\quad\text{with}\quad
H_k(\omega)=(-1)^k\e^{\omega^2}\frac{{\rm d}^k}{{\rm d}\omega^k}\big(e^{-\omega^2}\big),
\end{equation*}
where $\mu_H$ and $\sigma_H$ are the mean value and the standard deviation, respectively. One can show that the Hermite polynomials $\Phi_k$
are orthogonal with respect to $\mu(\omega)=\frac{1}{\sqrt{2\pi\sigma_H^2}}\exp\big(-\nicefrac{({\omega}-\mu_H)^2}{(2\sigma_H^2)}\big)$ and
the constants in \eqref{orthogonal_property} are
\begin{equation}
c_k=k!
\label{c_mu_hermite}
\end{equation}
\end{enumerate}

We use gPC expansions for the source term $S_\theta$ on the right-hand side (RHS) of \eqref{NS_equations_pert},
\begin{equation}
S_\theta(\x,t,\omega)=\sum_{k=0}^M(\widehat{S_\theta})_k(\x,t)\Phi_k(\omega),
\label{expansion_Stheta}
\end{equation}
for the terms on the RHS of \eqref{cloud_equations}
\begin{equation}
\begin{aligned}
\rho(\x,t,\omega)\(-C(\x,t,\omega)+E(\x,t,\omega)\)=:r_1(\x,t,\omega)&=\sum_{k=0}^M(\widehat{r_1})_k(\x,t)\Phi_k(\omega),\\
\rho(\x,t,\omega)\(C(\x,t,\omega)-A_1(\x,t,\omega)-A_2(\x,t,\omega)\)=:r_2(\x,t,\omega)&=\sum_{k=0}^M(\widehat{r_2})_k(\x,t)\Phi_k(\omega),
\\
\rho(\x,t,\omega)\(A_1(\x,t,\omega)+A_2(\x,t,\omega)-E(\x,t,\omega)\)=:r_3(\x,t,\omega)&=\sum_{k=0}^M(\widehat{r_3})_k(\x,t)\Phi_k(\omega),
\end{aligned}
\label{5.4}
\end{equation}
and for the raindrop fall velocity,
\begin{equation}
v_q(\x,t,\omega)=\sum_{k=0}^M(\widehat{v_q})_k(\x,t)\Phi_k(\omega).
\label{5.5}
\end{equation}

Substituting \eqref{expansion_NS_variables}, \eqref{expansion_Stheta} into \eqref{NS_equations_pert} and \eqref{expansion_rhoq},
\eqref{5.4}, \eqref{5.5} into \eqref{cloud_equations}, applying the Galerkin projection and using the orthogonality property
\eqref{orthogonal_property} yield the following deterministic system consisting of $(d+2)(M+1)$ equations for the gPC coefficients of the
fluid variables
\begin{equation}
\begin{aligned}
((\widehat{\rhoprime})_k)_t+\nabla\cdot(\widehat{\rho\uvec})_k&=0,\\
((\widehat{\rho \uvec})_k)_t+\nabla\cdot\(\widehat{\mathfrak{N}}_k+(\widehat{\pprime})_k\text{Id}\)-\mu_m(\widehat{\bm{d_1}})_k&=
-(\widehat{\rhoprime})_kg\bm{e_d},\\
((\widehat{\rhothetaprime})_k)_t+\nabla\cdot\(\thetabar(\widehat{\rho\uvec})_k+\widehat{\bm{\eta}}_k\)-\mu_h(\widehat{d_2})_k&=
(\widehat{S_\theta})_k,
\end{aligned}
\label{NS_equations_coefficients}
\end{equation}
and $3(M+1)$ equations for the gPC coefficients of the cloud variables
\begin{equation}
\begin{aligned}
((\widehat{\rho q_v})_k)_t+\nabla\cdot((\widehat{\bm{\eta_1^q}})_k)-\mu_q(\widehat{{d_1^q}})_k&=(\widehat{r_1})_k,\\
((\widehat{\rho q_c})_k)_t+\nabla\cdot((\widehat{\bm{\eta_2^q}})_k)-\mu_q(\widehat{{d_2^q}})_k&=(\widehat{r_2})_k,\\
((\widehat{\rho q_r})_k)_t+\nabla\cdot((\widehat{\bm{\eta_3^q}})_k)-\mu_q(\widehat{{d_3^q}})_k&=(\widehat{r_3})_k,
\end{aligned}
\label{cloud_equations_coefficients}
\end{equation}
for $k=0,\ldots,M$. It should be pointed that for the pressure term in the momentum equation in \eqref{NS_equations_coefficients} the
expected values for $q_v$, $q_c$ and $q_r$ are used for computing $R_m$, which gives
\begin{equation*}
\pprime=\frac{\gamma_mp_0}{\rhothetabar}\left(\frac{R\rhothetabar}{p_0}\right)^{\gamma_m}\rhothetaprime=
\frac{\gamma_mp_0}{\rhothetabar}\left(\frac{R\rhothetabar}{p_0}\right)^{\gamma_m}\sum_{k=0}^M(\widehat{\rhothetaprime})_k\Phi_k
\end{equation*}
with $R_m=(1-(\widehat{q_v})_0-(\widehat{q_c})_0-(\widehat{q_r})_0)R+(\widehat{q_v})_0R_v$. In this way we keep the linear formulation of
the pressure as in the deterministic case. In what follows, we explain how the expansion coefficients in \eqref{NS_equations_coefficients}
and \eqref{cloud_equations_coefficients} are obtained.

The coefficients $\{(\widehat{r_1})_k,(\widehat{r_2})_k,(\widehat{r_3})_k\}_{k=0}^M$ and $\{(\widehat S_\theta)_k\}_{k=0}^M$ are calculated
via discrete Legendre/Hermite transform (DT) and inverse discrete Legendre/Hermite transform (IDT).
\begin{itemize}
\item {\em Discrete transform (DT)}

The discrete transform starts with the expansion of a function $f$ in the stochastic space
\begin{equation}
f(\x,t,\omega)=\sum_{k=0}^M\widehat f_k(\x,t)\Phi_k(\omega)
\label{dt_exp}
\end{equation}
and by using the orthogonality property \eqref{orthogonal_property} ends up with the expansion coefficients
\begin{equation}
\begin{aligned}
\widehat f_k(\x,t)=\frac{1}{c_k}\int\limits_\Gamma f(\x,t,\omega)\Phi_k(\omega)\mu(\omega)\,{\rm d}\omega\quad\mbox{for}\quad0\le k\le M.
\end{aligned}
\label{dt_2}
\end{equation}
We approximate the above integral using an appropriate Gaussian quadrature rule. We distinguish between the two cases considered in this
paper---Legendre and Hermite polynomials.
\begin{enumerate}
\item For Legendre polynomials the coefficients $c_k$ are determined in \eqref{c_mu_legendre}, $\Gamma=[-1,1]$, and the expansion
coefficients in \eqref{dt_2} are given by
\begin{equation*}
\widehat f_k(\x,t)=\frac{2k+1}{2}\int\limits_{-1}^1f(\x,t,\omega)\Phi_k(\omega)\,{\rm d}\omega\quad\mbox{for}\quad0\le k\le M.
\end{equation*}
Approximating the above integrals using the Gauss-Legendre quadrature leads to
\begin{equation}
\begin{aligned}
\mbox{DT}\left[\left\{f(\x,t,\omega_n)\right\}_{n=0}^N\right]=\left\{\widehat f_k(\x,t)\right\}_{k=0}^M=
\bigg\{\frac{2k+1}{2}\sum_{n=0}^N\beta_n f(\x,t,\omega_n)\Phi_k(\omega_n)\bigg\}_{k=0}^M,
\end{aligned}
\label{DT_legendre}
\end{equation}
where $\omega_n$ and $\beta_n$ are the Gauss-Legendre nodes and weights, respectively.
\item For Hermite polynomials the coefficients $c_k$ are determined in \eqref{c_mu_hermite}, $\Gamma=(-\infty,\infty)$, and the expansion
coefficients in \eqref{dt_2} are given by
\begin{equation*}
\widehat f_k(\x,t)=\frac{1}{k!\sqrt{2\pi\sigma_H^2}}\int\limits_{-\infty}^\infty f(\x,t,\omega)\Phi_k(\omega)
\e^{-\frac{(\omega-\mu_H)^2}{2\sigma_H^2}}{\rm d}\omega\quad\mbox{for}\quad0\le k\le M.
\end{equation*}
Approximating the above integral using the Gauss-Hermite quadrature leads to
\begin{equation}
\hspace*{-0.2cm}\begin{aligned}
\mbox{DT}&\left[\left\{f(\x,t,\omega_n)\right\}_{n=0}^N\right]=\left\{\widehat f_k(\x,t)\right\}_{k=0}^M
=\bigg\{\frac{1}{k!\sqrt{2\pi\sigma_H^2}}\sum_{n=0}^N\beta_nf(\x,t,\omega_n)\Phi_k(\omega_n)
e^{-\frac{(\omega_n-\mu_H)^2}{2\sigma_H^2}+\omega_n^2}\bigg\}_{k=0}^M\!,
\end{aligned}
\label{DT_hermite}
\end{equation}
where $\omega_n$ and $\beta_n$ are the Gauss-Hermite nodes and weights, respectively.
\end{enumerate}
\item {\em Inverse discrete transform (IDT)}

The inverse discrete transform maps the expansion coefficients $\{\widehat f_k(\x,t)\}_{k=0}^M$ to the point values $f(\x,t,\omega_n)$,
$0\le n\le N$. To this end we simply compute the point values of $f$ using the gPC expansion \eqref{dt_exp}:
\begin{equation}
\mbox{IDT}\left[\left\{\widehat f_k(\x,t)\right\}_{k=0}^M\right]=\big\{f(\x,t,\omega_n)\big\}_{n=0}^N=
\left\{\sum_{k=0}^M\widehat f_k(\x,t)\Phi_k(\omega_n)\right\}_{n=0}^N.
\label{IDT}
\end{equation}
\end{itemize}
\begin{remark}
The number of quadrature points $N$ can be chosen equal to the number of expansion coefficients $M$, or even higher for a more accurate
approximation.
\end{remark}
\begin{remark}
We stress that the quadrature weights $\beta_n$ and the values $\Phi_k(\omega_n)$, $0\le k\le M$, $0\le n\le N$, which are used in
\eqref{DT_legendre}--\eqref{IDT}, can be pre-computed for the code efficiency.
\end{remark}

The coefficients $\{(\widehat{q_\ell})_k(\x,t)\}_{k=0}^M$ are computed in the following way:
\begin{equation}
\big\{(\widehat{q_\ell})_k(\x,t)\big\}_{k=0}^M=\mbox{DT}\left[\frac{\mbox{IDT}\left[\{(\widehat{\rho q_\ell})_k(\x,t)\}_{k=0}^M\right]}
{\mbox{IDT}\left[\{(\widehat{\rhoprime})_k(\x,t)\}_{k=0}^M\right]+\rhobar(\x)}\right],\quad\ell\in\{v,c,r\},
\label{IDT_DT}
\end{equation}
where the DT and IDT operators are defined in \eqref{DT_legendre}--\eqref{IDT}. Additionally, the coefficients
$\{\widehat{\mathfrak{N}}_k\}_{k=0}^M$, $\{\bm{\widehat{\eta}}_k\}_{k=0}^M$, $\{(\bm{\widehat{d_1}})_k\}_{k=0}^M$ and
$\{(\widehat{d_2})_k\}_{k=0}^M$ appearing in \eqref{NS_equations_coefficients} and the coefficients
$\{({\widehat{\bm{\eta_1^q}}})_k,({\widehat{\bm{\eta_2^q}}})_k,(\widehat{\bm{\eta_3^q}})_k\}_{k=0}^M$ and
$\{(\widehat{{d_1^q}})_k,(\widehat{{d_2^q}})_k,(\widehat{{d_3^q}})_k\}_{k=0}^M$ appearing in \eqref{cloud_equations_coefficients} are
obtained through the DT and IDT in a similar manner as in \eqref{IDT_DT}. Here, the nonlinear Navier-Stokes advection coefficients
$\{\widehat{\mathfrak{N}}_k\}_{k=0}^M$ and $\{\bm{\widehat{\eta}}_k\}_{k=0}^M$ are obtained from
\begin{equation}
\nabla\cdot\(\rho\bm u\otimes\bm u\)=\sum_{k=0}^M\nabla\cdot\widehat{\mathfrak{N}}_k\Phi_k,\quad
\nabla\cdot\(\thetaprime\rho\uvec\)=\sum_{k=0}^M\nabla\cdot\widehat{\bm\eta}_k\Phi_k,
\label{eta_i}
\end{equation}
and the nonlinear diffusion coefficients $\{(\widehat{\bm{d_1}})_k\}_{k=0}^M$ and $\{(\widehat{d_2})_k\}_{k=0}^M$ are obtained from
\begin{equation}
\nabla\cdot\left(\rho\left(\nabla\bm u+(\nabla\bm u)^T\right)\right)=\sum_{k=0}^M(\widehat{\bm{d_1}})_k\Phi_k,\quad
\nabla\cdot\left(\rho\nabla\theta\right)=\sum_{k=0}^M(\widehat{d_2})_k\Phi_k.
\label{d_i}
\end{equation}
Finally, the nonlinear cloud dynamics advection coefficients
$\{({\widehat{\bm{\eta_1^q}}})_k,({\widehat{\bm{\eta_2^q}}})_k,(\widehat{\bm{\eta_3^q}})_k\}_{k=0}^M$ are obtained from
\begin{equation}
\begin{aligned}
&\nabla\cdot(\rho q_v\uvec)=\sum_{k=0}^M\nabla\cdot(\widehat{\bm{\eta_1^q}})_k\Phi_k,\quad
\nabla\cdot(\rho q_c\uvec)=\sum_{k=0}^M\nabla\cdot(\widehat{\bm{\eta_2^q}})_k\Phi_k,\\
&\nabla\cdot(\rho q_r\uvec-\rho q_rv_q\bm{e_d})=\sum_{k=0}^M\nabla\cdot(\widehat{\bm{\eta_3^q}})_k\Phi_k,
\end{aligned}
\label{cloud_eta_i}
\end{equation}
and the nonlinear diffusion coefficients $\{(\widehat{{d_1^q}})_k,(\widehat{{d_2^q}})_k,(\widehat{{d_3^q}})_k\}_{k=0}^M$ are obtained from
\begin{equation}
\begin{aligned}
\nabla\cdot(\rho\nabla q_v)=\sum_{k=0}^M(\widehat{{d_1^q}})_k\Phi_k,\quad
\nabla\cdot(\rho\nabla q_c)=\sum_{k=0}^M(\widehat{{d_2^q}})_k\Phi_k,\quad
\nabla\cdot(\rho\nabla q_r)=\sum_{k=0}^M(\widehat{{d_3^q}})_k\Phi_k.
\end{aligned}
\label{cloud_d_i}
\end{equation}

\subsection{Discretization of the gPC system}\label{SUBSEC_discr_gpc}
In this section we describe the numerical methods used to solve the resulting gPC system \eqref{NS_equations_coefficients},
\eqref{cloud_equations_coefficients}. To this end we denote by
\begin{equation*}
\widehat{\w}:=\Big(\underline{\widehat{\rhoprime}},\underline{\widehat{\rho \uvec}},\underline{\widehat{\rhothetaprime}}\Big)^\top\quad
\mbox{and}\quad
\widehat{\w_q}:=\big(\underline{\widehat{\rho q_v}},\underline{\widehat{\rho q_c}},\underline{\widehat{\rho q_r}}\big)^\top
\end{equation*}
the solution vectors of \eqref{NS_equations_coefficients} and \eqref{cloud_equations_coefficients}, respectively. Here, the underline
$\underline{(\cdot)}$ denotes the vector of the respective coefficients. For instance, for the solution coefficients we have
\begin{equation*}
\begin{aligned}
&\underline{\widehat{\rhoprime}}:=\big((\widehat{\rhoprime})_0,\ldots,(\widehat{\rhoprime})_M\big),\quad
\underline{\widehat{\rho\uvec}}:=\big((\widehat{\rho u_1})_0,\ldots,(\widehat{\rho u_d})_0;\ldots;
(\widehat{\rho u_1})_M,\ldots,(\widehat{\rho u_d})_M\big),\\
&\underline{\widehat{\rhothetaprime}}:=\big((\widehat{\rhothetaprime})_0,\ldots,(\widehat{\rhothetaprime})_M\big),\quad
\underline{\widehat{\rho q_\ell}}=\big((\widehat{\rho q_\ell})_0,\ldots,(\widehat{\rho q_\ell})_M\big),~~\ell\in\{v,c,r\}.
\end{aligned}
\end{equation*}
Then, the coupled system can be written as
\begin{align}
&\widehat{\w}_t=-\nabla\cdot\F(\widehat{\w})+\D(\widehat{\w})+\R(\widehat{\w}),\label{operator_NS_coeff}\\
&(\widehat{\w_q})_t=-\nabla\cdot\F_q(\widehat{\w_q})+\D_q(\widehat{\w_q})+\R_q(\widehat{\w_q}),\label{operator_cloud_coeff}
\end{align}
where $\F$ and $\F_q$ are convective fluxes and $\D$, $\R$ and $\D_q$, $\R_q$ denote the diffusion and reaction operators of the respective
systems. They are given by
\begin{align}
&\F(\widehat{\w}):=\big(\underline{\widehat{\rho\uvec}},\underline{\widehat{\pprime}\text{Id}}+\underline{\widehat{\mathfrak{N}}},
\thetabar\underline{\widehat{\rho\uvec}}+\underline{\bm{\widehat{\eta}}}\big)^\top,\quad
\D(\widehat{\w}):=\big(0,\mu_m\underline{\widehat{\bm{d_1}}},\mu_h\underline{\widehat{d_2}}\big)^\top,\quad
\R(\widehat{\w}):=\big(0,-\underline{\widehat{\rhoprime} g\bm{e_d}},\underline{\widehat{S_\theta}}\big)^\top,\label{NS_operator}\\
&\F_q(\widehat{\w_q}):=\big(\underline{\widehat{\bm{\eta_1^q}}},\underline{\widehat{\bm{\eta_2^q}}},
\underline{\widehat{\bm{\eta_3^q}}}\big)^\top,\quad
\D_q(\widehat{\w_q}):=\mu_q\big(\underline{\widehat{{d_1^q}}},\underline{\widehat{{d_2^q}}},\underline{\widehat{{d_3^q}}}\big)^\top,\quad
\R_q(\widehat{\w_q}):=\big(\underline{\widehat{r_1}},\underline{\widehat{r_2}},\underline{\widehat{r_3}}\big)^\top,\label{cloud_operator}
\end{align}
where the respective components of the above vectors were defined in \eqref{eta_i}--\eqref{cloud_d_i}. It should be observed that
\eqref{operator_NS_coeff}--\eqref{cloud_operator} is a deterministic system for the expansion coefficients. This system has the same
structure as the deterministic system studied in \cite{SG2019}. Therefore, one can directly apply the finite volume method from
\cite{SG2019} for the spatial discretization of the system \eqref{operator_NS_coeff}--\eqref{cloud_operator} taking into account that the
number of equations has increased by a factor of $M+1$ and that the DT and IDT should be applied for each evaluation of the nonlinear terms
on the RHS of \eqref{operator_NS_coeff} and \eqref{operator_cloud_coeff}.

Concerning the time discretization, we also follow the approach presented in \cite{SG2019} and implement a Strang splitting approach between
the Navier-Stokes and cloud dynamics systems. We then approximate the cloud dynamics \eqref{operator_cloud_coeff}, \eqref{cloud_operator}
with a third-order large stability domain explicit Runge-Kutta method (DUMKA3 see \cite{Dumka3,Dumka3_code}) and the fluid dynamics
\eqref{operator_NS_coeff}, \eqref{NS_operator} with the IMEX AP ARS(2,2,2) method from \cite{ARS222}. For the latter, we split the
Navier-Stokes operators in \eqref{NS_operator} into two parts
\begin{equation*}
\begin{aligned}
&\F(\widehat{\w})=\F_L(\widehat{\w})+\F_N(\widehat{\w})\quad\mbox{with}\quad\F_L(\widehat{\w}):=
\big(\underline{\widehat{\rho\uvec}},\underline{\widehat{\pprime}\text{Id}},\thetabar(\underline{\widehat{\rho\uvec}})\big)^\top~~\mbox{and}
~~\F_N(\widehat{\w}):=\big(0,\underline{\widehat{\mathfrak{N}}},\underline{\bm{\widehat{\eta}}}\big)^\top,\\
&\R(\widehat{\w})=\R_L(\widehat{\w})+\R_N(\widehat{\w})\quad\mbox{with}\quad\R_L(\widehat{\w}):=
\big(0,-\underline{\widehat{\rhoprime}g\bm{e_d}},0\big)^\top~~\mbox{and}~~
\R_N(\widehat{\w}):=\big(0,0,\underline{\widehat{S_\theta}}\big)^\top,
\end{aligned}
\end{equation*}
and then define the stiff linear operator $\mathcal{L}:=-\nabla\cdot\F_L(\widehat{\w})+\R_L(\widehat{\w})$ and the nonstiff nonlinear
operator $\mathcal{N}:=-\nabla\cdot\F_N(\widehat{\w})+\D(\widehat{\w})+\R_N(\widehat{\w})$, which are treated implicitly and explicitely,
respectively.

Finally, we note that as in \cite{SG2019}, the time steps for both the Navier-Stokes and cloud dynamics splitting substeps are chosen
adaptively at every time level.
%and for the case considered here should satisfy the following conditions:
%\begin{equation*}
%\begin{aligned}
%&\max\(\frac{\max(\mu_h,\mu_m)}{h^2},\max_{n=0,\ldots,N}\,\max_{s=1,.., d}\,\max_{i=1,\ldots,N_x}(|(u_s)_i(\omega_n)|)\,\frac{d}{h}\)\,\Delta t_{\rm NS}<0.5,\\
%&\max\(\frac{\mu_q}{h^2},\max_{n=0,\ldots,N}\,\max_{s=1,.., d-1}\,\max_{i=1,\ldots,N_x}(|(u_s)_i(\omega_n)|,|(u_d)_i(\omega_n)+v_q(\omega_n)|)\,\frac{d}{h}\)\,\Delta t_{\rm cloud}<0.5,
%\end{aligned}
%\end{equation*}
%where $h$ is the spatial mesh size and $N_x$ is the total number of mesh cells.

\section{Numerical experiments}\label{numerical_experiments}
In this section, we present experimental results for the fully random Navier-Stokes-cloud model \eqref{NS_equations_pert},
\eqref{cloud_equations}. In Examples \ref{example6} and \ref{example7}, we investigate the experimental convergence of our numerical scheme
using the well-known meteorological benchmark describing the free convection of a smooth warm and moist air bubble; see, e.g.,
\cite{BF02,DT50}. In Example \ref{example6}, we demonstrate the spatio-temporal convergence as well as the convergence in the stochastic
space for the case in which the initial vapor concentration $q_v$ is perturbed by $10\%$ which is realized with a uniform distribution of
the randomness. In Example \ref{example7}, we investigate the convergence for the same setup as in Example \ref{example6} but with normally
distributed randomness. Since we use the same numerical method for the space and time discretization as in Example \ref{example6}, in
Example \ref{example7} we just investigate the convergence in the stochastic space. In Examples \ref{RB_stoch_example2} and
\ref{RB_stoch_example4}, we present the results of the uncertainty study for the Rayleigh-B\'{e}nard convection in both 2-D and 3-D. We also
compare the results obtained in this work for the fully random model with the deterministic one, in which both the Navier-Stokes euqations
\eqref{NS_equations_pert} and the cloud equations \eqref{cloud_equations} are deterministic, and with the semi-random one, in which the
deterministic Navier-Stokes equations are coupled with the random cloud dynamics (this semi-random model was studied in the first part of
this work in \cite{SG2019}). We note that the parameters used in the numerical experiments presented below are slightly different from the
one used in \cite{SG2019}; however, the main qualitative features remain the same. In both Rayleigh-B\'{e}nard experiments (Examples
\ref{RB_stoch_example2} and \ref{RB_stoch_example4}), we investigate uncertainty propagation, which is triggered by the initial data of the
water vapor concentration $q_v$ which we perturbed uniformly. A comparison with a normally distributed initial perturbation or even
perturbations of certain parameters is beyond the scope of this work and is left for future study.

In all of the following examples we set $\mu_m=10^{-3}$ and $\mu_h=\mu_q=10^{-2}$ in \eqref{NS_equations_pert} and \eqref{cloud_equations}.

\begin{example}[\textbf{Stochastic initial data with uniformly distributed perturbation}]\label{example6}
\phantom{.}

In this experiment, we simulate free convection of a smooth warm and moist air bubble in 2-D. Due to the shear friction with the surrounding
air at the warm/cold air interface, the warm air bubble rises and deforms axisymmetrically and gradually forms a mushroom-like shape. The
bubble is placed at $(2500\,\mathrm{m},2000\,\mathrm{m})$ in a domain $\Omega=[0,5000]\times[0,5000]\,m^2$. We consider a 10\% perturbation
of the initial water vapor concentration. This is realized through the following initial conditions in the case of a uniformly distributed
randomness for the cloud variables:
$$
\begin{aligned}
(\widehat{q_v})_0(\x,0)&=0.005\,\thetaprime(\x,0),~~(\widehat{q_v})_1(\x,0)=0.1,~~(\widehat{q_v})_0(\x,0),~~(\widehat{q_v})_k(\x,0)=0~\,
\mbox{for}\;2\le k\le M,\\
(\widehat{q_c})_0(\x,0)&=10^{-4}\,\thetaprime(\x,0),~~(\widehat{q_c})_k(\x,0)=0~\,\mbox{for}\;1\le k\le M,\\
(\widehat{q_r})_0(\x,0)&=10^{-6}\,\thetaprime(\x,0),~~(\widehat{q_r})_k(\x,0)=0~\,\mbox{for}\;1\le k\le M,
\end{aligned}
$$
and for the Navier-Stokes variables:
\begin{equation*}
\begin{aligned}
&(\widehat{\rhoprime})_0(\x,0)=-\rhobar(\x)\frac{(\widehat{\thetaprime})_0(\x,0)}{\thetabar(\x)+(\widehat{\thetaprime})_0(\x,0)},~~
(\widehat{\rhoprime})_k(\x,0)=0~\,\mbox{for}\;1\le k\le M,\\
&(\widehat{\rho\bm{u}})_k(\x,0)=0~\,\mbox{for}\;0\le k\le M,\\
&(\widehat{\rhothetaprime})_0(\x,0)=\bar{\rho}(\x)(\widehat{\thetaprime})_0(\x,0)+\bar{\theta}(\widehat{\rhoprime})_0(\x,0)+
(\widehat{\thetaprime})_0(\x,0)(\widehat{\rhoprime})_0(\x,0),~~(\widehat{\rhothetaprime})_k(\x,0)=0~\,\mbox{for}\;1\le k\le M,
\end{aligned}
\end{equation*}
where
$$
(\widehat{\thetaprime})_0(\x,0)=\begin{cases}
2\cos^2\left(\dfrac{\pi r}{2}\right),&r:=\sqrt{(x_1-2500)^2+(x_2-2000)^2}\le2000,\\
0,&\text{otherwise},\\
\end{cases}\quad
(\widehat{\thetaprime})_k(\x,0)=0~\,\mbox{for}\;1\le k\le M.
$$
Additionally, we set $\thetabar=285\,\mathrm{K}$, $p_0=\pbar=10^5\,\mathrm{Pa}$ and
\begin{align*}
\rhobar(\x)=
\frac{p_0}{R\thetabar(\x)}\left(1-\frac{g x_2}{c_p\thetabar}\right)^\frac{1}{\gamma-1}
\end{align*}
with $c_p=1005\,\mathrm{\nicefrac{J}{(kg \cdot K)}}$, $c_v=718\,\mathrm{\nicefrac{J}{(kg \cdot K)}}$ and $\gamma=c_p/c_v$. We start here
with nonzero values for the cloud drops concentration $q_c$ and the rain concentration $q_r$ to avoid values close to the machine precision
since the main purpose of the test is the convergence study. Furthermore, we apply the no-slip boundary conditions for the velocities and
homogeneous Neumann boundary conditions for the remaining variables, that is, $\nabla\rhoprime\cdot\bm n=0$,
$\nabla\rhothetaprime\cdot\bm n=0$ and $\nabla(\rho q_\ell)\cdot\bm n=0$, $\,\ell\in\{v,c,r\}$.

In Figure \ref{SWAB_stoch_full_legendre_T150_T200}, we depict the expected values of the potential temperature $\theta$ and the cloud
variables $q_v$, $q_c$ and $q_r$, computed using a $160\times160$ uniform mesh at time $200s$ with $M=L=3$. For comparison purposes, in
Figure \ref{comp_fullu1} we show the potential temperature $\theta$ and the water vapor concentration $q_v$ computed using the deterministic
Navier-Stokes-cloud model and the potential temperature $\theta$ and the expected values of the water vapor concentration $q_v$ computed
using the semi-random Navier-Stokes-cloud model. Note that for a better comparison, we have used the same vertical scales for presenting the
results in Figures \ref{SWAB_stoch_full_legendre_T150_T200} and \ref{comp_fullu1}. It can be observed that the fully random results are more
smeared compared to the deterministic ones and that in the fully random experiment no additional vortices beneath the bubble have been
developed and the results are slight variations of the deterministic ones, which is to be expected. In order to investigate the appearance
of the vortices in the semi-random case (see the second row Figure \ref{comp_fullu1}), we depict in Figure \ref{P1P5} the semi-random
results obtained for the smaller initial water vapor perturbation taken as 1\%, 5\% and 7\%. One can clearly see, the vortices develop
gradually with higher perturbation and that the 1\% perturbation results are very close to the deterministic ones. Thus, the vortex features
of the solutions obtained with the semi-random model seem to result from the missing feedback to the dynamics of the fluid and are not a
defect of the numerical method. We also note that in the semi-random model, the energy conservation is (slightly) violated. This might lead
to the differences in the simulations, since the dynamics is then driven by the averaged latent heat release and not by the one in the
realization.
\begin{figure}[ht!]
\centerline{\includegraphics[trim=0.0cm 0.1cm 0.2cm 0.1cm,clip,width=6.0cm]{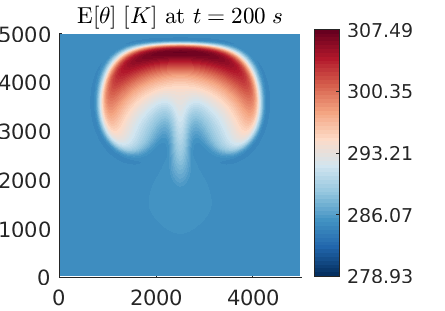}
\hspace*{1.0cm}
\includegraphics[trim=0.0cm 0.1cm 0.2cm 0.1cm,clip,width=6.0cm]{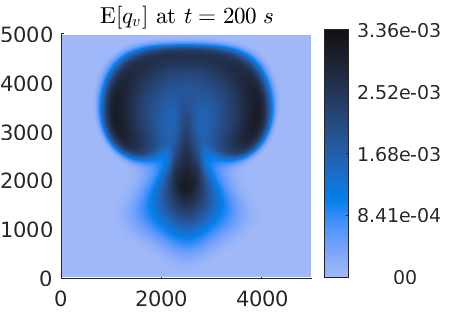}}
\vskip8pt
\centerline{\includegraphics[trim=0.0cm 0.1cm 0.2cm 0.1cm,clip,width=6.0cm]{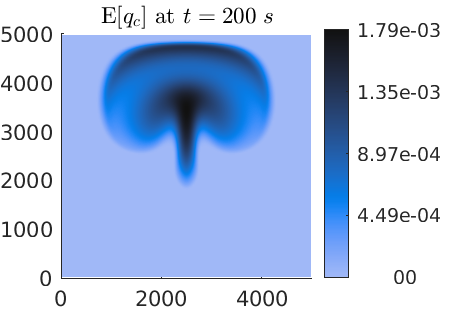}\hspace*{1.0cm}
\includegraphics[trim=0.0cm 0.1cm 0.2cm 0.1cm,clip,width=6.0cm]{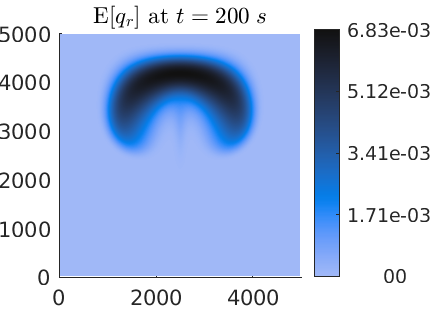}}
\caption {Example \ref{example6}: The expected values of the potential temperature $\theta$, the water vapor concentration $q_v$, cloud
drops concentration $q_c$ and rain concentration $q_r$ computed using the fully random model.\label{SWAB_stoch_full_legendre_T150_T200}}
\end{figure}
\begin{figure}[ht!]
\centerline{\includegraphics[trim=0.0cm 0.1cm 0.2cm 0.1cm,clip,width=6.0cm]{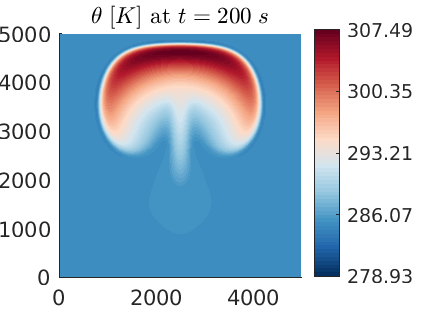}\hspace*{1.0cm}
\includegraphics[trim=0.0cm 0.1cm 0.2cm 0.1cm,clip,width=6.0cm]{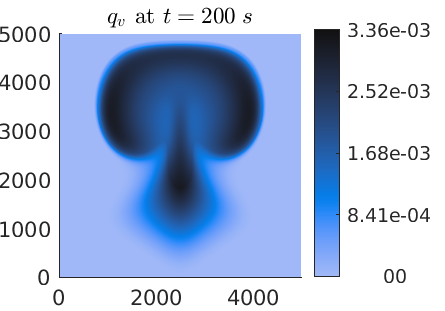}}
\vskip8pt
\centerline{\includegraphics[trim=0.0cm 0.1cm 0.2cm 0.1cm,clip,width=6.0cm]{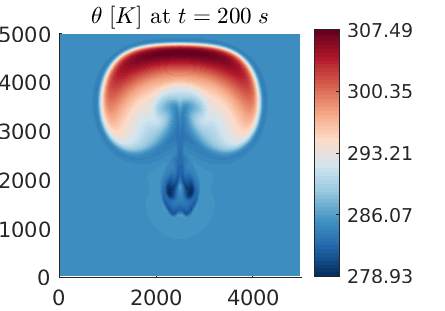}
\hspace*{1.0cm}
\includegraphics[trim=0.0cm 0.1cm 0.2cm 0.1cm,clip,width=6.0cm]{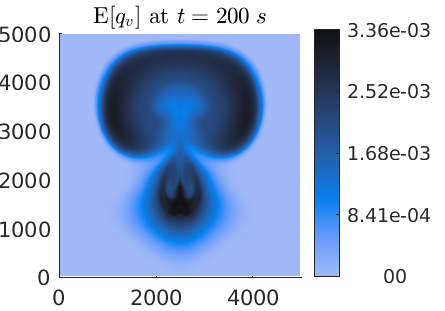}}
\caption {Example \ref{example6}: The potential temperature $\theta$ and the water vapor concentration $q_v$ computed using the
deterministic model (first row) and the potential temperature $\theta$ and the expected values of the water vapor concentration $q_v$
computed using the semi-random model (second row).\label{comp_fullu1}}
\end{figure}
\begin{figure}[ht!]
\centerline{\includegraphics[trim=0.0cm 0.1cm 0.2cm 0.1cm,clip,width=6.0cm]{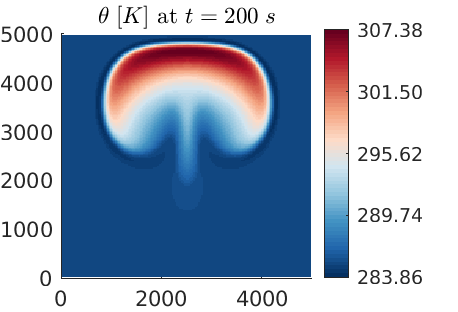}\hspace*{1.0cm}
\includegraphics[trim=0.0cm 0.1cm 0.2cm 0.1cm,clip,width=6.0cm]{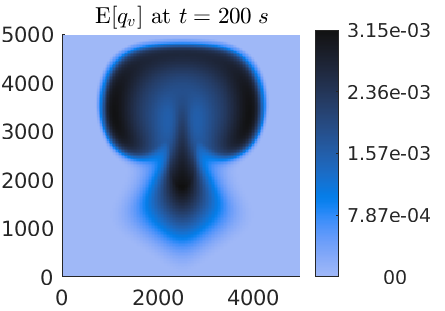}}
\vskip8pt
\centerline{\includegraphics[trim=0.0cm 0.1cm 0.2cm 0.1cm,clip,width=6.0cm]{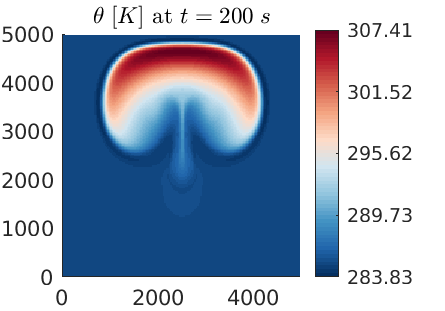}\hspace*{1.0cm}
\includegraphics[trim=0.0cm 0.1cm 0.2cm 0.1cm,clip,width=6.0cm]{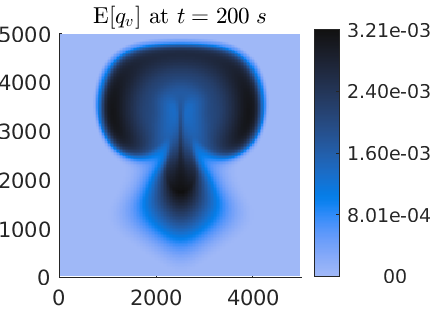}}
\vskip8pt
\centerline{\includegraphics[trim=0.0cm 0.1cm 0.2cm 0.1cm,clip,width=6.0cm]{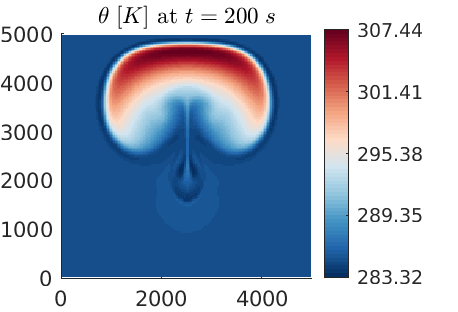}\hspace*{1.0cm}
\includegraphics[trim=0.0cm 0.1cm 0.2cm 0.1cm,clip,width=6.0cm]{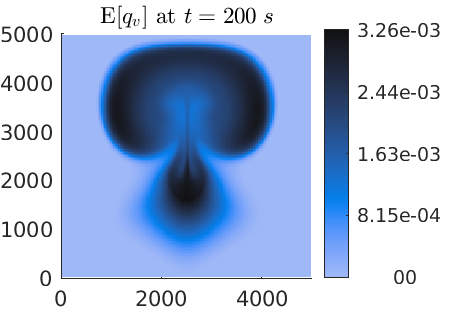}}
\caption {Example \ref{example6}: The potential temperature $\theta$ (left column) and the expected value of the water vapor concentration
$q_v$ (right column) computed using the semi-random model for 1\% (first row), 5\% (second row) and 7\% (third row) of initial water vapor
perturbations.\label{P1P5}}
\end{figure}

In Figure \ref{stoch_full_Ndt_T10}, we present the spatio-temporal convergence study for the expected values of the cloud and flow
variables at the time $t=10s$. We compute the solutions on different $N\times N$ uniform meshes with $M=L=3$. As in the deterministic case
presented in \cite{SG2019}, one can clearly see a second-order convergence for the studied fully random model.
\begin{figure}[ht!]
\centerline{\includegraphics[scale=0.80]{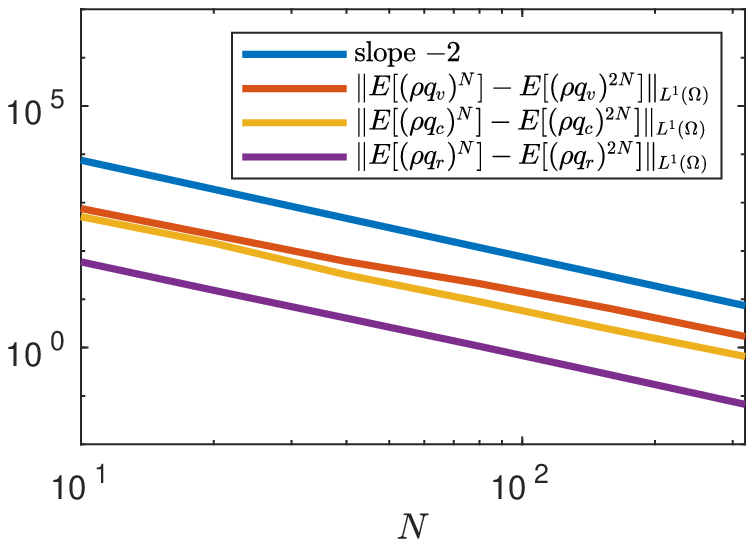}\hspace*{1.0cm}
\includegraphics[scale=0.80]{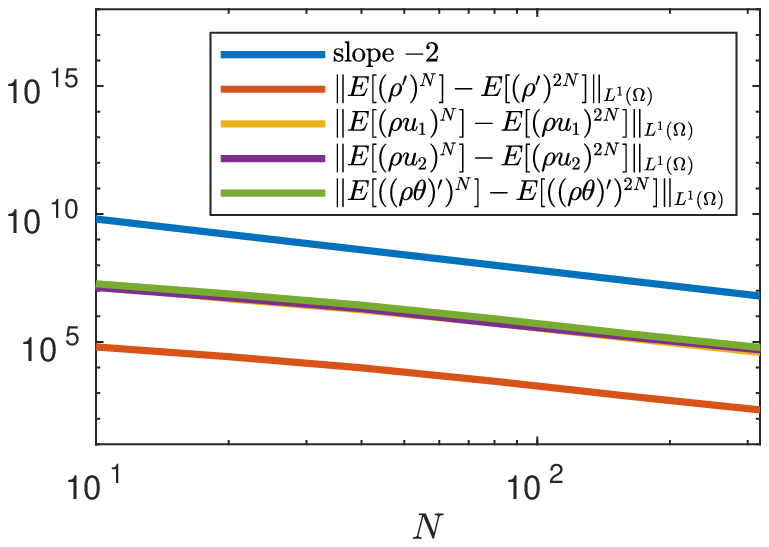}}
\caption{Example \ref{example6}: Spatio-temporal $L^1$ convergence study for the expected values of the cloud variables $q_v$ $q_c$ and
$q_r$ (left) and the flow variables $\rhoprime$, $\rho u_1$, $\rho u_2$ and $\rhothetaprime$ (right) computed at time $t=10s$ using the
fully random model with the constant time step $\Delta t=\nicefrac{256}{100N}$.\label{stoch_full_Ndt_T10}}
\end{figure}

The stochastic convergence studies are presented in Figures \ref{stoch_full_NS_M_L1_T10} and \ref{stoch_full_cloud_M_L1_T10} for the cloud
and Navier-Stokes variables, respectively, at time $t=10s$ using a $160\times160$ uniform mesh. We compute the difference between the
approximate solutions with different numbers of modes $M$ and $L=M$ and the reference solution obtained with 20 stochastic modes and $L=19$.
One can observe a spectral convergence with an approximate rate of $e^{-0.3M}$. One can also see that the error of the rain drops in Figure
\ref{stoch_full_cloud_M_L1_T10} (right) basically stays constant at some point because in this case it approaches the machine precision.
\begin{figure}[ht!]
\centerline{\includegraphics[scale=0.80]{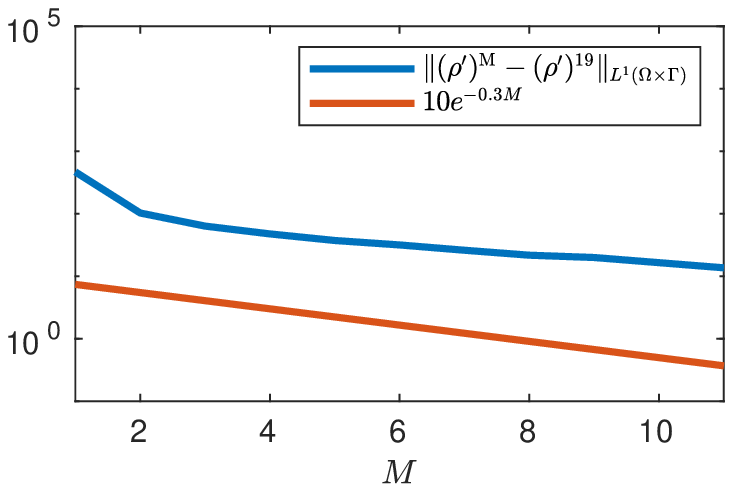}\hspace*{1.0cm}
\includegraphics[scale=0.80]{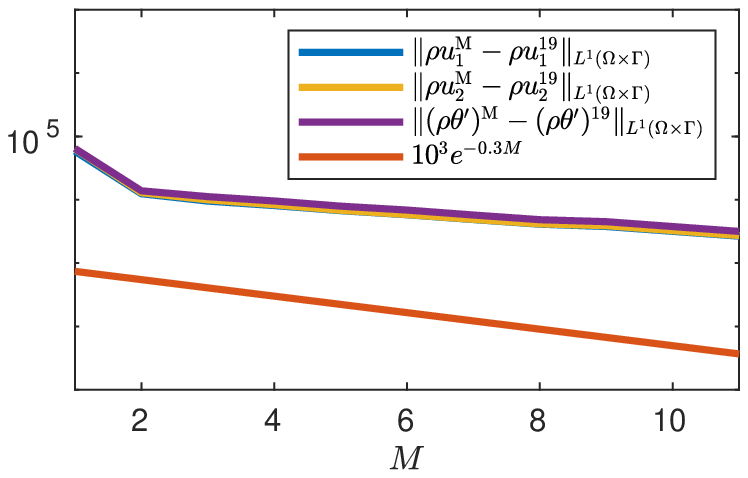}}
\caption{Example \ref{example6}: $L^1$ convergence study for the Navier-Stokes variables $\rhoprime$ (left) and $\rho u_1$, $\rho u_2$ and
$\rhothetaprime$ (right) in the stochastic space computed at time $t=10s$ using the fully random model with the constant time step
$\Delta t=0.01$.\label{stoch_full_NS_M_L1_T10}}
\end{figure}
\begin{figure}[ht!]
\centerline{\includegraphics[scale=0.80]{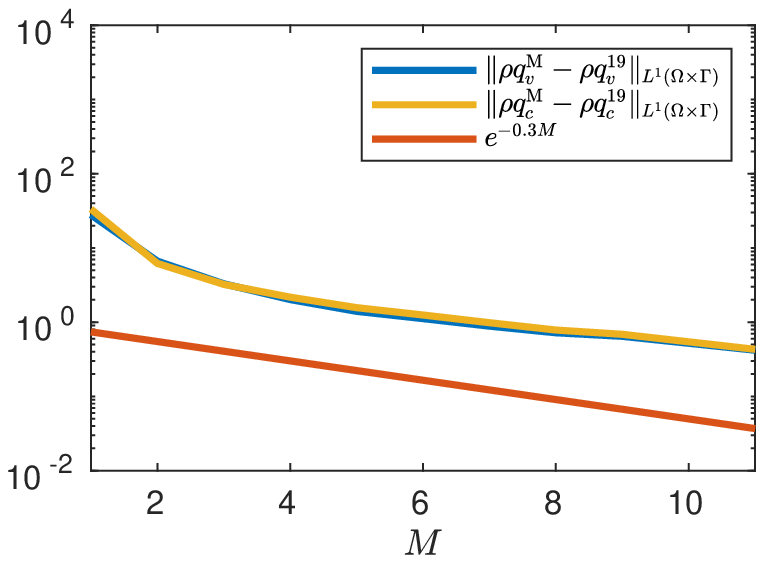}\hspace*{1.0cm}
\hspace*{0.1cm}\includegraphics[scale=0.80]{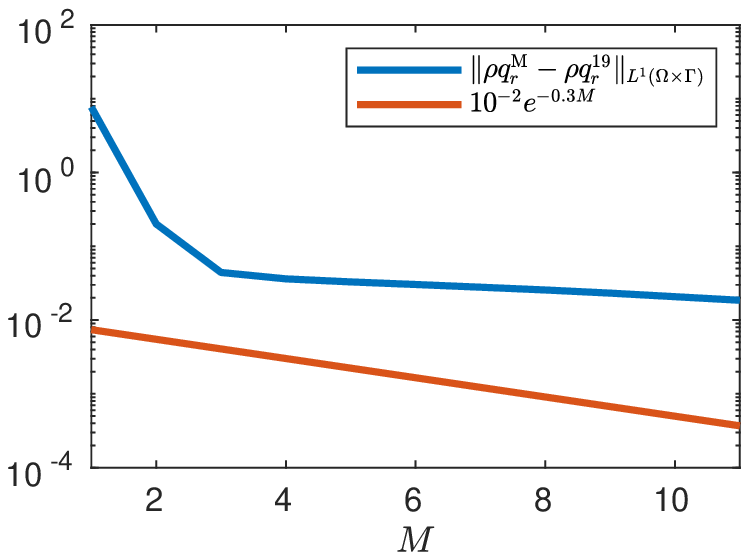}}
\caption{Example \ref{example6}: $L^1$ convergence study for the cloud variables $\rho q_v$ and $q_c$ (left) and $q_r$ (right) in the
stochastic space computed at time $t=10s$ using the fully random model with the constant time step $\Delta t=0.01$.
\label{stoch_full_cloud_M_L1_T10}}
\end{figure}
\end{example}

\begin{example}[\textbf{Stochastic initial data with normally distributed perturbation}]\label{example7}
\phantom{.}

In this experiment, we demonstrate that the convergence of the stochastic Galerkin method for the fully stochastic model does not depend on
the choice of the distribution of the randomness. For this purpose, we choose the same initial conditions as in Example \ref{example6}, but
this time with a normally $\mathcal{N}(0,1)$ distributed perturbation.

In Figure \ref{comp_fulln1}, we compare the solutions (the potential temperature $\theta$ and the water vapor concentration $q_v$) computed
using the deterministic, semi-random and fully random Navier-Stokes-cloud models. For a better comparison, we have used the same vertical
scales for presenting the results. As in the previous example, one can observe that in the fully random experiment no additional vortices
beneath the bubble have been developed and the results are slight variations of the deterministic ones. Thus, the vortex features of the
semi-random results are independent of the distributions of the initial perturbation and caused by the missing feedback to the dynamics of
the fluid.
\begin{figure}[ht!]
\centerline{\includegraphics[trim=0.0cm 0.1cm 0.2cm 0.1cm,clip,width=6.0cm]{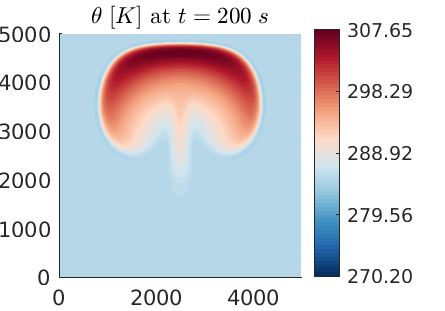}\hspace*{1.0cm}
\includegraphics[trim=0.0cm 0.1cm 0.2cm 0.1cm,clip,width=6.0cm]{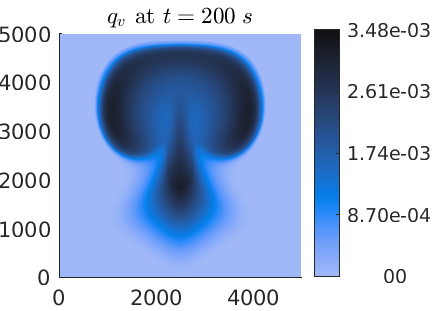}}
\vskip8pt
\centerline{\includegraphics[trim=0.0cm 0.1cm 0.2cm 0.1cm,clip,width=6.0cm]{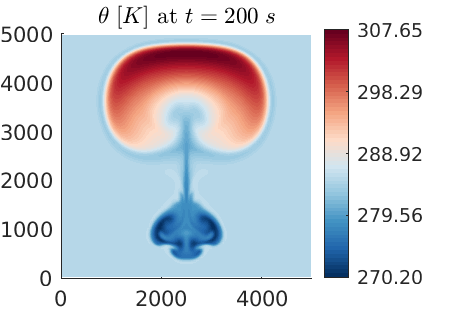}
\hspace*{1.0cm}
\includegraphics[trim=0.0cm 0.1cm 0.2cm 0.1cm,clip,width=6.0cm]{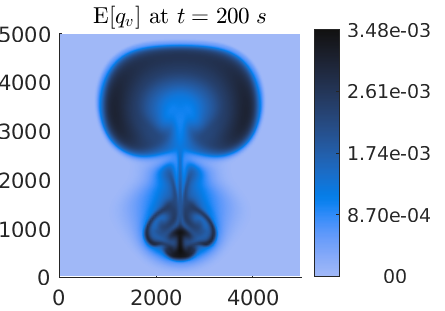}}
\vskip8pt
\centerline{\includegraphics[trim=0.0cm 0.1cm 0.2cm 0.1cm,clip,width=6.0cm]{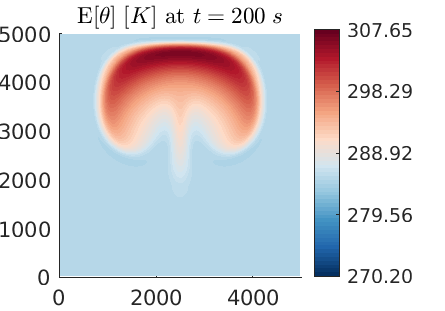}
\hspace*{1.0cm}
\includegraphics[trim=0.0cm 0.1cm 0.2cm 0.1cm,clip,width=6.0cm]{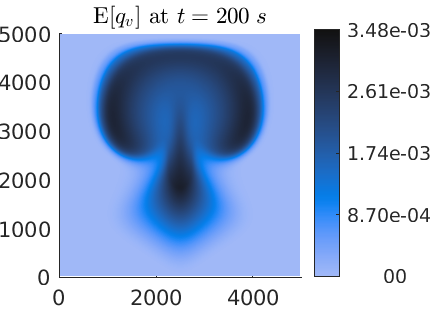}}
\caption {Example \ref{example7}: The potential temperature $\theta$ and the water vapor concentration $q_v$ computed using the
deterministic model (first row); the potential temperature $\theta$ and the expected values of the water vapor concentration $q_v$
computed using the semi-random model (second row); the expected values of the potential temperature $\theta$ and the water vapor
concentration $q_v$ computed using the fully random model (third row).\label{comp_fulln1}}
\end{figure}

Next, we investigate the influence of the choice of distribution for the initial perturbation. In Figure \ref{comp_un}, we depict the
he cloud drops concentration $q_c$ computed using the fully random model with the initial normally and uniformly distributed perturbations;
the latter one was computed in Example \ref{example6}. For a better comparison, we have used the same vertical scales for presenting the
results. Since the initial perturbation is rather small, the results look very alike. In general, both simulations smear the boundaries of
the bubble. However, the smearing with the normal distribution is not as strong as with the uniform distribution. This effect is due to the
concentrated shape of the normal distribution around the expected value; thus, the different realizations are closer to the averaged
potential temperature as a feedback to the energy equation.
\begin{figure}[ht!]
\centerline{\includegraphics[trim=0.0cm 0.1cm 0.2cm 0.1cm,clip,width=6.0cm]{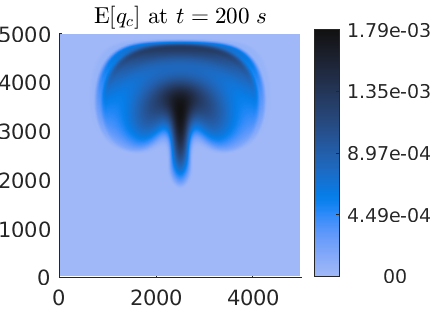}
\hspace*{1.0cm}
\includegraphics[trim=0.0cm 0.1cm 0.2cm 0.1cm,clip,width=6.0cm]{SWAB_2D_stoch_full_legendre_qc_T200}}
\caption {Examples \ref{example6} and \ref{example7}: The expected values of the cloud drops concentration $q_c$ computed using the fully
random model with the initial normally (left) and uniformly (right) distributed perturbations.\label{comp_un}}
\end{figure}

The convergence studies in the stochastic space are presented in Figures \ref{stoch_hermite_full_NS_M_L1_T10} and
\ref{stoch_hermite_full_cloud_M_L1_T10} for the cloud and Navier-Stokes variables, respectively, at time $t=10s$ using a $160\times160$
uniform mesh. We computed the difference between the approximate solutions with different numbers of modes $M$ and $L=M$ and the reference
solution obtained with 12 stochastic modes and $L=11$. As in the case with a uniform distribution studied in Example \ref{example6}, one can
observe a spectral convergence with an approximate rate of $e^{-0.3M}$. This demonstrates that the experimental convergence rate is
independent of the chosen distribution.
\begin{figure}[ht!]
\centerline{\includegraphics[scale=0.80]{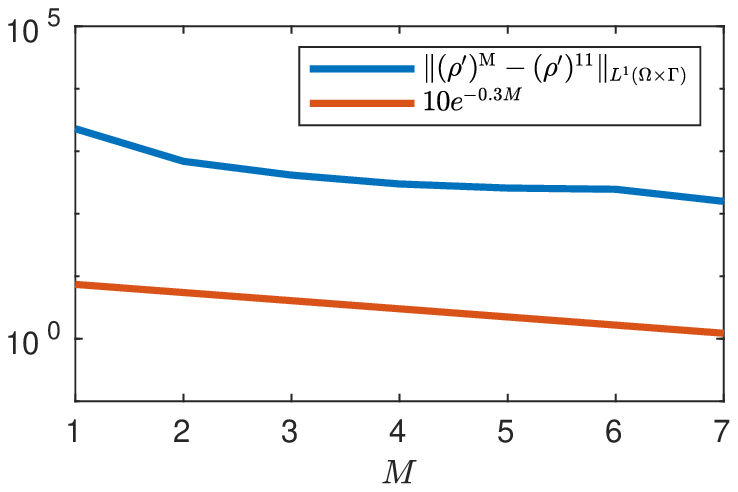}\hspace*{1.0cm}
\includegraphics[scale=0.80]{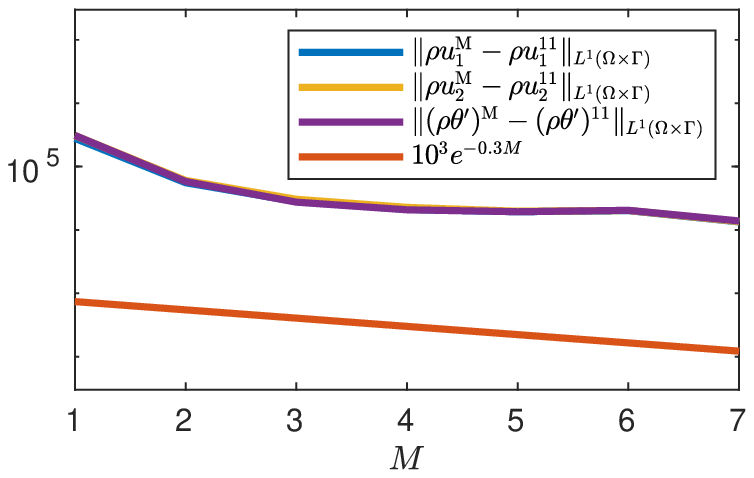}}
\caption{Example \ref{example7}: $L^1$ convergence study for the Navier-Stokes variables $\rhoprime$ (left) and $\rho u_1$, $\rho u_2$ and
$\rhothetaprime$ (right) in the stochastic space computed at time $t=10s$ using the fully random model with the constant time step
$\Delta t=0.01$.\label{stoch_hermite_full_NS_M_L1_T10}}
\end{figure}
\begin{figure}[ht!]
\centerline{\includegraphics[scale=0.80]{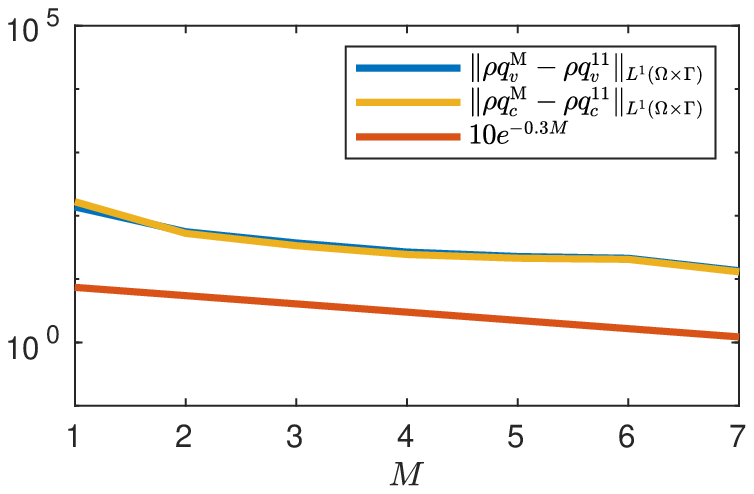}\hspace*{1.0cm}
\hspace*{0.1cm}\includegraphics[scale=0.80]{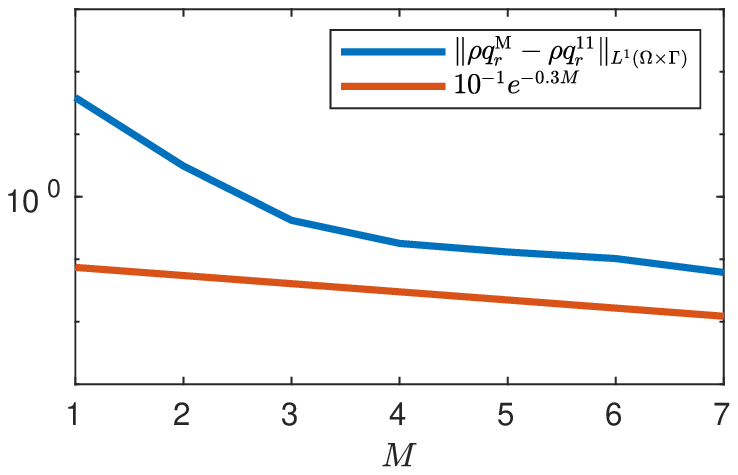}}
\caption{Example \ref{example7}: $L^1$ convergence study for the cloud variables $\rho q_v$ and $q_c$ (left) and $q_r$ (right) in the
stochastic space computed at time $t=10s$ using the fully random model with the constant time step $\Delta t=0.01$.
\label{stoch_hermite_full_cloud_M_L1_T10}}
\end{figure}
\end{example}

\newpage
$~$

\newpage
\begin{example}[\textbf{Rayleigh-B\'{e}nard convection in 2-D}]\label{RB_stoch_example2}
\phantom{.}

We consider a 2-D stochastic Rayleigh-B\'{e}nard convection simulated on a domain $\Omega=[0,5000]\times[0,1000]\,m^2$ that has been
discretized using a uniform $160\times160$ mesh. The initial data for the cloud variables are
\begin{equation*}
\begin{aligned}
(\widehat{q_v})_0(\x,0)&=0.025(\thetaprime(\x,0))_+,~~(\widehat{q_v})_1(\x,0)=0.1(\widehat{q_v})_0(\x,0),~~(\widehat{q_v})_k(\x,0)=0~\,
\mbox{for}\;2\le k\le M,\\
(\widehat{q_c})_0(\x,0)&=10^{-4}(\thetaprime(\x,0))_+,~~(\widehat{q_c})_k(\x,0)=0~\,\mbox{for}\;1\le k\le M,\\
(\widehat{q_r})_0(\x,0)&=10^{-6}(\thetaprime(\x,0))_+,~~(\widehat{q_r})_k(\x,0)=0~\,\mbox{for}\;1\le k\le M,
\end{aligned}
\end{equation*}
and for the Navier-Stokes variables are
\begin{equation*}
\begin{aligned}
&(\widehat{\rhoprime})_0(\x,0)=-\rhobar(\x)\frac{(\widehat{\thetaprime})_0(\x,0)}{\thetabar(\x)+(\widehat{\thetaprime})_0(\x,0)},~~
(\widehat{\rhoprime})_k(\x,0)=0~\,\mbox{for}\;1\le k\le M,\\
&(\widehat{\rho u_1})_0(\x,0)=0.001\big[(\widehat{\rhoprime})_0(\x,0)+\rhobar(\x)\big],~~(\widehat{\rho u_1})_k(\x,0)=0~\,
\mbox{for}\;1\le k\le M,\\
&(\widehat{\rho u_2})_0(\x,0)=\sin\left(\frac{\pi x_2}{500}\right)\big[(\widehat{\rhoprime})_0(\x,0)+\rhobar(\x)\big],~~
(\widehat{\rho u_2})_k(\x,0)=0~\,\mbox{for}\;1\le k\le M,\\
&(\widehat{\rhothetaprime})_0(\x,0)=\bar\rho(\x)(\widehat{\thetaprime})_0(\x,0)+\bar\theta(\widehat{\rhoprime})_0(\x,0)+
(\widehat{\thetaprime})_0(\x,0)(\widehat{\rhoprime})_0(\x,0),~~(\widehat{\rhothetaprime})_k(\x,0)=0~\,\mbox{for}\;1\le k\le M,
\end{aligned}
\end{equation*}
where
\begin{equation}
\begin{aligned}
&(\widehat{\thetaprime})_0(\x,0)=0.6\sin\left(\frac{\pi x_2}{500}\right),~~(\widehat{\thetaprime})_k(\x,0)=0~\,\mbox{for}\;1\le k\le M,\\
&\thetabar(\x)=284-\frac{x_2}{1000},\quad\rhobar(\x)=\frac{p_0}{R\thetabar(\x)}\pi_e(\x)^\frac{1}{\gamma-1},\quad
\pi_e(\x)=1-\frac{gx_2}{c_p\thetabar(\x)}.
\end{aligned}
\label{IC_stoch_NS2}
\end{equation}

\newpage
$~$

\newpage
We implement the following Dirichlet boundary conditions for the potential temperature:
\begin{equation*}
\theta(x_2=0)=284\,\mathrm{K}\quad\mbox{and}\quad\theta(x_2=1000)=283\,\mathrm{K},
\end{equation*}
as well as the periodic boundary conditions for all of the variables in the horizontal direction, no-slip boundary conditions for the
velocities at the vertical boundaries, and zero Neumann conditions for the remaining variables in the vertical direction, that is,
$\nabla\rhoprime\cdot\bm n=0$. Also, these boundary conditions have to be projected onto the stochastic space. The projections of the
periodic, no-slip and Neumann boundary conditions are straightforward and lead to the same conditions as in the deterministic case for all
of the expansion coefficients of the respective variable. Here, we briefly explain how the projection of the Dirichlet boundary conditions
works. We implement the Dirichlet boundary conditions for the potential temperature using
$\rho\theta(\x,t,\omega)=\rhothetaprime(\x,t,\omega)+\rhothetabar(\x)$. Rearranging and inserting the expansion for $\rhoprime(\x,t,\omega)$
and $\rhothetaprime(\x,t,\omega)$ gives
\begin{equation*}
\begin{aligned}
&\rhothetaprime(\x,t,\omega)-\rhoprime(\x,t,\omega)\theta(\x,t,\omega)=\rhobar(\x)\theta(\x,t,\omega)-\rhothetabar(\x)\\
\Longleftrightarrow&\sum_{k=0}^M(\widehat{\rhothetaprime})_k(\x,t)\Phi_k(\omega)-
\bigg(\sum_{k=0}^M(\widehat{\rhoprime})_k(\x,t)\Phi_k(\omega)\bigg)\theta(\x,t,\omega)=\rhobar(\x)\theta(\x,t,\omega)-\rhothetabar(\x).
\end{aligned}
\end{equation*}
At $\theta$ is constant at the boundary, applying the projection leads to
\begin{equation*}
\begin{aligned}
&(\widehat{\rhothetaprime})_0(x_2=0,t)-(\widehat{\rhoprime})_0(x_2=0,t)\theta(x_2=0)=\rhobar(x_2=0)\theta(x_2=0)-\rhothetabar(x_2=0),\\
&(\widehat{\rhothetaprime})_k(x_2=0,t)-(\widehat{\rhoprime})_k(x_2=0,t)\theta(x_2=0)=0~\,\mbox{for}\;1\le k\le M,
\end{aligned}
\end{equation*}
and analogously for $x_2=1000$.

In Figures \ref{RB_stoch_full_theta} and \ref{RB_stoch_full_qc}, we present snapshots of the expected values and standard deviations of the
potential temperature and the cloud variables at times $t=1000$ and $6000s$, respectively. Additionally, in Figure \ref{RB_stoch_full_qv1},
we plot the the differences between the expected value of the water vapor concentration and the saturation mixing ratio $E[q_v]-q_*$ at the
same times. As one can observe, the potential temperature exhibits a strong vertical gradient at time $t=1000s$. Similarly to the
deterministic and semi-random cases, at a later time $t=6000s$, supersaturated regions are formed in the rolls where the convection takes
place leading to the overall roll-like cloud flow structure.
\begin{figure}[ht!]
\centerline{\includegraphics[trim=0.9cm 0.1cm 0.4cm 0.0cm,clip,width=11.0cm]{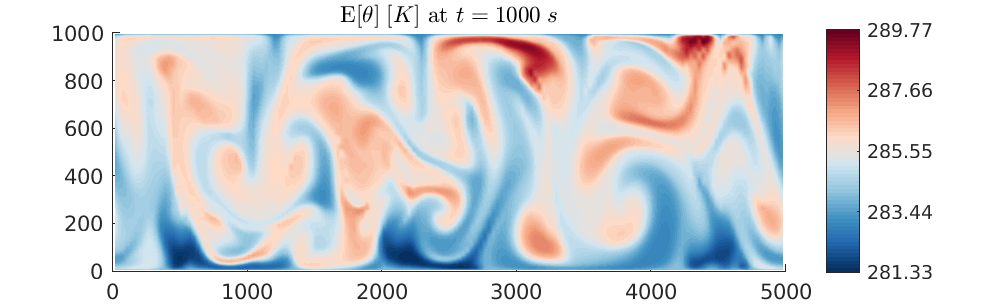}}
\vskip6pt
\centerline{\includegraphics[trim=0.9cm 0.1cm 0.4cm 0.0cm,clip,width=11.0cm]{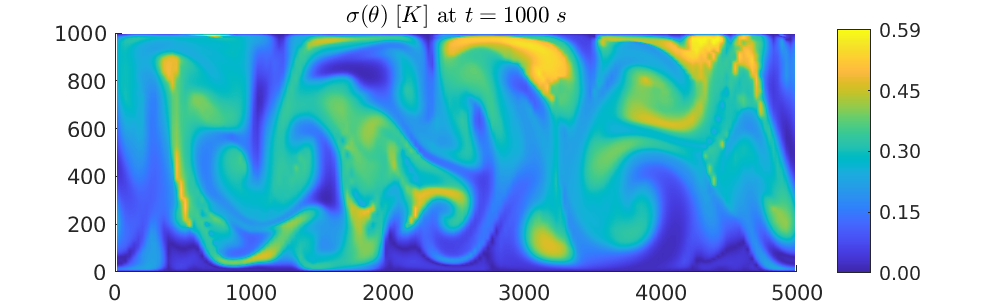}}
\vskip8pt
\centerline{\includegraphics[trim=0.9cm 0.1cm 0.4cm 0.0cm,clip,width=11.0cm]{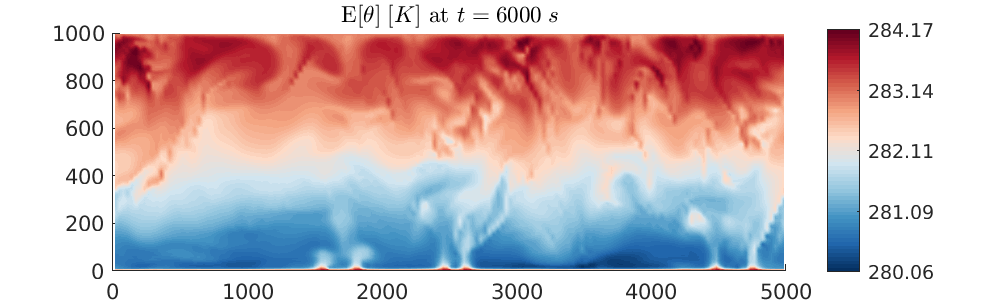}}
\vskip6pt
\centerline{\includegraphics[trim=0.9cm 0.1cm 0.4cm 0.0cm,clip,width=11.0cm]{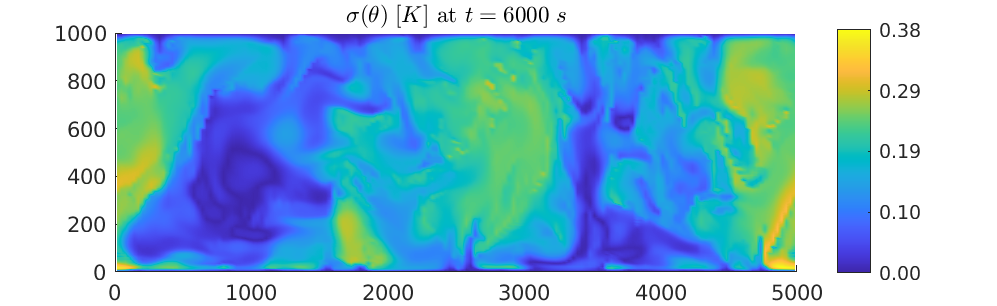}}
\caption{Example \ref{RB_stoch_example2}: Expected value and standard deviation of the potential temperature $\theta$ at $t=1000s$ and
$6000s$.\label{RB_stoch_full_theta}}
\end{figure}
\begin{figure}[ht!]
\centerline{\includegraphics[trim=0.9cm 0.1cm 0.4cm 0.0cm,clip,width=11.0cm]{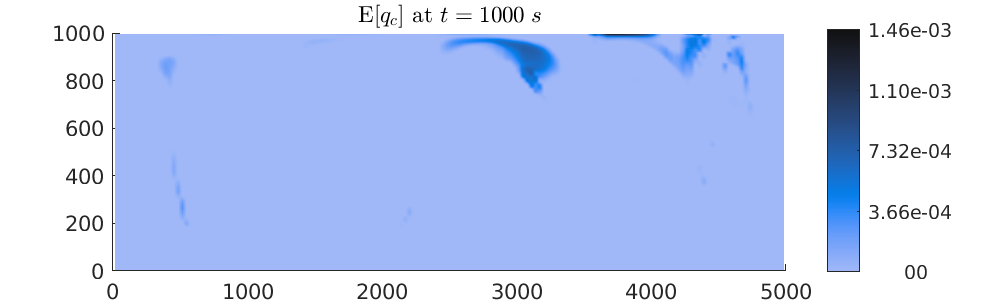}}
\vskip6pt
\centerline{\includegraphics[trim=0.9cm 0.1cm 0.4cm 0.0cm,clip,width=11.0cm]{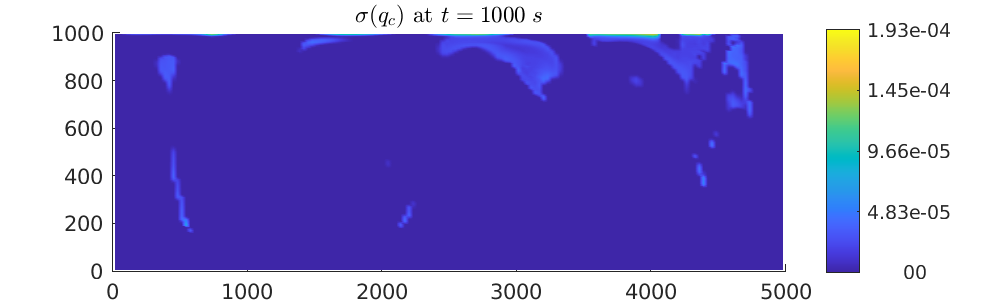}}
\vskip8pt
\centerline{\includegraphics[trim=0.9cm 0.1cm 0.4cm 0.0cm,clip,width=11.0cm]{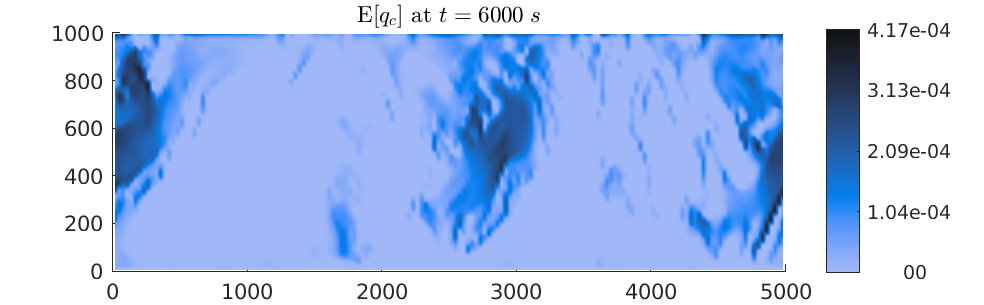}}
\vskip6pt
\centerline{\includegraphics[trim=0.9cm 0.1cm 0.4cm 0.0cm,clip,width=11.0cm]{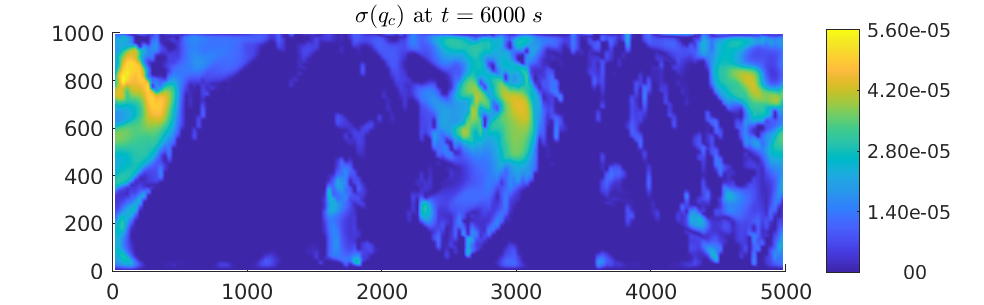}}
\caption{Example \ref{RB_stoch_example2}: Expected value and standard deviation of the cloud drops concentration $q_c$ at $t=1000s$ and
$6000s$.\label{RB_stoch_full_qc}}
\end{figure}
\begin{figure}[ht!]
\centerline{\includegraphics[trim=0.9cm 0.1cm 0.4cm 0.0cm,clip,width=11.0cm]{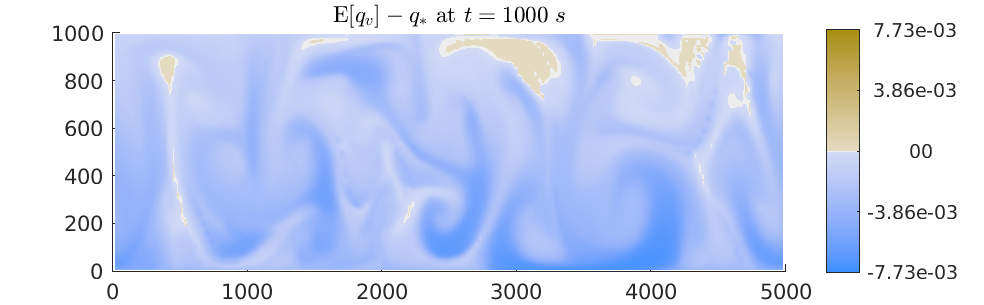}}
\vskip6pt
\centerline{\includegraphics[trim=0.9cm 0.1cm 0.4cm 0.0cm,clip,width=11.0cm]{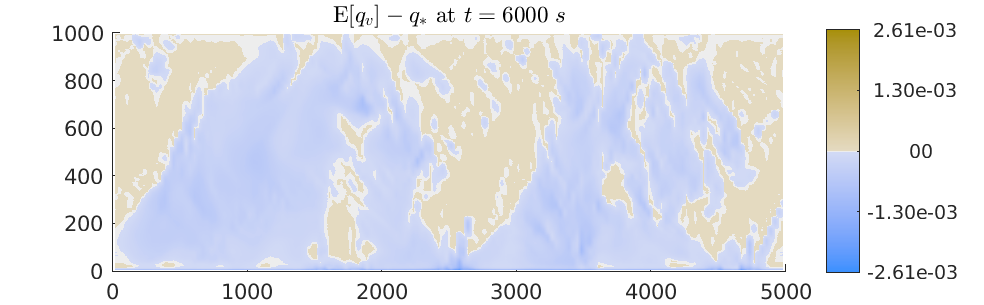}}
\caption{Example \ref{RB_stoch_example2}: Difference between the expected value of the water vapor concentration and the saturation mixing
ratio ($\,\mathbb E[q_v]-q_*$) at $t=1000s$ and $6000s$.\label{RB_stoch_full_qv1}}
\end{figure}

In Figures \ref{2D_RB_stoch_full_E_and_sigma_theta_comp} and \ref{2D_RB_stoch_full_E_and_sigma_comp}, we present the time evolution of the
mean expected value per $m^2$ as well as the mean standard deviation per $m^2$ for the potential temperature and the cloud variables. In $d$
space dimensions these quantities can be computed for uniformly distributed perturbations in the following way:
\begin{equation*}
\begin{aligned}
\mathbb{E}\Bigg[\frac{h^d}{|\Omega|}\sum_{i=1}^{N^d}(q_\ell)_i\Bigg]=\frac{h^d}{|\Omega|}\sum_{i=1}^{N^d}\mathbb{E}\left[(q_\ell)_i\right]=
\frac{h^d}{|\Omega|}\sum_{i=1}^{N^d}(\widehat{(q_\ell)}_0)_i,~~\sigma\Bigg(\frac{h^d}{|\Omega|}\sum_{i=1}^{N^d}(q_\ell)_i\Bigg)=
\frac{h^d}{|\Omega|}\sqrt{\sum_{k=1}^M\frac{1}{2k+1}\left(\sum_{i=1}^{N^d}(\widehat{(q_\ell)}_k)_i\right)^2},
\end{aligned}
\end{equation*}
where $N^d$ is the number of mesh cells and $\ell\in\{v,c,r\}$. We compare the solutions using 0\% (purely deterministic model) and 10\% of
perturbation of the initial data in $q_v$, where for 10\% of perturbation the solutions are added in both fully- and semi-random
Navier-Stokes-cloud models. The time evolution of the averaged quantities clearly shows the differences between the semi-random and
fully random models. In all shown cases (including the purely deterministic one), the time evolution starts with cloud formation and thus
increase of cloud water on the expense of water vapor and also latent heat release (increase of $\theta$). However, for the semi-random
model the rain formation starts earlier than in the deterministic and fully random simulations. Since rain is falling into subsaturated
regions which induces evaporation, this leads to a different time evolution in all variables. Generally, we observe a much stronger cooling
effect of the system due to evaporation of rain in the semi-random model than in the other simulations. This is probably due to the use of
the expected values of the terms for phase changes in the energy equation. Although the general qualitative behavior in the time evolution
of the expected values is quite similar, the absolute values differ quite substantially. The same is true for the standard deviations of the
cloud variables $q_v$ $q_c$ and $q_r$ as shown in the right column of Figure \ref{2D_RB_stoch_full_E_and_sigma_comp}; the variation in the
water variables (and also in $\theta$) is much larger for the fully random model than for the semi-random one. This is also reasonable,
since the fully random model can capture the correct feedback of the latent heat release in the phase changes for the ``different
realizations'', whereas the semi-random model only feeds back the averaged potential temperature, leading to a smaller variability. Some
examples of the expected values and the related standard deviations for $\theta$, $q_c$ and the super/subsaturation (in terms of
$\mathbb E[q_v]-q_*$) are shown in Figures \ref{RB_stoch_full_theta}--\ref{RB_stoch_full_qv1}.
\begin{figure}[ht!]
\centerline{\includegraphics[trim=0.8cm 0.1cm 1.2cm 0.3cm,clip,width=7.5cm]{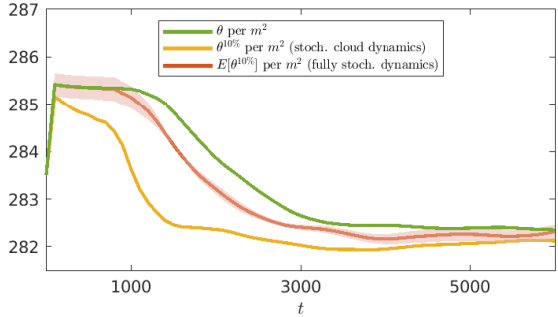}}
\caption {Example \ref{RB_stoch_example2}: Time evolution of the expected values with their standard deviations (shaded region) for the
potential temperature $\theta$ per $m^2$ using 0\% (purely deterministic case) and 10\% perturbation of the initial data in $q_v$, where the
latter was simulated using both fully- and semi-random models.\label{2D_RB_stoch_full_E_and_sigma_theta_comp}}
\end{figure}
\begin{figure}[ht!]
\centerline{\includegraphics[trim=0.3cm 0.0cm 0.7cm 0.0cm,clip,width=5.05cm]{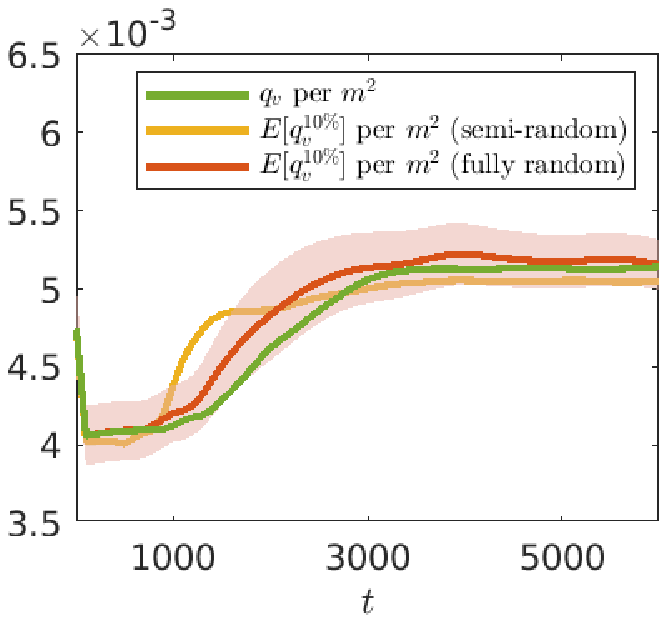}\hspace*{1.2cm}
\includegraphics[trim=0.0cm 0.0cm 0.0cm 0.0cm,clip,width=4.85cm]{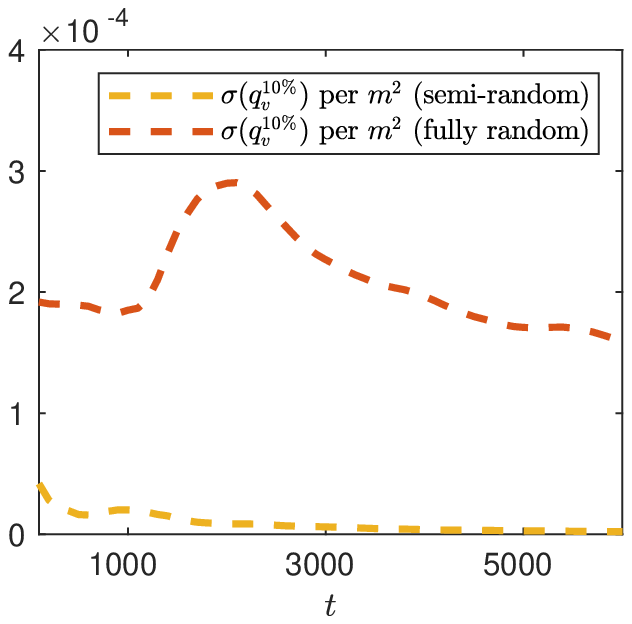}}
\vskip1pt
\centerline{\includegraphics[trim=0.3cm 0.0cm 0.7cm 0.0cm,clip,width=5.0cm]{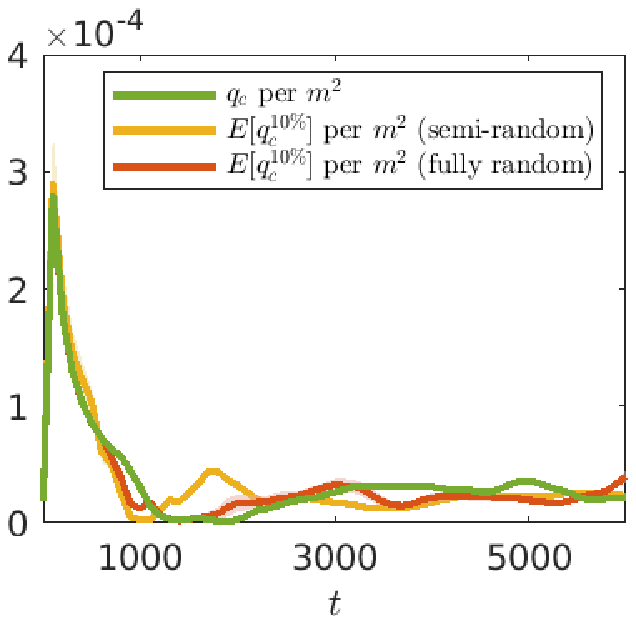}\hspace*{1.2cm}
\includegraphics[trim=0.0cm 0.0cm 0.0cm 0.0cm,clip,width=4.85cm]{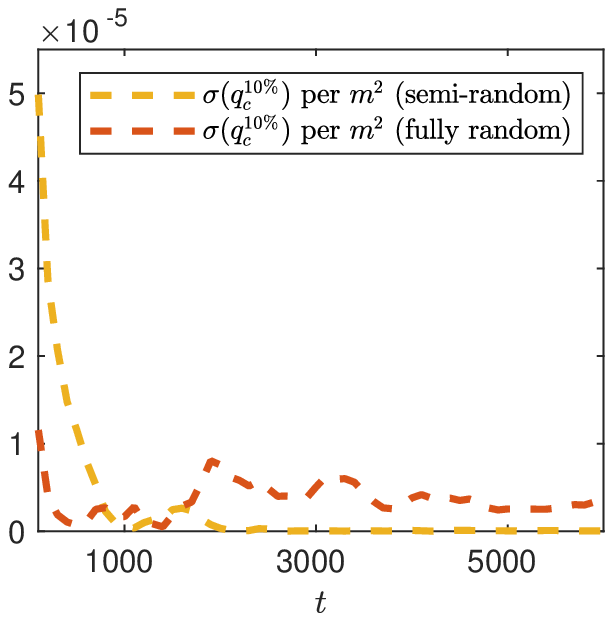}}
\vskip1pt
\centerline{\includegraphics[trim=0.3cm 0.0cm 0.7cm 0.0cm,clip,width=5.0cm]{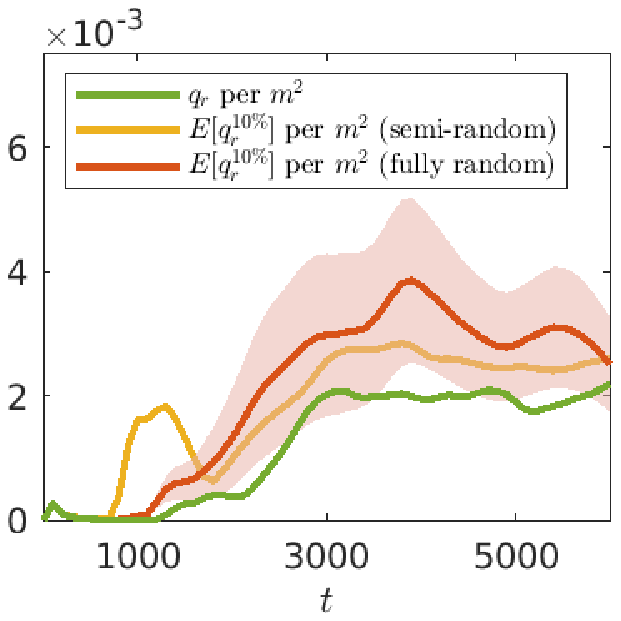}\hspace*{0.95cm}
\includegraphics[trim=0.0cm 0.0cm 0.0cm 0.0cm,clip,width=5.05cm]{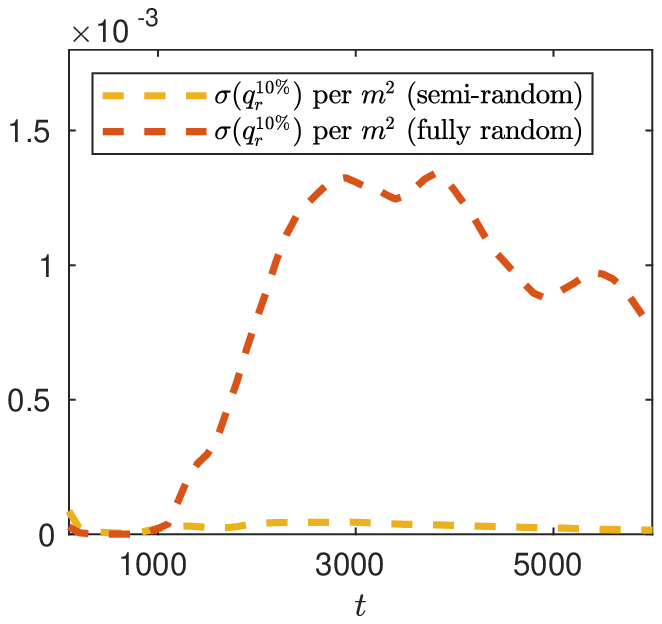}}
\caption {Example \ref{RB_stoch_example2}: Time evolution of the expected values with their standard deviations for the cloud variables per
$m^2$ (shaded region, left column) and standard deviation (right column) using 0\% (purely deterministic case) and 10\% perturbation of the
initial data in $q_v$, where the latter was simulated using both fully- and semi-random models.\label{2D_RB_stoch_full_E_and_sigma_comp}}
\end{figure}
\end{example}

\begin{example}[\textbf{Rayleigh-B\'{e}nard convection in 3-D}]\label{RB_stoch_example4}
\phantom{.}

In the final example, we consider a 3-D stochastic Rayleigh-B\'{e}nard convection. The initial data for the cloud variables are
\begin{equation*}
\begin{aligned}
(\widehat{q_v})_0(\x,0)&=0.025(\thetaprime(\x,0))_+,~~(\widehat{q_v})_1(\x,0)=\nu(\widehat{q_v})_0(\x,0),~~(\widehat{q_v})_k(\x,0)=0~\,
\mbox{for}\;2\le k\le M,\\
(\widehat{q_c})_0(\x,0)&=10^{-4}(\thetaprime(\x,0))_+,~~(\widehat{q_c})_k(\x,0)=0~\,\mbox{for}\;1\le k\le M,\\
(\widehat{q_r})_0(\x,0)&=10^{-6}(\thetaprime(\x,0))_+,~~(\widehat{q_r})_k(\x,0)=0~\,\mbox{for}\;1\le k\le M,
\end{aligned}
\end{equation*}
with $\nu=0$, 0.1, 0.2 or 0.5, which correspond to 0\% (pure deterministic case), 10\%, 20\% or 50\% perturbation of the initial water vapor
concentration. For the Navier-Stokes variables we take purely deterministic initial data, which in terms of their expansion coefficients
read as
\begin{equation*}
\begin{aligned}
&(\widehat{\rhoprime})_0(\x,0)=-\rhobar(\x)\frac{(\widehat{\thetaprime})_0(\x,0)}{\thetabar(\x)+(\widehat{\thetaprime})_0(\x,0)},\quad
(\widehat{\rho u_1})_0(\x,0)=0.001\big[(\widehat{\rhoprime})_0(\x,0)+\rhobar(\x)\big],\\
&(\widehat{\rho u_2})_0(\x,0)=0.001\big[(\widehat{\rhoprime})_0(\x,0)+\rhobar(\x)\big],\quad
(\widehat{\rho u_3})_0(\x,0)=\sin\left(\frac{\pi x_3}{500}\right)\big[(\widehat{\rhoprime})_0(\x,0)+\rhobar(\x)\big],\\
&(\widehat{\rhothetaprime})_0(\x,0)=\bar{\rho}(\x)(\widehat{\thetaprime})_0(\x,0)+\bar{\theta}(\widehat{\rhoprime})_0(\x,0)+
(\widehat{\thetaprime})_0(\x,0)(\widehat{\rhoprime})_0(\x,0),\\
&(\widehat{\rhoprime})_k(\x,0)=(\widehat{\rho u_1})_k(\x,0)=(\widehat{\rho u_2})_k(\x,0)=(\widehat{\rho u_3})_k(\x,0)=
(\widehat{\rhothetaprime})_k(\x,0)=0~\,\mbox{for}\;1\le k\le M,
\end{aligned}
\end{equation*}
where
\begin{equation*}
(\widehat{\thetaprime})_0(\x,0)=0.6\sin\left(\frac{\pi x_3}{500}\right),
\end{equation*}
and $\thetabar(\x)$ and $\rhobar(\x)$ are chosen as in \eqref{IC_stoch_NS2}. The solution is computed in the domain
$\Omega=[0,5000]\times[0,5000]\times[0,1000]\,m^3$ which is discretized using a uniform $50\times50\times50$ mesh.

In Figures \ref{3D_RB_stoch_full_theta}--\ref{3D_RB_stoch_full_qr}, we present the influence of the 10\%, 20\% and 50\% initial water vapor
perturbation on the expected values of the potential temperature, cloud droplets and rain drops concentration at times $t=1000s$ and
$6000s$. The influence on the supersaturated and subsaturated regions is highlighted as a 2-D slice along $x_1=3000$ in Figure
\ref{3D_RB_stoch_full_qstar_slice2}, where we depict the difference between the expected water vapor concentration and the saturation mixing
ratio. For a better comparison, we have used the same vertical scales in all of the plots. Here, one can clearly observe a different
behavior compared with the semi-random case. The vertical gradient of the potential temperature increases as the size of the initial
perturbations increases (see Figure \ref{3D_RB_stoch_full_theta}), while the pattern of the developed convection cells is similar for
different perturbations. The latent heat release increases the vertical motions in the convective cells, which leads to additional feedback,
such as stronger and more cloud formation (see Figure \ref{3D_RB_stoch_full_qc}), which in turn leads to the formation of a much larger
amount of rain water, especially at a later time $t=6000s$ (see Figure \ref{3D_RB_stoch_full_qr}). At the time $t=1000s$ one can see that
the roll-like structure of the clouds in the deterministic case (that is, with 0\% perturbation) again end up in a more cell-like structure
in the initially perturbed cases.
\begin{figure}[ht!]
\centerline{\includegraphics[trim=1.2cm 2.2cm 0.7cm 2.7cm,clip,width=6.0cm]{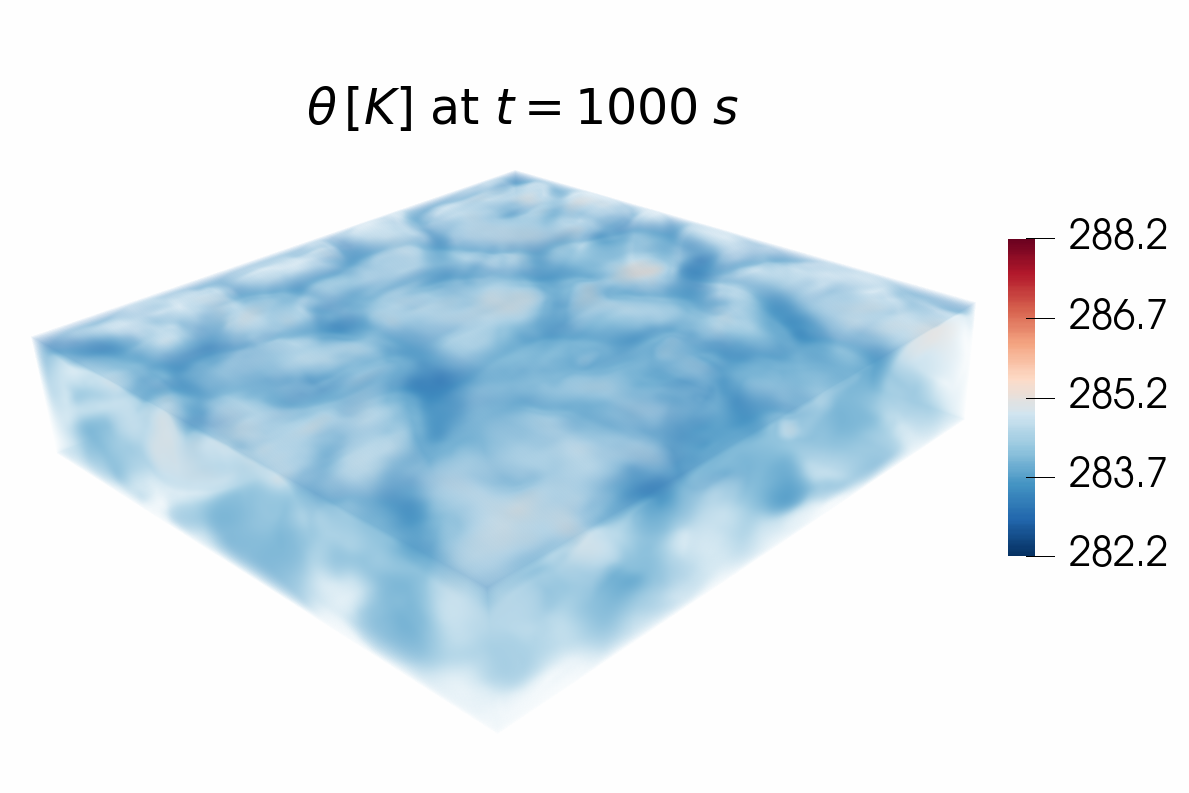}\hspace*{1.0cm}
\includegraphics[trim=1.2cm 2.2cm 0.7cm 2.7cm,clip,width=6.0cm]{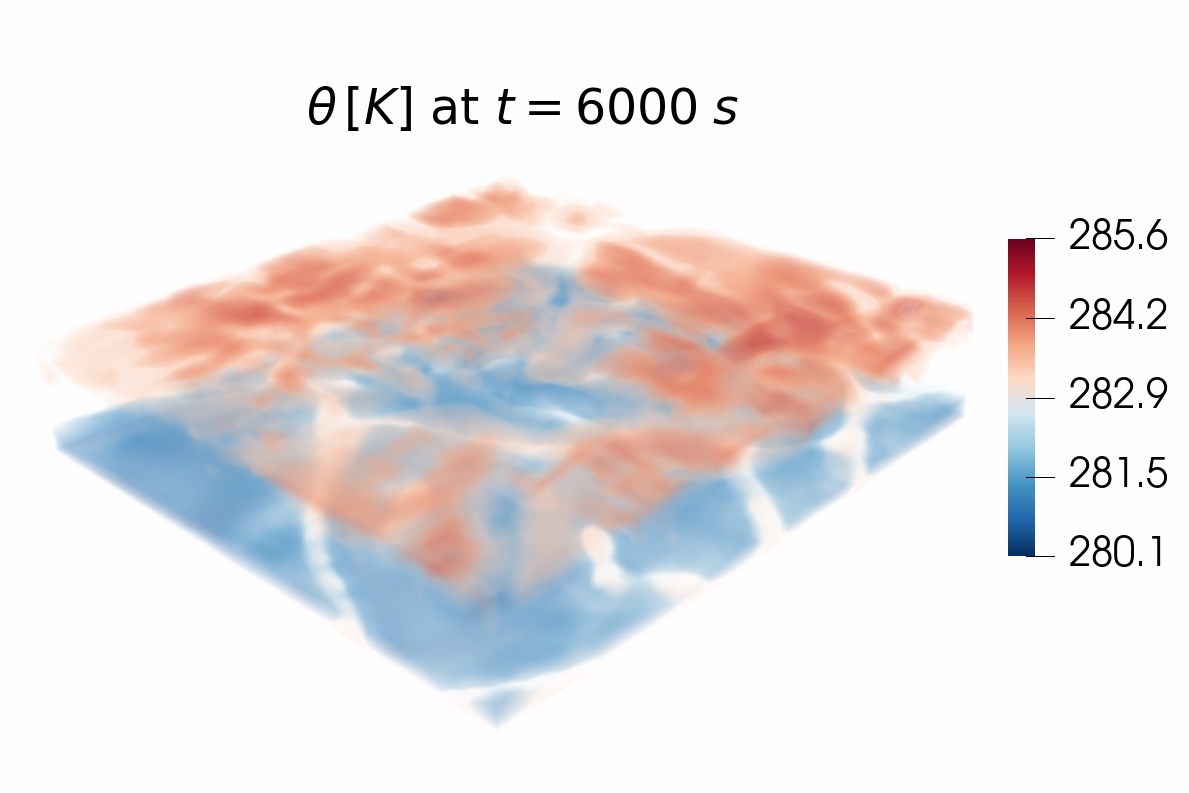}}
\vskip8pt
\centerline{\includegraphics[trim=1.2cm 2.2cm 0.7cm 2.7cm,clip,width=6.0cm]{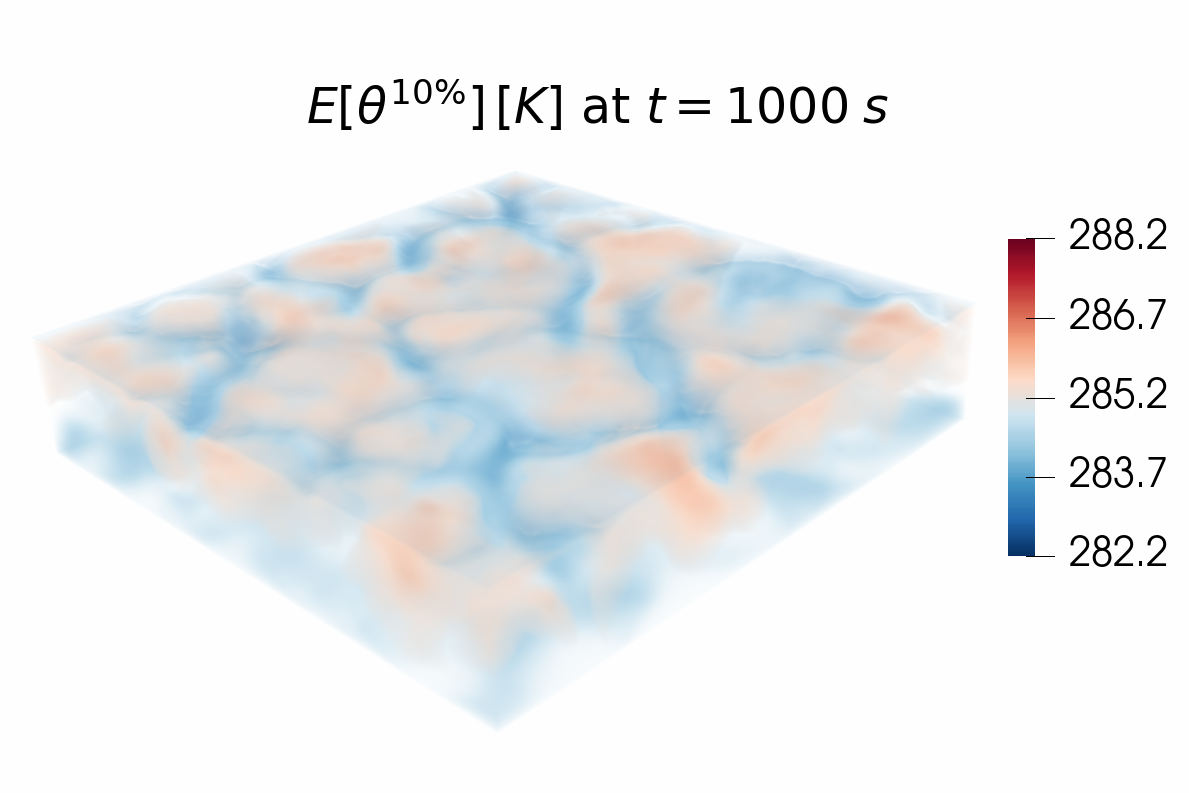}\hspace*{1.0cm}
\includegraphics[trim=1.2cm 2.2cm 0.7cm 2.7cm,clip,width=6.0cm]{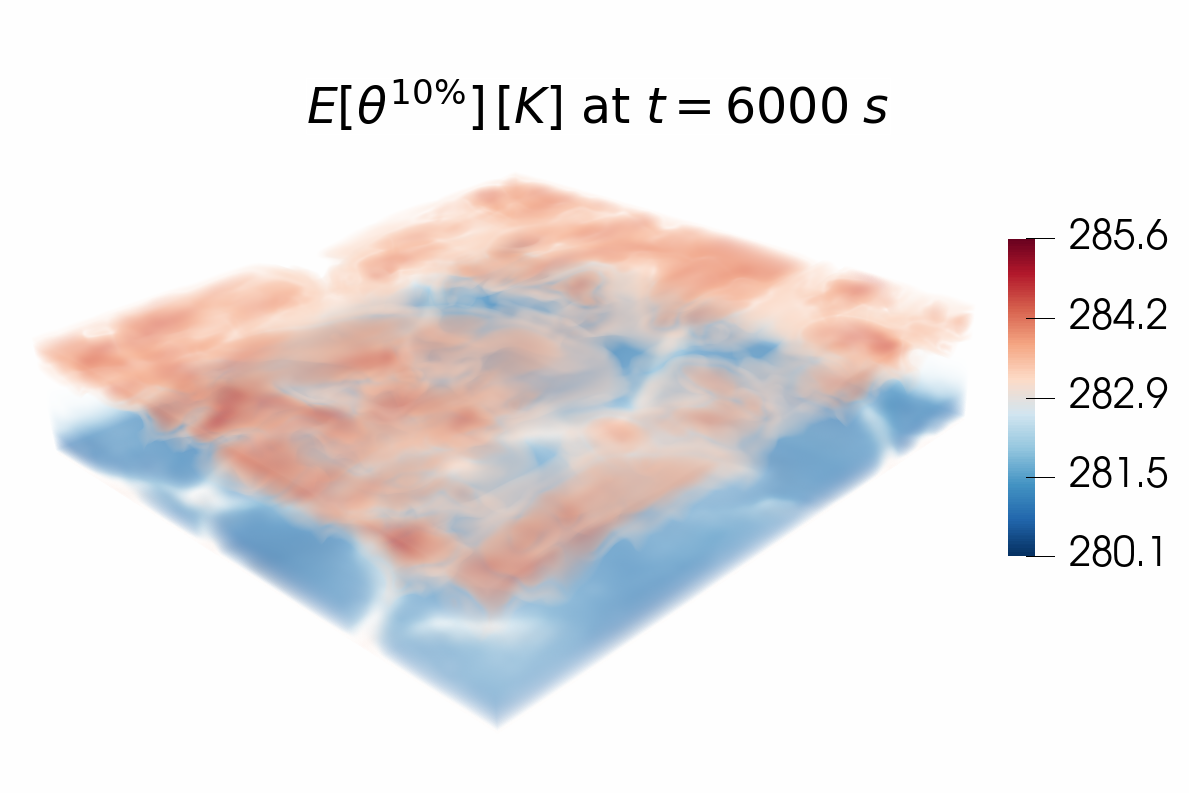}}
\vskip8pt
\centerline{\includegraphics[trim=1.2cm 2.2cm 0.7cm 2.7cm,clip,width=6.0cm]{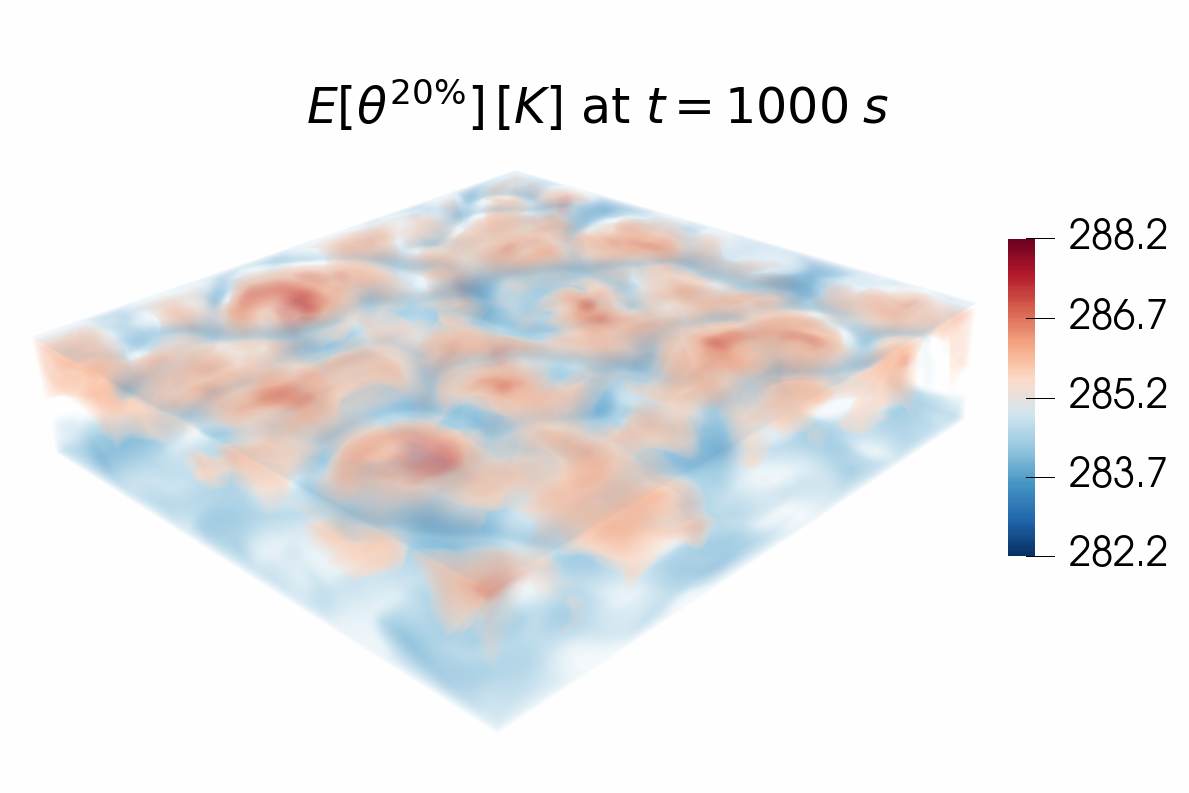}\hspace*{1.0cm}
\includegraphics[trim=1.2cm 2.2cm 0.7cm 2.7cm,clip,width=6.0cm]{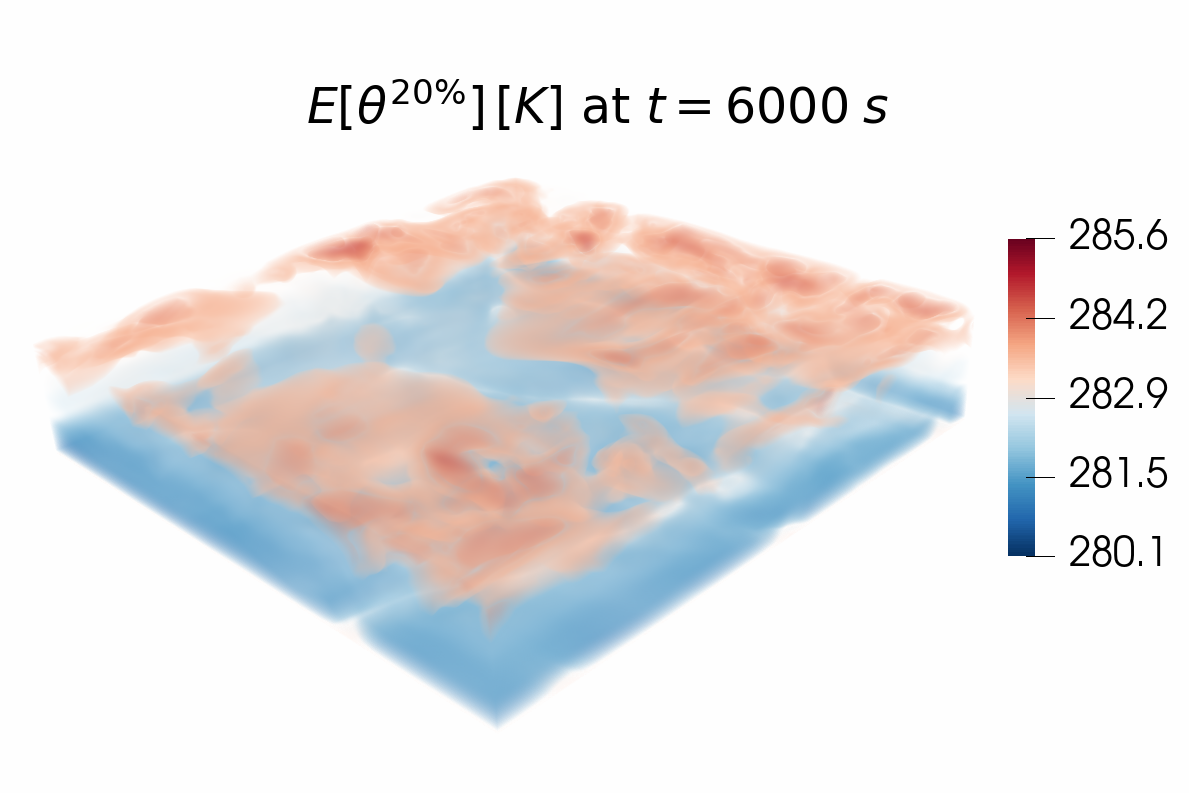}}
\vskip8pt
\centerline{\includegraphics[trim=1.2cm 2.2cm 0.7cm 2.7cm,clip,width=6.0cm]{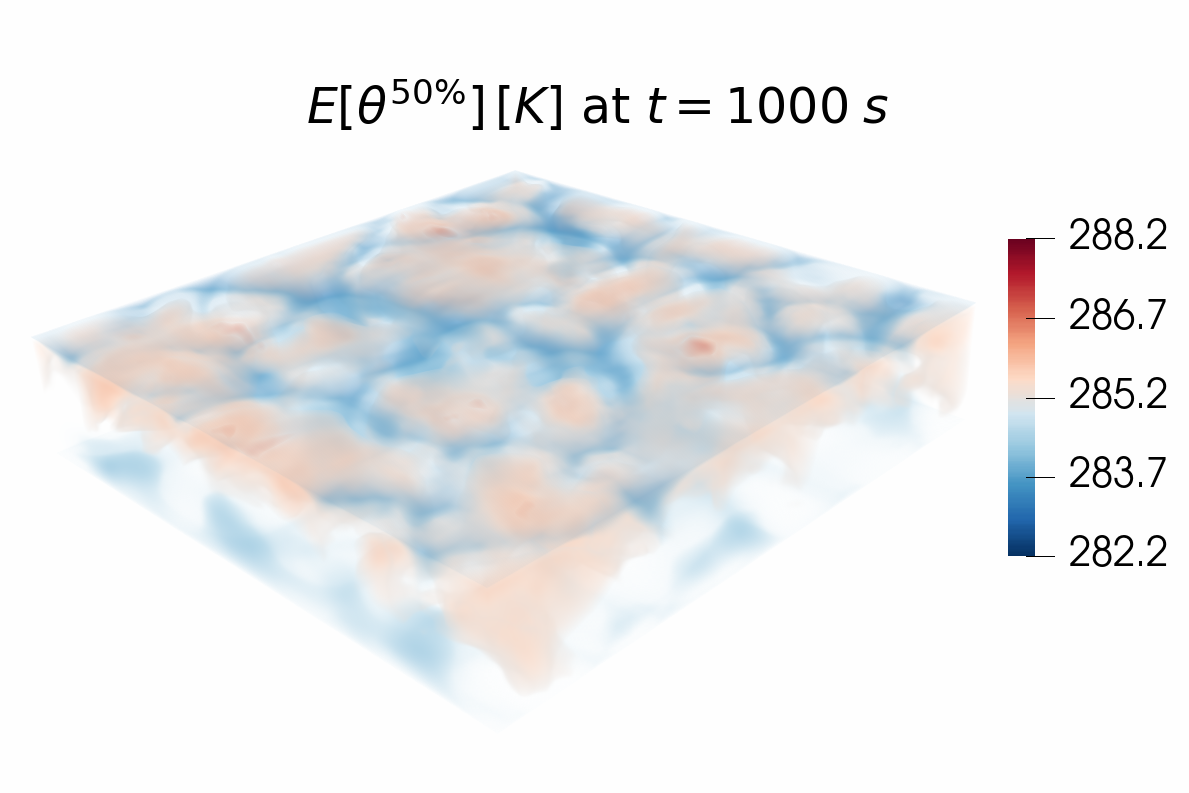}\hspace*{1.0cm}
\includegraphics[trim=1.2cm 2.2cm 0.7cm 2.7cm,clip,width=6.0cm]{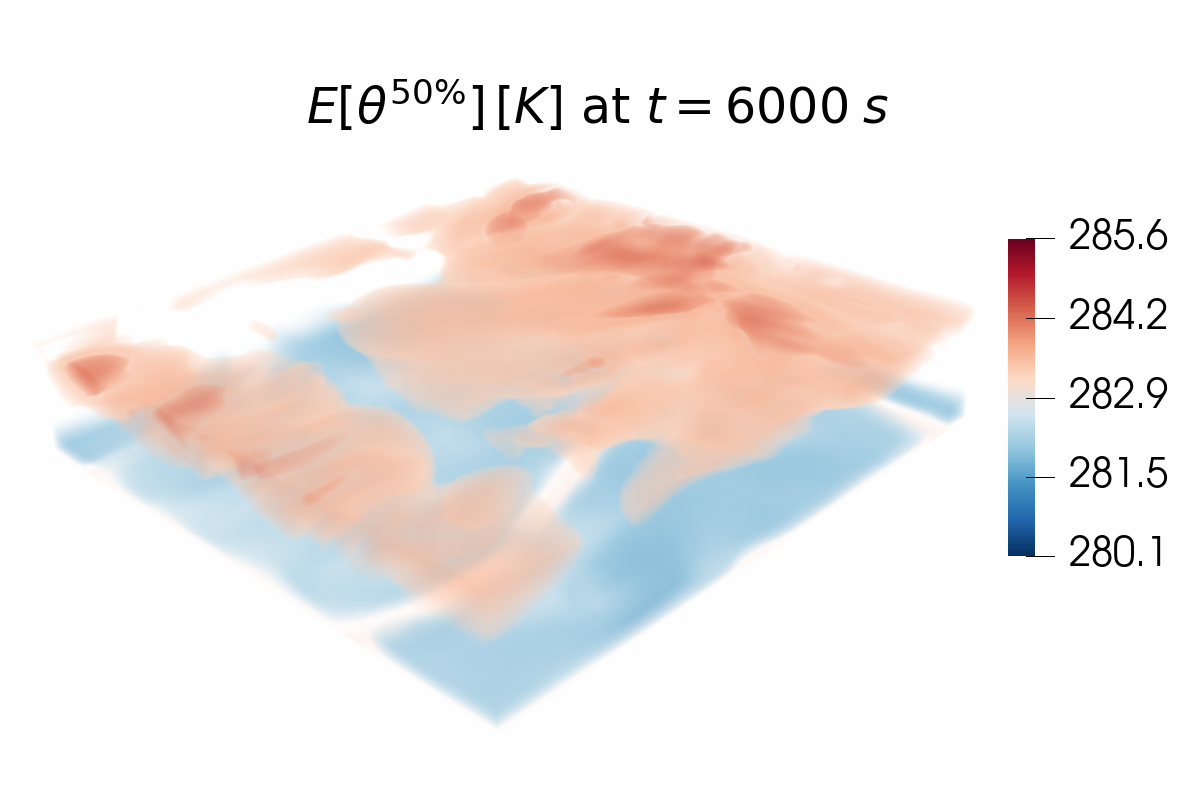}}
\caption{Example \ref{RB_stoch_example4}: Expected value of the potential temperature $\theta$ at times $t=1000s$ and $6000s$ with 0\%,
10\%, 20\% and 50\% perturbations of the initial water vapor concentration.\label{3D_RB_stoch_full_theta}}
\end{figure}
\begin{figure}[ht!]
\centerline{\includegraphics[trim=1.2cm 2.2cm 0.1cm 2.5cm,clip,width=6.0cm]{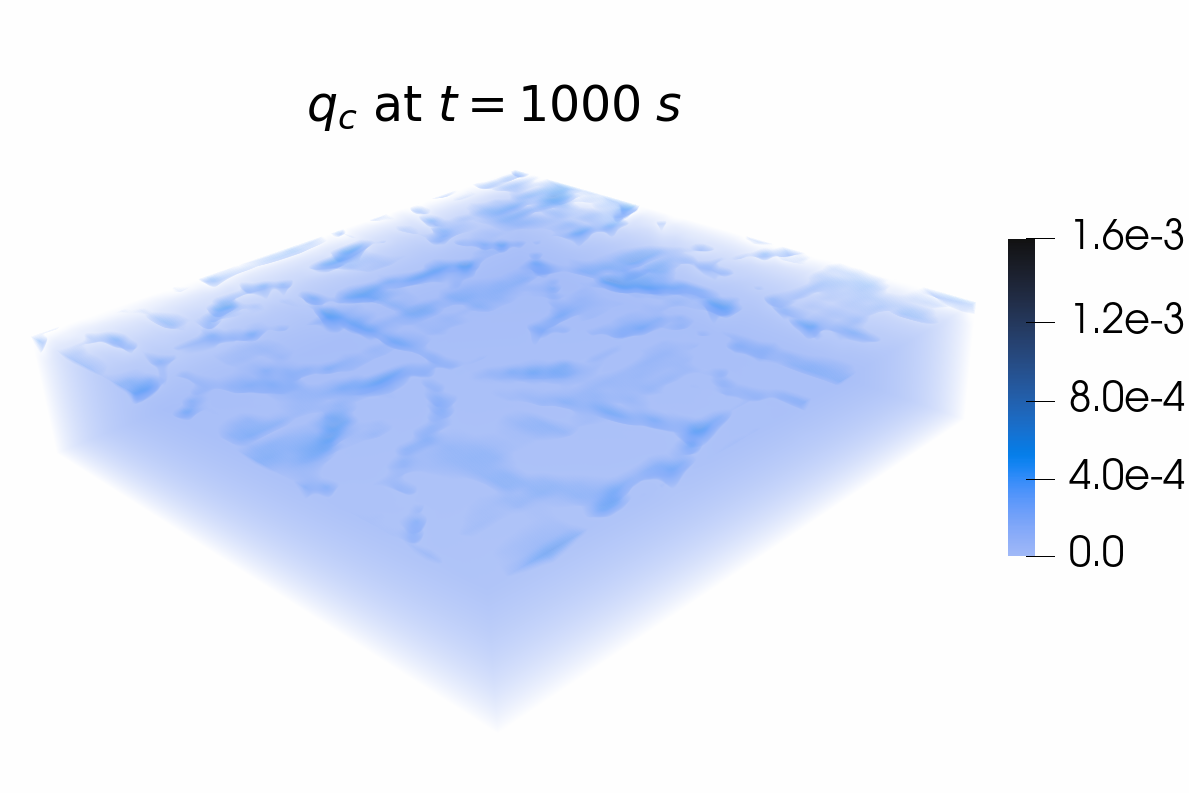}\hspace*{1.0cm}
\includegraphics[trim=1.2cm 2.2cm 0.1cm 2.5cm,clip,width=6.0cm]{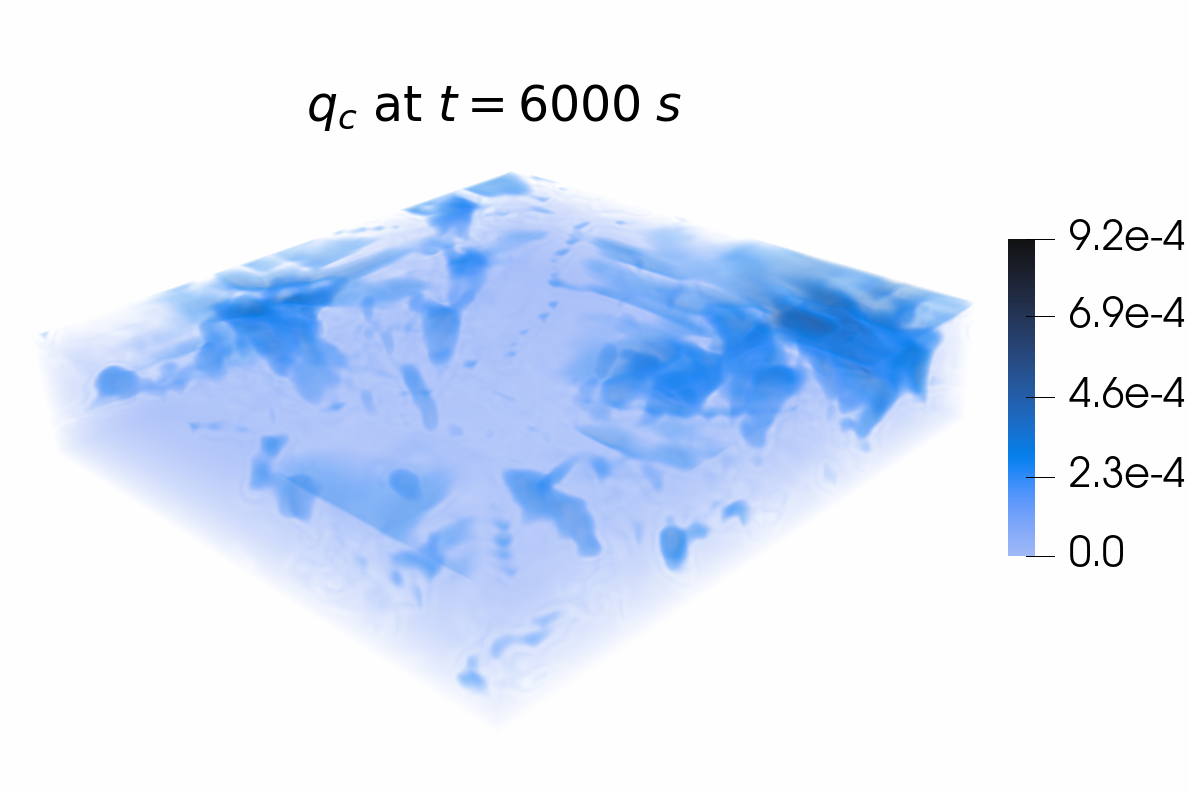}}
\vskip8pt
\centerline{\includegraphics[trim=1.2cm 2.2cm 0.1cm 2.5cm,clip,width=6.0cm]{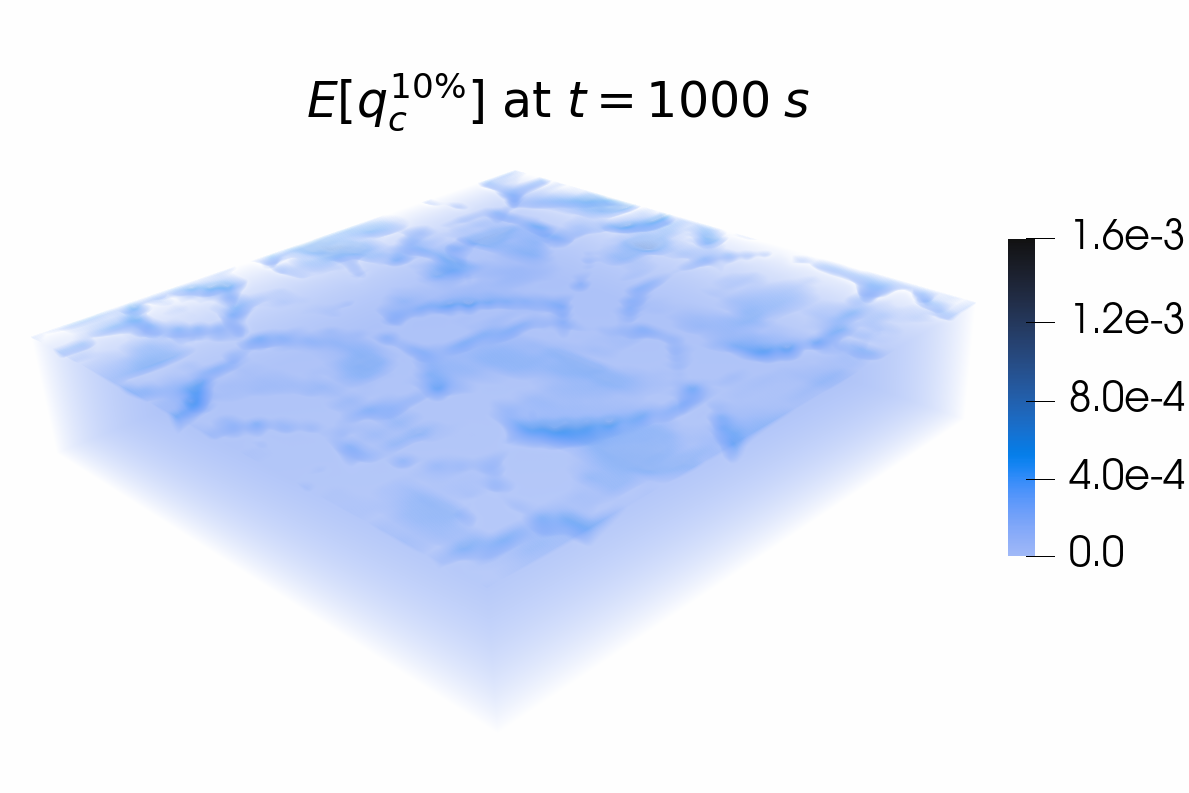}\hspace*{1.0cm}
\includegraphics[trim=1.2cm 2.2cm 0.1cm 2.5cm,clip,width=6.0cm]{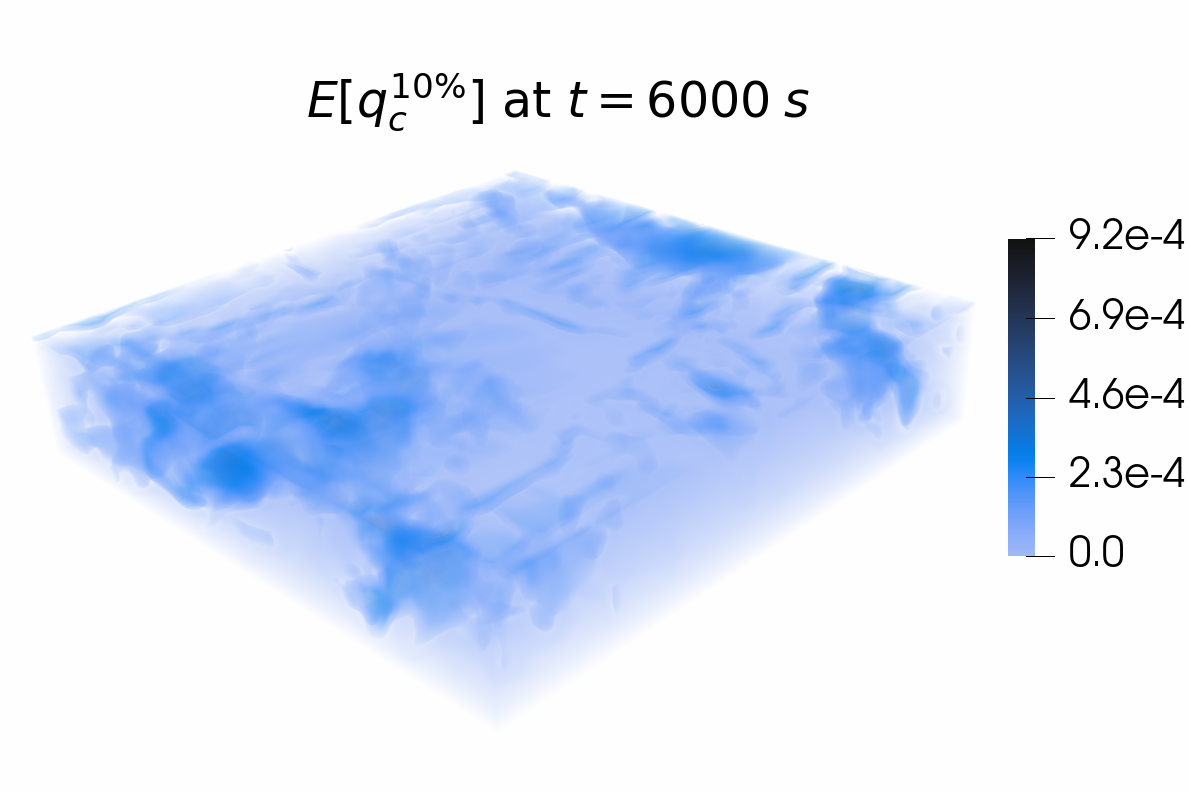}}
\vskip8pt
\centerline{\includegraphics[trim=1.2cm 2.2cm 0.1cm 2.5cm,clip,width=6.0cm]{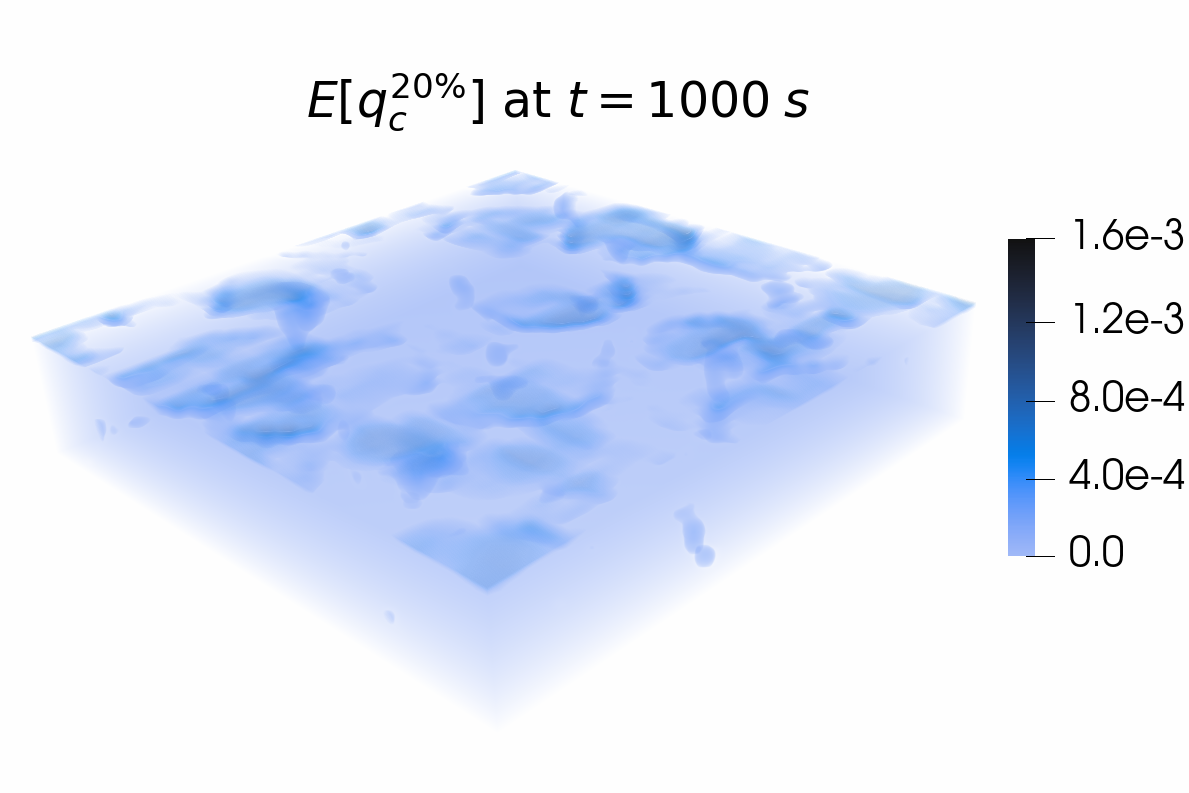}\hspace*{1.0cm}
\includegraphics[trim=1.2cm 2.2cm 0.1cm 2.5cm,clip,width=6.0cm]{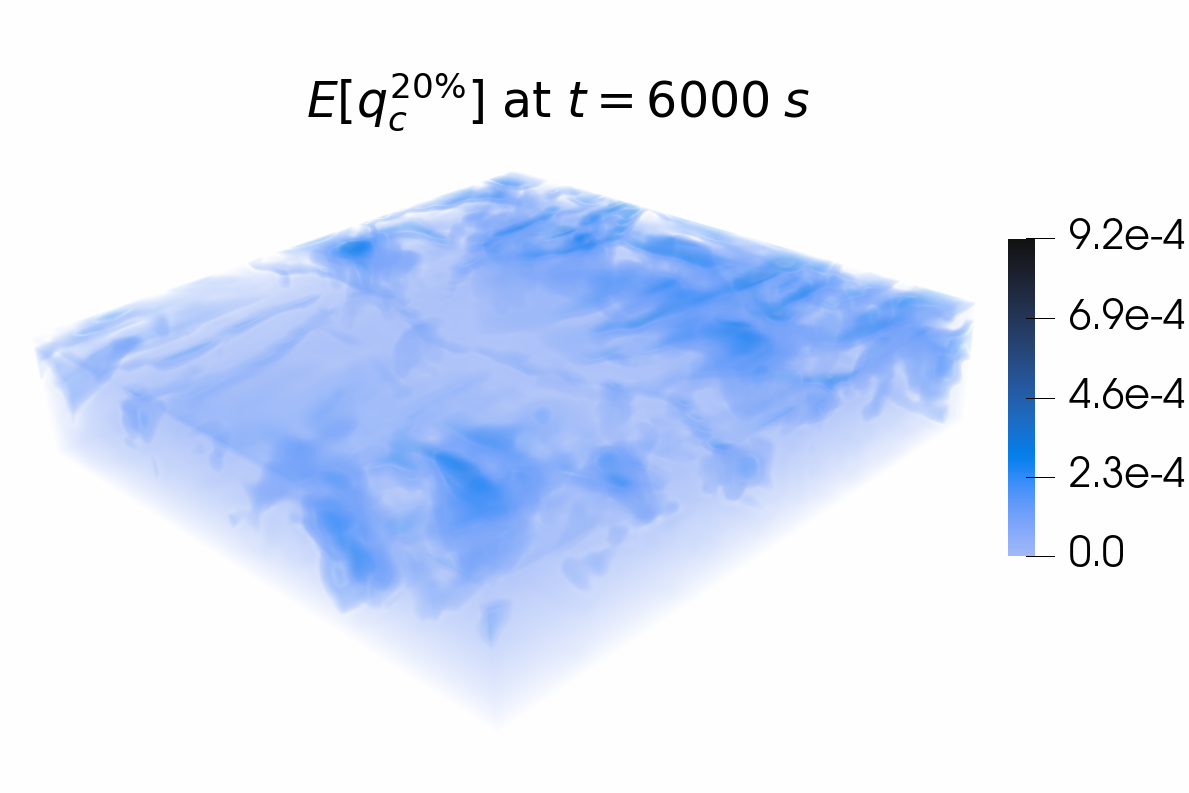}}
\vskip8pt
\centerline{\includegraphics[trim=1.2cm 2.2cm 0.1cm 2.5cm,clip,width=6.0cm]{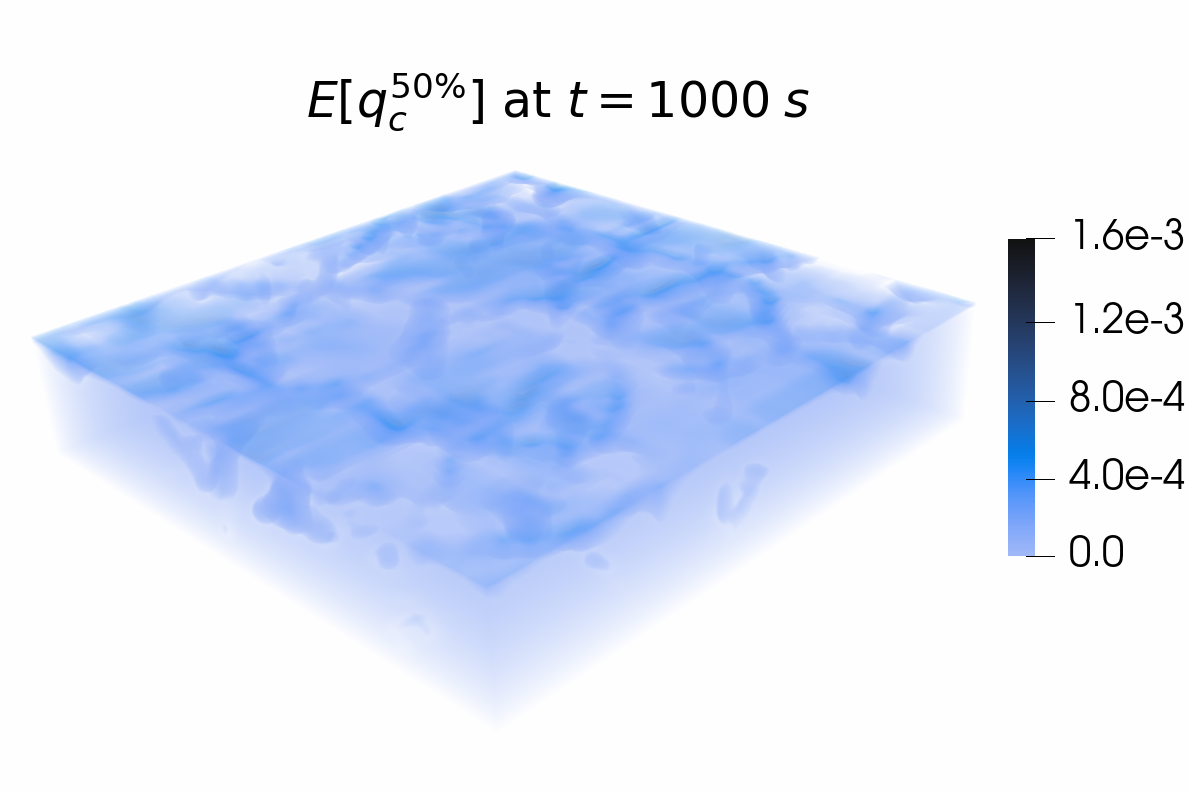}\hspace*{1.0cm}
\includegraphics[trim=1.2cm 2.2cm 0.1cm 2.5cm,clip,width=6.0cm]{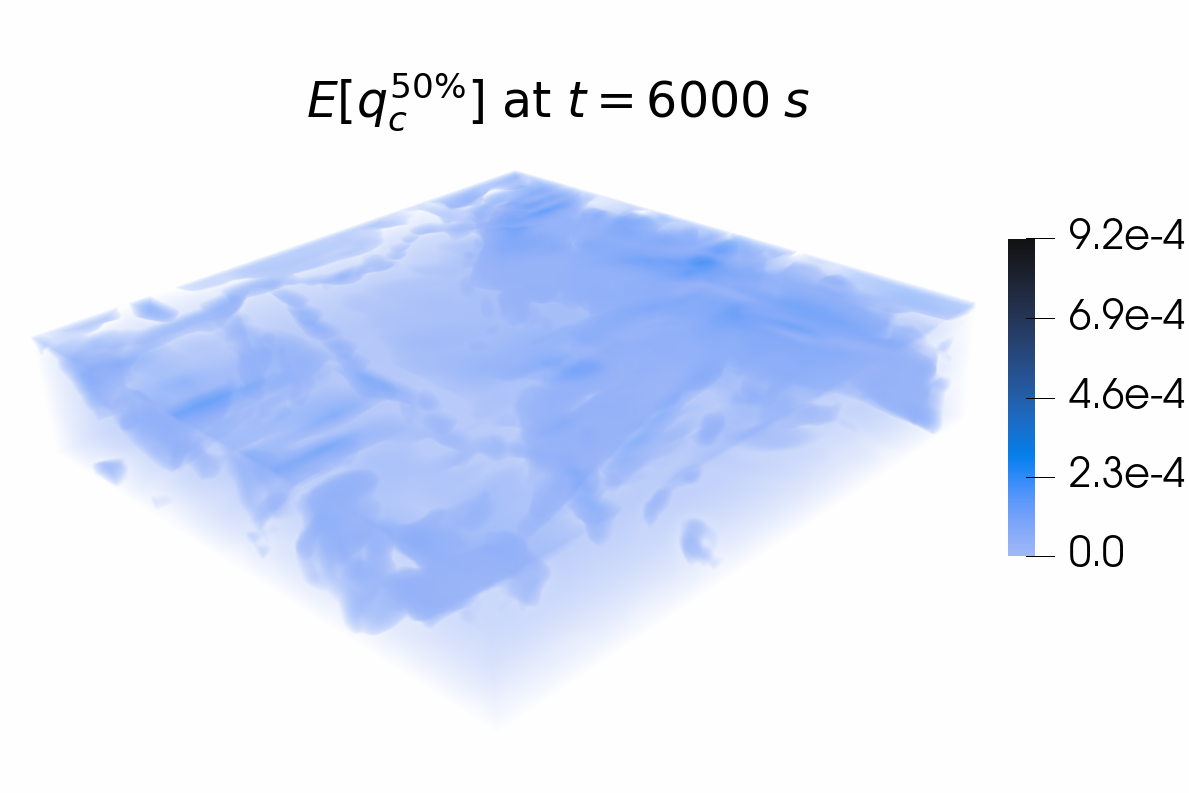}}
\caption{Example \ref{RB_stoch_example4}: Expected value of the cloud drops concentration $q_c$ at times $t=1000s$ and $6000s$ with 0\%,
10\%, 20\% and 50\% perturbations of the initial water vapor concentration.\label{3D_RB_stoch_full_qc}}
\end{figure}
\begin{figure}[ht!]
\centerline{\includegraphics[trim=1.2cm 2.2cm 0.1cm 2.5cm,clip,width=6.0cm]{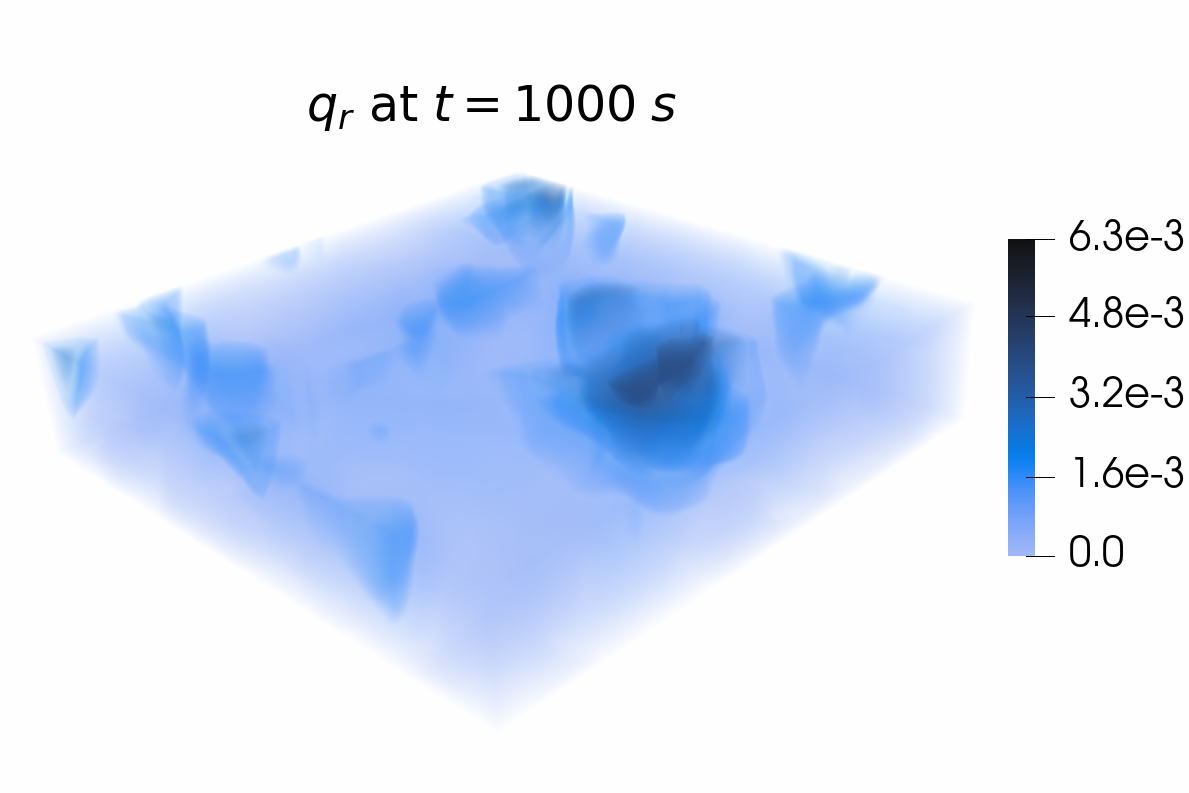}\hspace*{1.0cm}
\includegraphics[trim=1.2cm 2.2cm 0.1cm 2.5cm,clip,width=6.0cm]{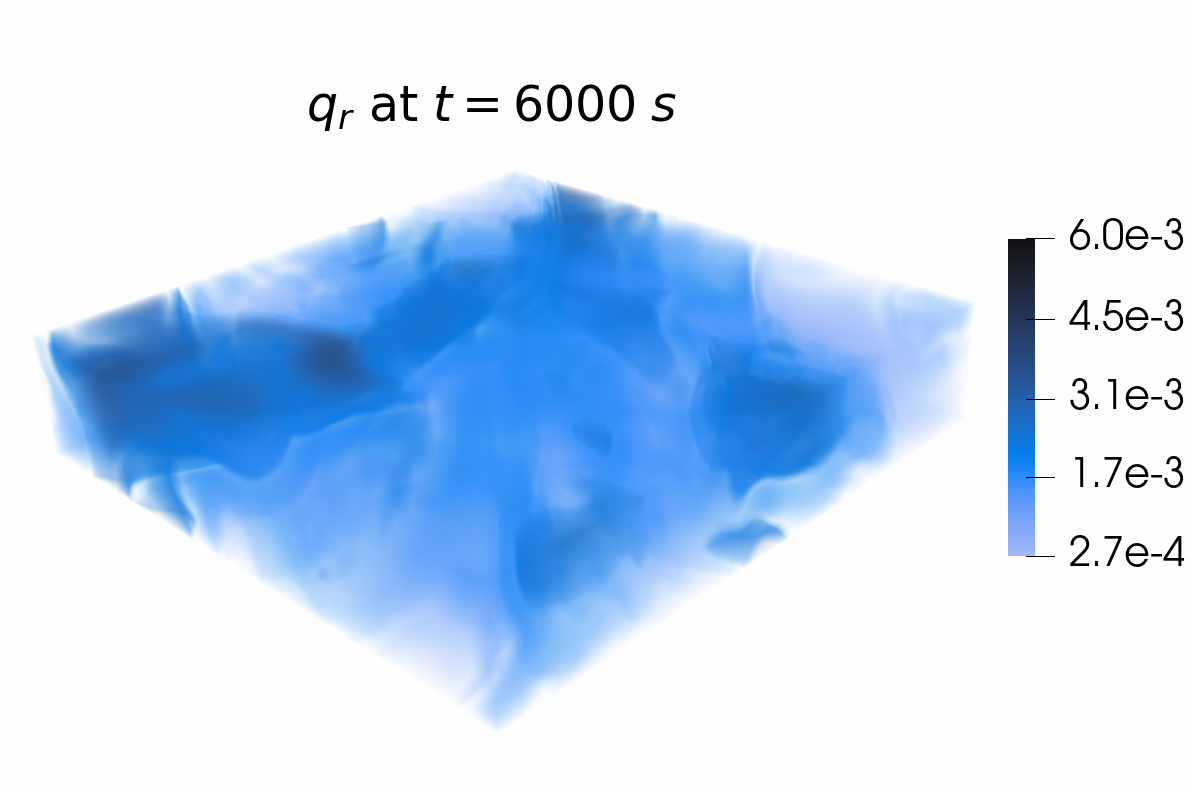}}
\vskip8pt
\centerline{\includegraphics[trim=1.2cm 2.2cm 0.1cm 2.5cm,clip,width=6.0cm]{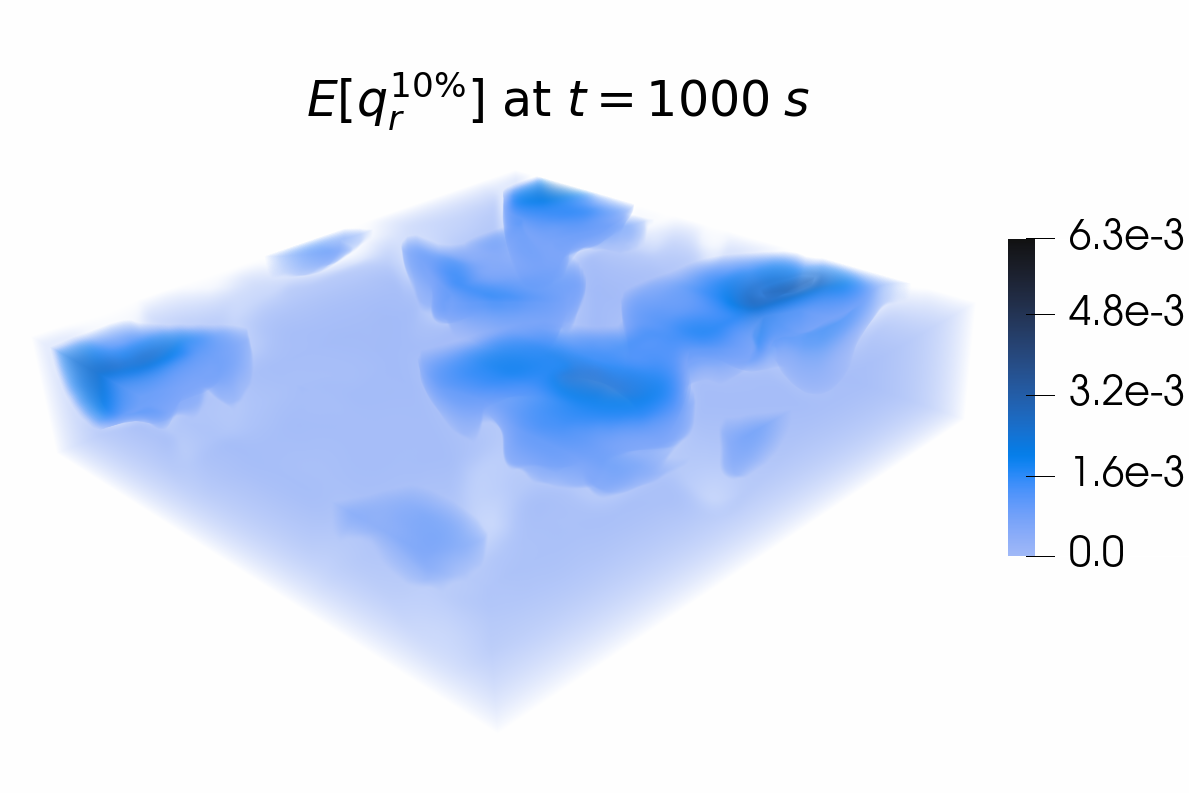}\hspace*{1.0cm}
\includegraphics[trim=1.2cm 2.2cm 0.1cm 2.5cm,clip,width=6.0cm]{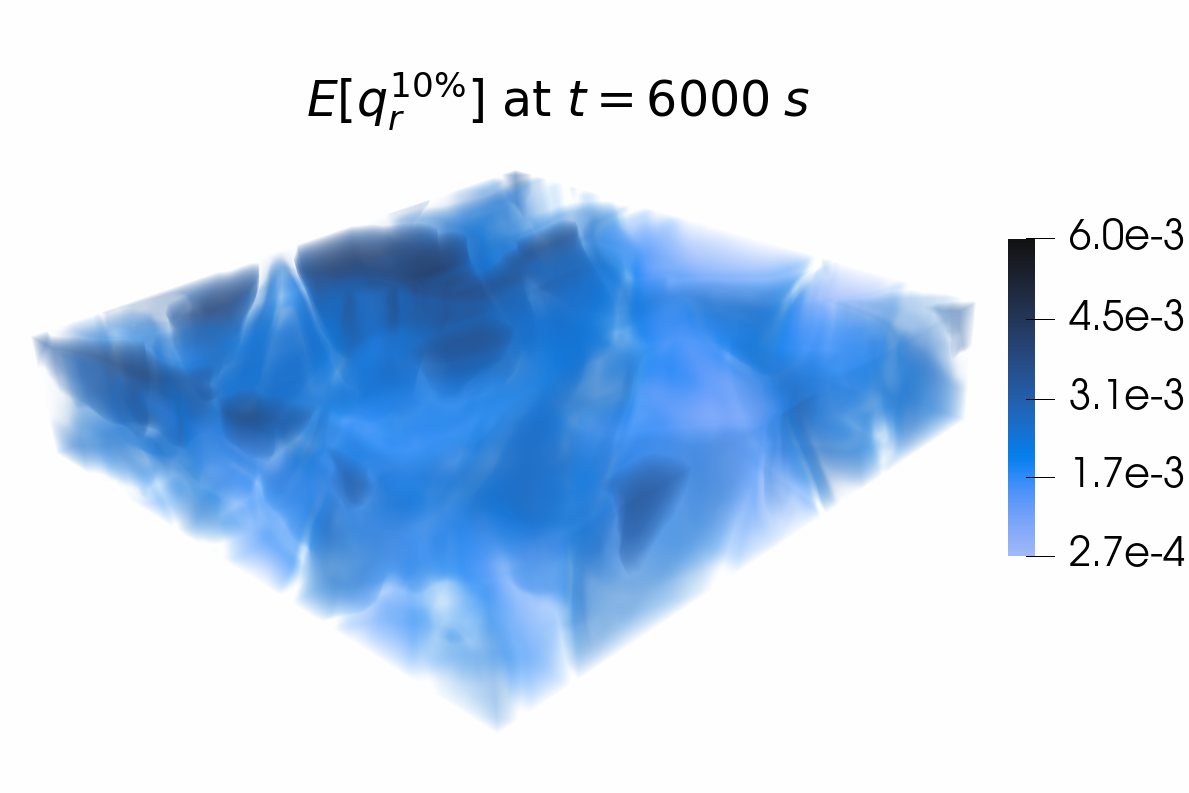}}
\vskip8pt
\centerline{\includegraphics[trim=1.2cm 2.2cm 0.1cm 2.5cm,clip,width=6.0cm]{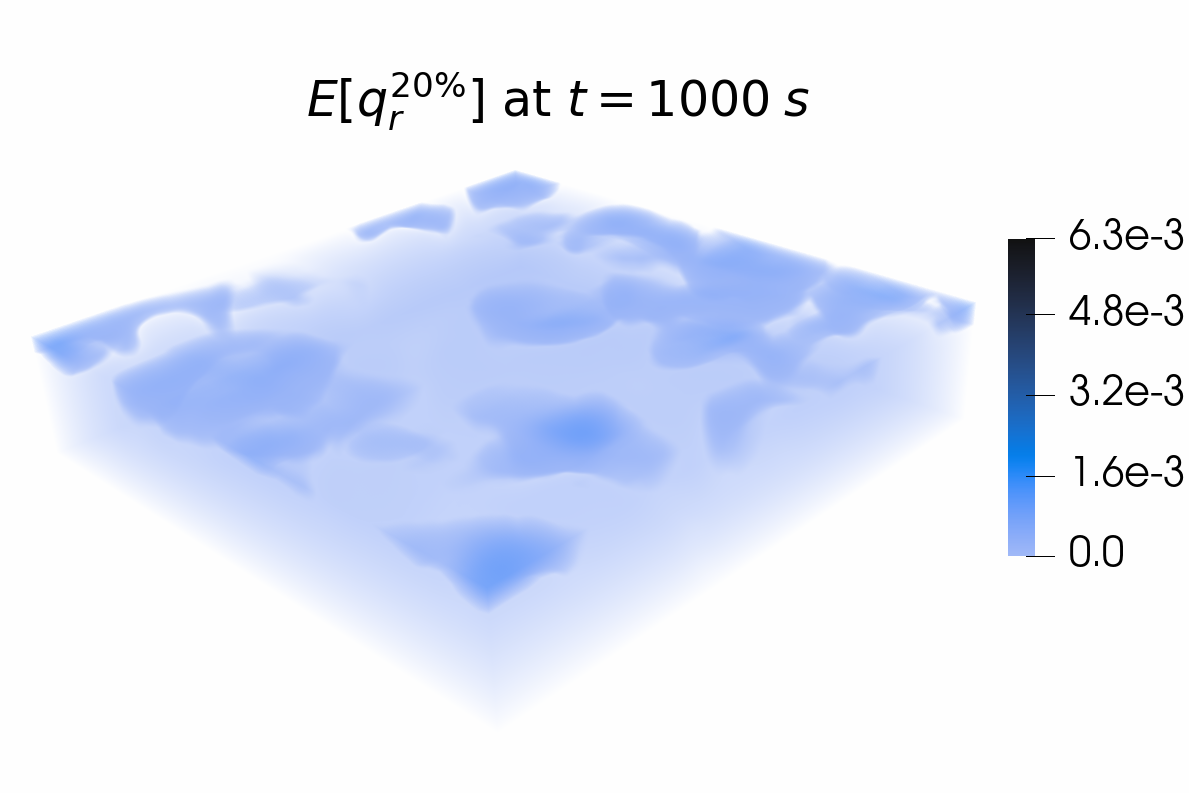}\hspace*{1.0cm}
\includegraphics[trim=1.2cm 2.2cm 0.1cm 2.5cm,clip,width=6.0cm]{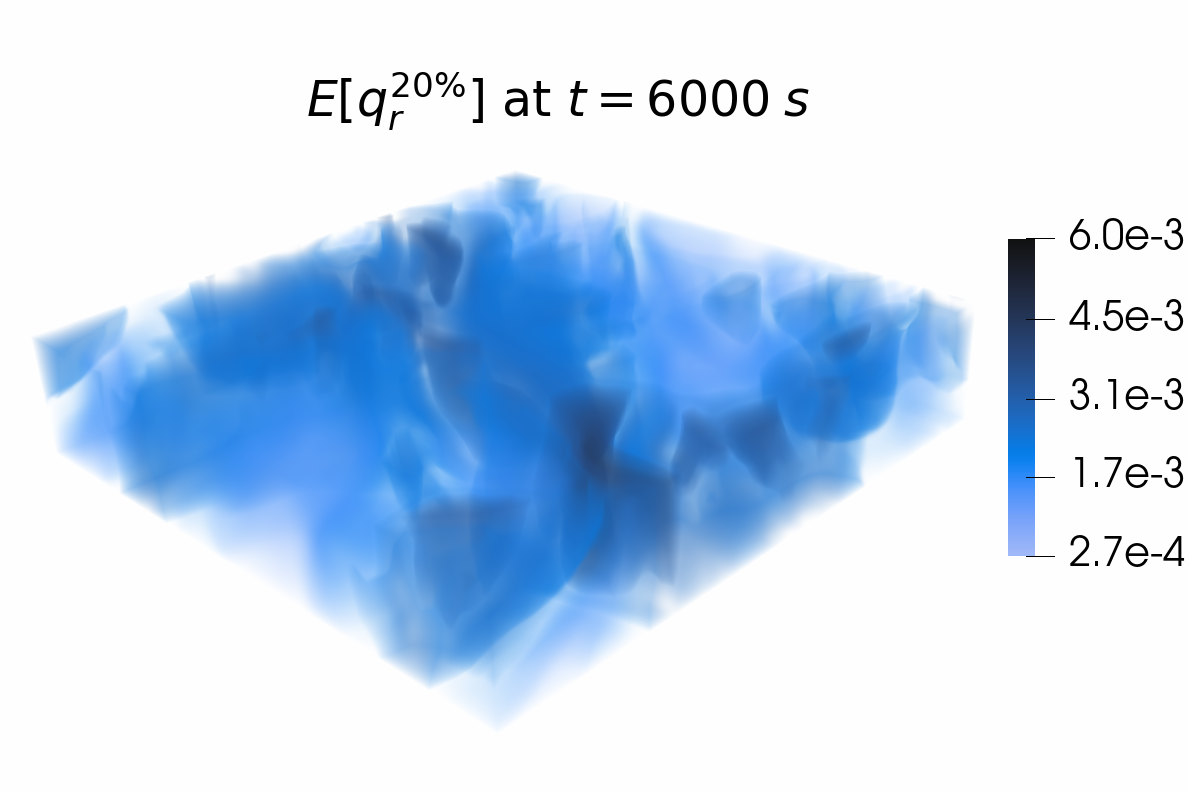}}
\vskip8pt
\centerline{\includegraphics[trim=1.2cm 2.2cm 0.1cm 2.5cm,clip,width=6.0cm]{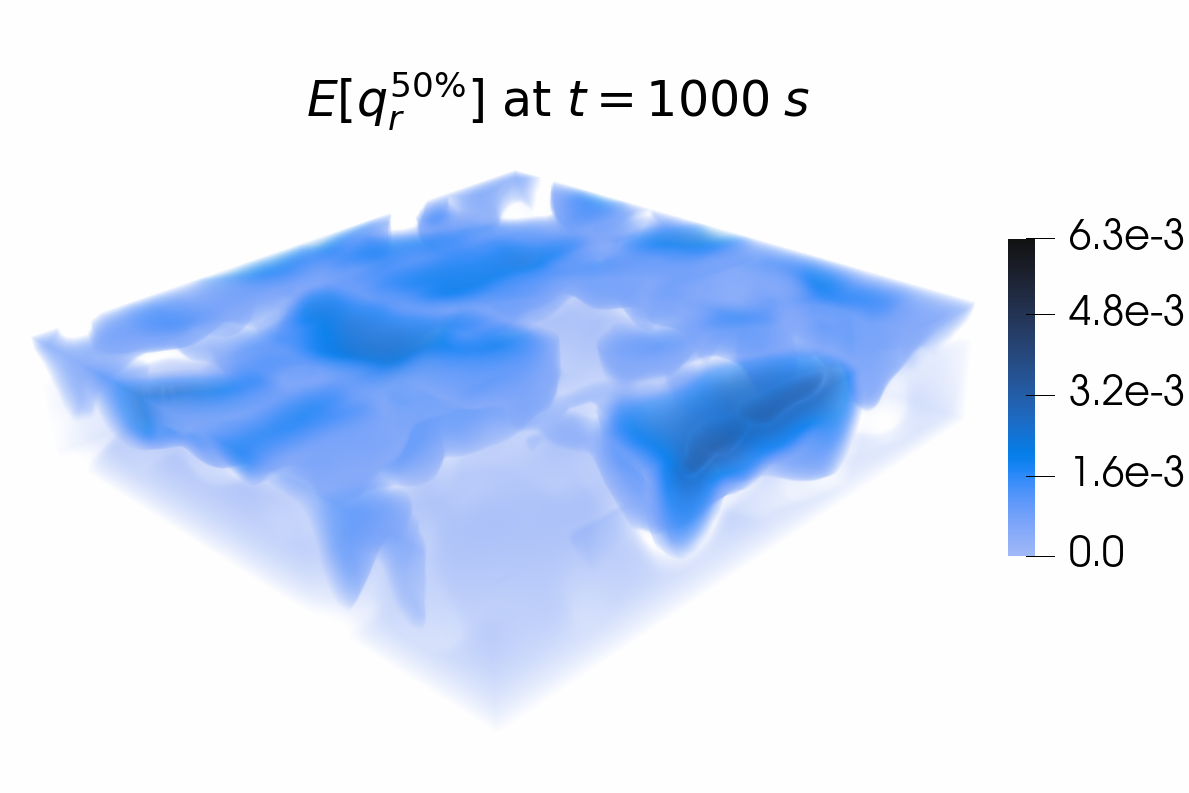}\hspace*{1.0cm}
\includegraphics[trim=1.2cm 2.2cm 0.1cm 2.5cm,clip,width=6.0cm]{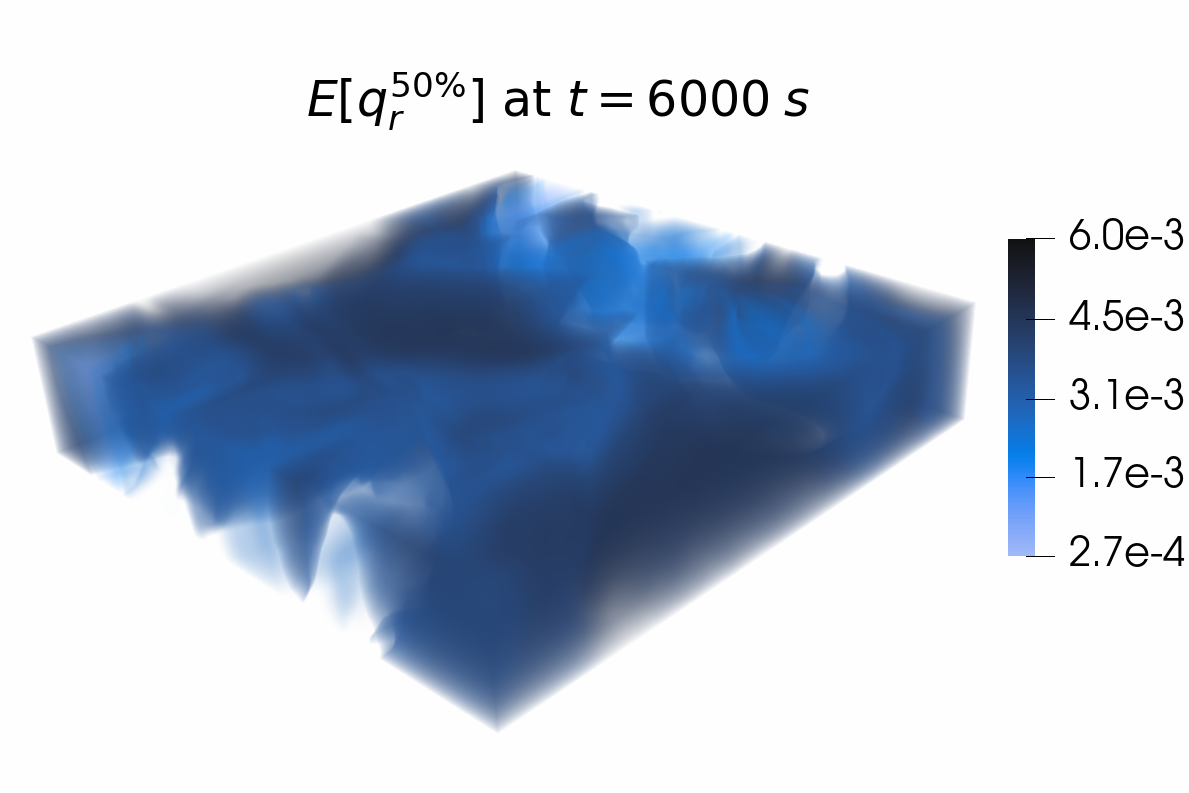}}
\caption{Example \ref{RB_stoch_example4}: Expected value of the rain concentration $q_r$ at times $t=1000s$ and $6000s$ with 0\%, 10\%,
20\% and 50\% perturbations of the initial water vapor concentration.\label{3D_RB_stoch_full_qr}}
\end{figure}
\begin{figure}[ht!]
\centerline{\includegraphics[trim=0.9cm 0.1cm 0.4cm 0.0cm,clip,width=11.0cm]{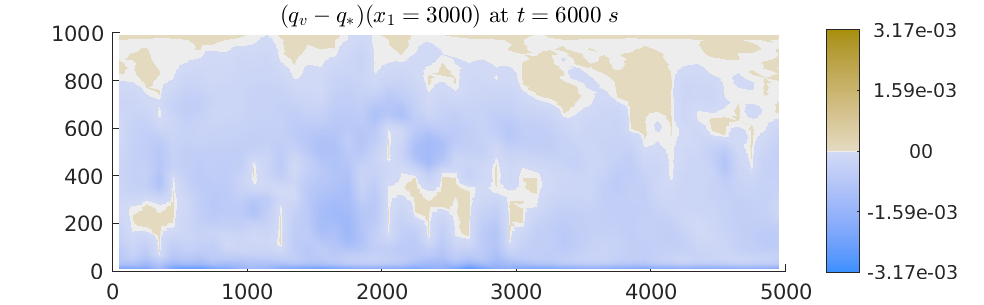}}
\vskip6pt
\centerline{\includegraphics[trim=0.9cm 0.1cm 0.4cm 0.0cm,clip,width=11.0cm]{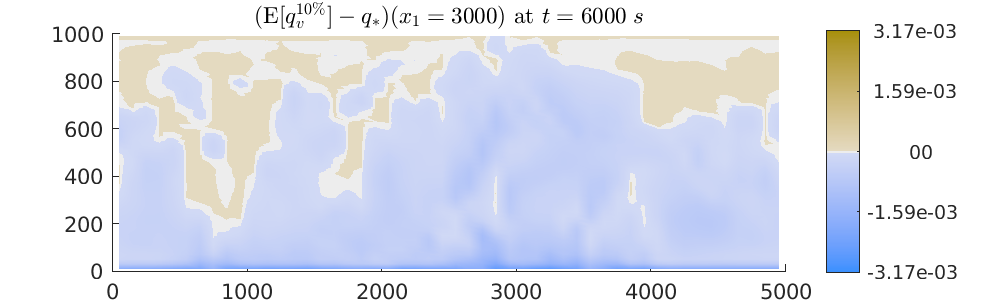}}
\vskip6pt
\centerline{\includegraphics[trim=0.9cm 0.1cm 0.4cm 0.0cm,clip,width=11.0cm]{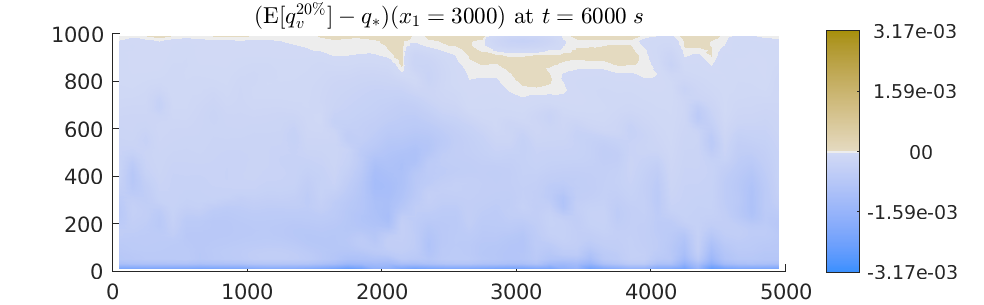}}
\vskip6pt
\centerline{\includegraphics[trim=0.9cm 0.1cm 0.4cm 0.0cm,clip,width=11.0cm]{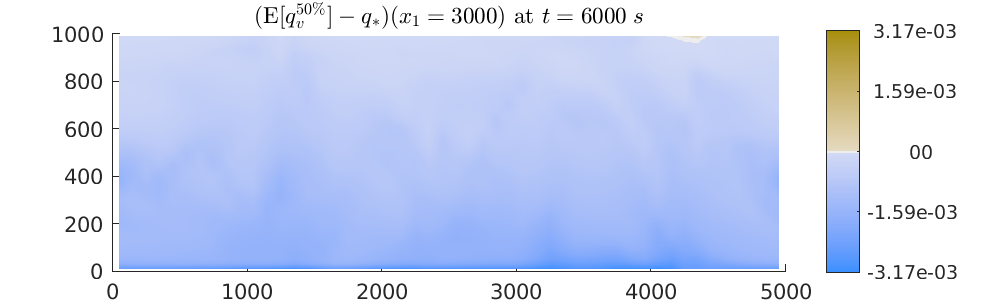}}
\caption{Example \ref{RB_stoch_example4}: Slices of the difference between the water vapor and the saturation mixing ratio
($\,\mathbb E[q_v]-q_*$) along $x_1=3000$ at time $t=6000s$ with 0\%, 10\%, 20\% and 50\% perturbations of the initial water vapor
concentration.\label{3D_RB_stoch_full_qstar_slice2}}
\end{figure}

In Figures \ref{3D_RB_stoch_full_E_and_sigma_theta} and \ref{3D_RB_stoch_full_E_and_sigma}, we show the time evolution of the mean expected
value per $m^3$ as well as the mean standard deviation per $m^3$ for the potential temperature and cloud variables in the cases with 0\%
(purely deterministic), 10\%, 20\% and 50\% perturbation of the initial water vapor concentration. For increasing perturbations, the spread
is increased, mostly for the water vapor concentration $q_v$ and the rain concentration $q_r$. The averaged quantities are dominated by the
positive perturbations, leading to (i) earlier cloud formation, (ii) thicker clouds due to more available water vapor, and (iii) enhanced
rain formation. These three features can be clearly seen in the case of the largest initial perturbation (50\%), where a large spread in
water vapor concentration is accompanied by a strong increase in cloud water and an earlier onset of strong precipitation. Due to the strong
rain formation the cloud concentration decreases when the perturbation size increases and also the amount of supersaturated regions is much
smaller as can be observed in Figure \ref{3D_RB_stoch_full_qstar_slice2} which leads to less new formation of clouds. We would also like to
note that the spread is only given by the standard deviation, whereas the actual minima (for instance, almost no cloud formation) cannot be
seen directly, although these scenarios are possible. Overall, one can see that the time evolution for the deterministic simulation as well
as for perturbations with 10\% and 20\% behave quite similarly and the averaged quantities follow closely the same evolution, although the
standard deviations increase quite substantially. However, for larger perturbations (50\%), the time evolution of the expected values of
$q_c$ and $q_r$ is strongly disturbed and shows large deviations from the other simulations. This can also be seen in the 3-D panels at the
later time $t=6000s$.
\begin{figure}[ht!]
\centerline{\includegraphics[trim=0.2cm 0.1cm 0.7cm 0.4cm,clip,width=5.0cm]{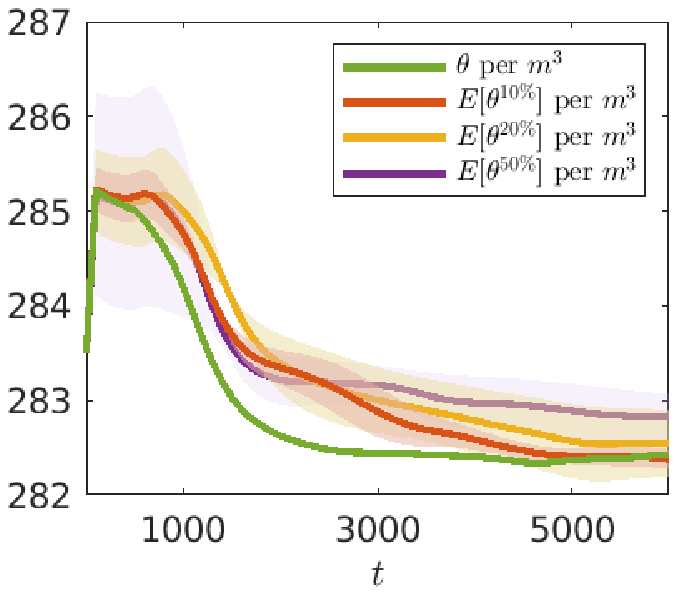}\hspace*{1.1cm}
\includegraphics[trim=0.0cm 0.0cm 0.0cm 0.0cm,clip,width=4.9cm]{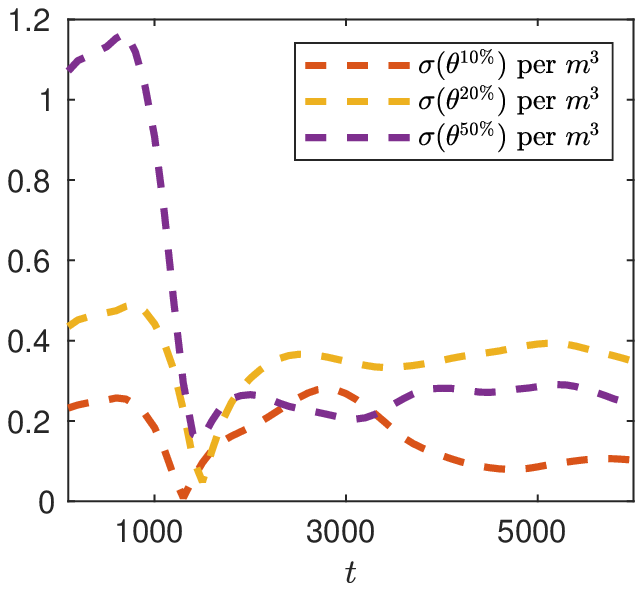}}
\caption{Example \ref{RB_stoch_example4}: Time evolution of the expected values with their standard deviations for the potential temperature
$\theta$ per $m^3$ (shaded region, left column) and standard deviations (right column) obtained using 0\%, 10\%, 20\% and 50\% perturbations
of the initial data in $q_v$.\label{3D_RB_stoch_full_E_and_sigma_theta}}
\end{figure}
\begin{figure}[ht!]
\hspace*{-0.1cm}\centerline{\includegraphics[trim=0.1cm 0.1cm 0.7cm 0.0cm,clip,width=5.1cm]{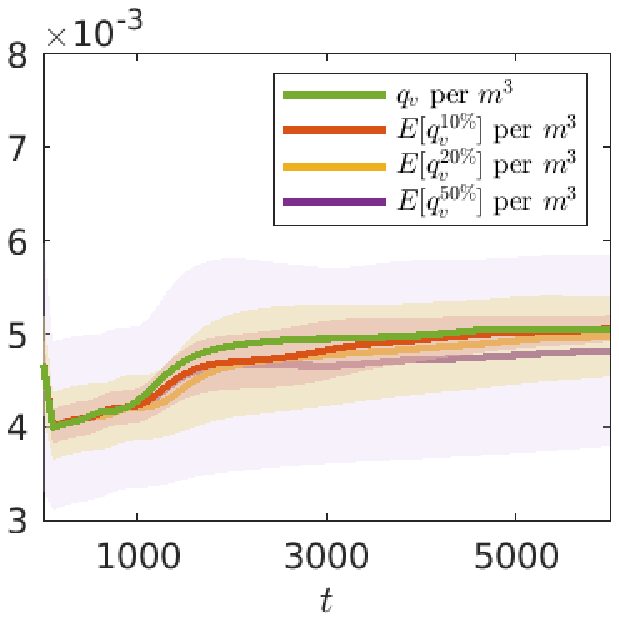}\hspace*{1.0cm}
\includegraphics[trim=0.0cm 0.0cm 0.0cm 0.0cm,clip,width=4.95cm]{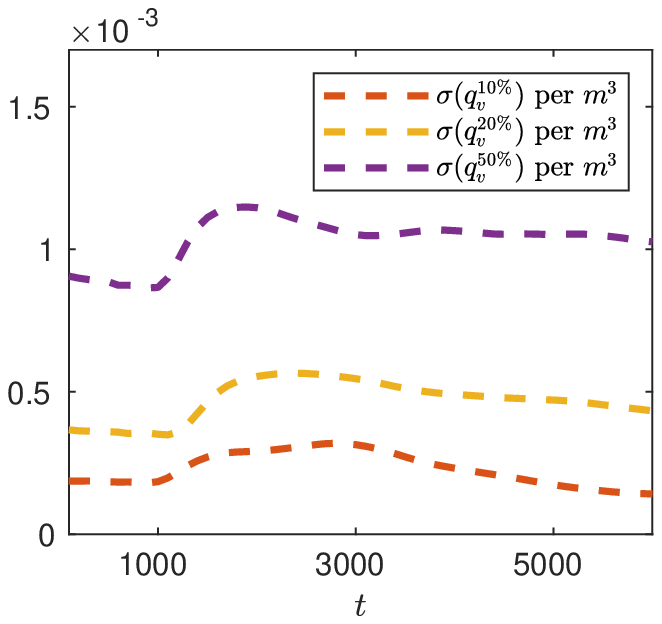}}
\vskip1pt
\centerline{\includegraphics[trim=0.3cm 0.1cm 0.7cm 0.0cm,clip,width=5.0cm]{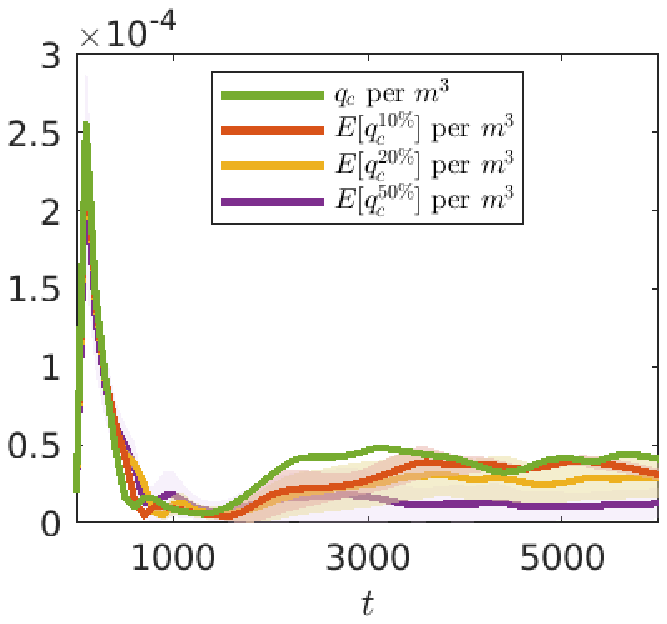}\hspace*{1.2cm}
\includegraphics[trim=0.0cm 0.0cm 0.0cm 0.0cm,clip,width=4.7cm]{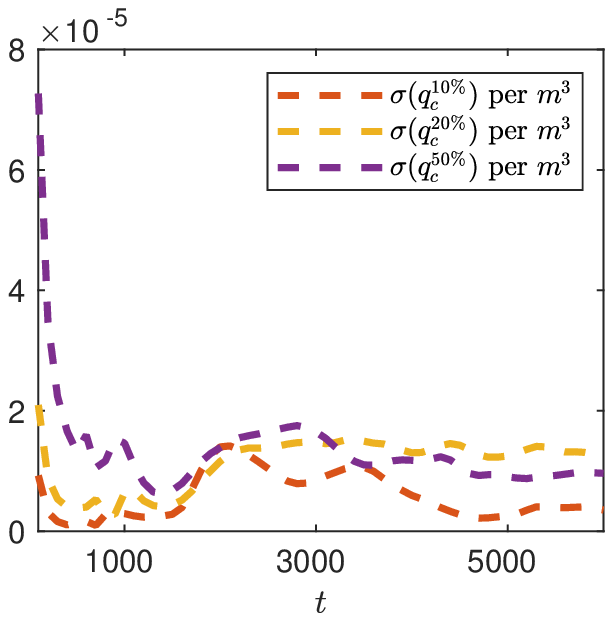}}
\vskip1pt
\centerline{\includegraphics[trim=0.3cm 0.1cm 0.7cm 0.0cm,clip,width=5.0cm]{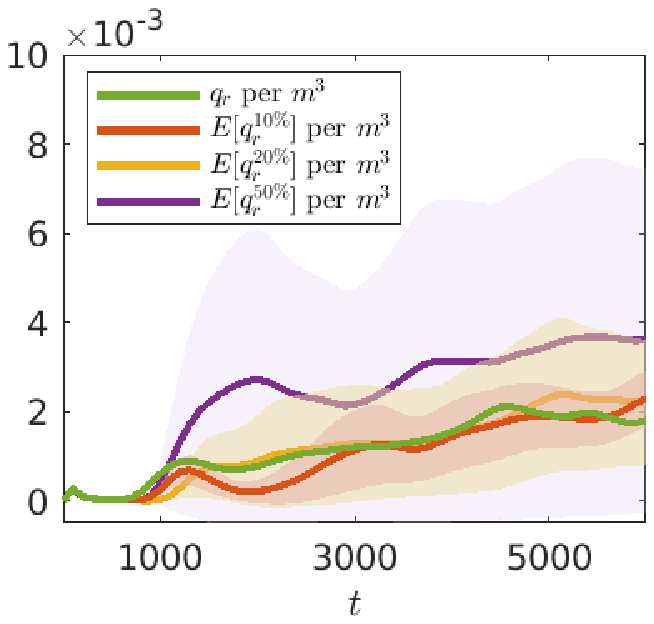}\hspace*{1.2cm}
\includegraphics[trim=0.0cm 0.0cm 0.0cm 0.0cm,clip,width=4.75cm]{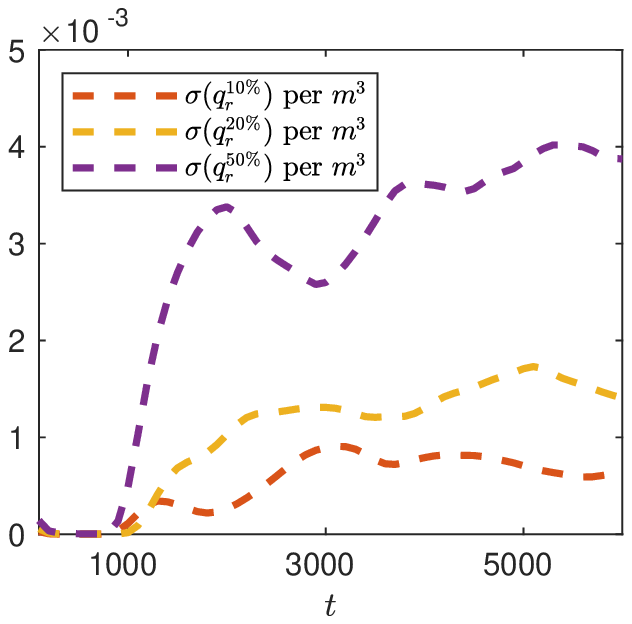}}
\caption{Example \ref{RB_stoch_example4}: Time evolution of the expected values with their standard deviations for the cloud variables per
$m^3$ (shaded region, left column) and standard deviations (right column) obtained using 0\%, 10\%, 20\% and 50\% perturbations of the
initial data in $q_v$.\label{3D_RB_stoch_full_E_and_sigma}}
\end{figure}
\end{example}

\section{Conclusion}
This paper is a continuation of our previous study on uncertainty propagation in atmospheric flows containing phase changes. In particular,
we consider warm cloud dynamics of weakly compressible fluids. This model consists of a multiscale system of PDEs in which the macroscopic
dynamics of the fluid is described by a weakly compressible Navier-Stokes system and the microscopic cloud dynamics is described by a system
of convection-diffusion-reaction equations. We have extended the gPC-SG method from \cite{SG2019}, where we considered a semi-random model
with the deterministic macroscopic dynamics of the fluid coupled with the random microscopic cloud dynamics, to the case of a fully random
multiscale system. To this end, we have first derived a system for the gPC coefficients and then presented a method we have used to
numerically solve the resulting system. The latter is an extension of the numerical method developed in \cite{SG2019} and approximates the
gPC coefficients for the dynamics of the fluid by an IMEX AP finite volume method and the gPC coefficients for the cloud dynamics by an
explicit finite volume method with an enlarged stability region.

The aim of this work is to demonstrate the applicability, accuracy and efficiency of the gPC-SG method for atmospheric flows. Comprehensive
studies of uncertainty propagation in these models considering different perturbation scenarios are left for a future work. Additionally, we
will investigate the accuracy and performance of different uncertainty quantification methods, for instance, stochastic Galerkin, stochastic
collocation and Monte Carlo method, in a review paper. Here, we have focused on numerical convergence and benchmark experiments as well as
comparison with the results of the previous semi-random model presented in \cite{SG2019}. We have demonstrated that the gPC-SG method for
the fully random model preserves the second-order spatial experimental convergence rate when the time increments are chosen according to the
time step restriction and additionally exhibited an experimental exponential convergence rate in the stochastic space. This experimental
convergence rate has been observed for both perturbation scenarios: uniform and normal distribution of the initial data perturbation.
Additionally, we have studied the numerical solutions of the fully random cloud model for both the 2-D and 3-D Rayleigh-B\'{e}nard
convection. By illustrating the behavior of clouds in different perturbed scenarios, we have demonstrated that perturbations of the initial
conditions of cloud variables can crucially change the time evolution. The results have also exhibited a clear difference of the solutions
of the semi- and fully random models in both the 2-D and 3-D Rayleigh-B\'{e}nard convection, which indicates that initial perturbations of
cloud variables propagate to the Navier-Stokes equations and have a significant effect on the fluid variables. Our numerical study clearly
demonstrates the applicability of the stochastic Galerkin method for the uncertainty quantification in complex atmospheric models and paves
the path for more extensive practically relevant numerical studies.

\begin{acknowledgment}
The work of A. Chertock was supported in part by NSF grant DMS-1818684. The work of A. Kurganov was supported in part by NSFC grants
12111530004 and 12171226, and by the fund of the Guangdong Provincial Key Laboratory of Computational Science and Material Design (No.
2019B030301001). The work of M.~Luk\'a\v{c}ov\'a-Medvi{\softd}ov\'a, P.~Spichtinger and B.~Wiebe was supported by the German Science Foundation
under the grant SFB/TRR 165 Waves to Weather (Project A2). M.~Luk\'a\v{c}ov\'a-Medvi{\softd}ov\'a and P.~Spichtinger gratefully ackonwledge the support of the Mainz Institute for Multiscale Modeling. M.~Luk\'a\v{c}ov\'a-Medvi{\softd}ov\'a thanks the Gutenberg Research College for supporting her research.
\end{acknowledgment}

\bibliographystyle{acm}
\bibliography{biblio}	
	
\end{document}